# SEMIPARAMETRICALLY EFFICIENT RANK-BASED INFERENCE FOR SHAPE
# I. OPTIMAL RANK-BASED TESTS FOR SPHERICITY

By Marc Hallin and Davy Paindaveine[1]

*Université Libre de Bruxelles*

We propose a class of rank-based procedures for testing that the shape matrix $\mathbf{V}$ of an elliptical distribution (with unspecified center of symmetry, scale and radial density) has some fixed value $\mathbf{V}_0$; this includes, for $\mathbf{V}_0 = \mathbf{I}_k$, the problem of testing for sphericity as an important particular case. The proposed tests are invariant under translations, monotone radial transformations, rotations and reflections with respect to the estimated center of symmetry. They are valid without any moment assumption. For adequately chosen scores, they are locally asymptotically maximin (in the Le Cam sense) at given radial densities. They are strictly distribution-free when the center of symmetry is specified, and asymptotically so when it must be estimated. The multivariate ranks used throughout are those of the distances—in the metric associated with the null value $\mathbf{V}_0$ of the shape matrix—between the observations and the (estimated) center of the distribution. Local powers (against elliptical alternatives) and asymptotic relative efficiencies (AREs) are derived with respect to the adjusted Mauchly test (a modified version of the Gaussian likelihood ratio procedure proposed by Muirhead and Waternaux [*Biometrika* **67** (1980) 31–43]) or, equivalently, with respect to (an extension of) the test for sphericity introduced by John [*Biometrika* **58** (1971) 169–174]. For Gaussian scores, these AREs are uniformly larger than one, irrespective of the actual radial density. Necessary and/or sufficient conditions for consistency under nonlocal, possibly nonelliptical alternatives are given. Finite sample performances are investigated via a Monte Carlo study.

Received February 2004; revised August 2005.

[1]Supported by a P.A.I. contract of the Belgian Federal Government and an Action de Recherche Concertée of the Communauté française de Belgique.

## 1. Introduction.

1.1. *The hypothesis of sphericity.* The distribution of a $k$-dimensional random vector $\mathbf{X}$ is called *spherical* if for some $\boldsymbol{\theta} \in \mathbb{R}^k$, the distribution of $\mathbf{X} - \boldsymbol{\theta}$ is invariant under orthogonal transformations. For multinormal variables, sphericity is equivalent to the covariance matrix $\boldsymbol{\Sigma}$ of $\mathbf{X}$ being proportional to the identity matrix $\mathbf{I}_k$. Under the assumption of ellipticity, finite second order moments need not exist and sphericity is equivalent to the shape matrix $\mathbf{V}$ being equal to the unit matrix $\mathbf{I}_k$; see Section 1.2 for precise definitions. Under more general nonelliptical distributions, however, this equivalence no longer holds: $\mathbf{V} = \mathbf{I}_k$ (the hypothesis of unit shape) does not imply sphericity, nor even that the directions $\mathbf{U} := (\mathbf{X} - \boldsymbol{\theta})/\|\mathbf{X} - \boldsymbol{\theta}\|$ be uniform over the unit sphere (the hypothesis of *isotropy*), and $\boldsymbol{\Sigma}$ (when finite) and $\mathbf{V}$ no longer coincide up to a positive scalar factor. The hypothesis of sphericity thus is a strict subhypothesis of the hypothesis of isotropy, itself a strict subhypothesis of the unit shape hypothesis. While isotropy and unit shape only deal with the $\mathbf{U}$'s, that is, with the directional features of $\mathbf{X} - \boldsymbol{\theta}$, sphericity also imposes a strong symmetry structure on the moduli $\|\mathbf{X} - \boldsymbol{\theta}\|$. This symmetry plays a crucial role in the approach we are adopting here and the null hypothesis of interest throughout is the hypothesis of spherical symmetry rather than that of unit shape; a detailed discussion of this issue is provided in Section 5.

Sphericity assumptions play a key role in a number of statistical problems, although the distinction between sphericity, isotropy, unit shape and a covariance matrix proportional to $\mathbf{I}_k$ is far from clear in the literature. Indeed, additional assumptions (Gaussian or elliptical densities, finite second-order moments, etc.)—the necessity of which is all too rarely questioned—in general cause the aforementioned assumptions to coincide. Besides this role as a technical assumption underlying some of the most frequently used statistical methods, null hypotheses of the type $\mathbf{V} = \mathbf{V}_0$—which (depending on the assumptions) reduce to the hypotheses of sphericity, isotropy or unit shape by substituting $\mathbf{V}_0^{-1/2}(\mathbf{X} - \boldsymbol{\theta})$ for $(\mathbf{X} - \boldsymbol{\theta})$—are also of direct interest in several specific domains of application such as geostatistics, paleomagnetic studies in geology, animal navigation, astronomy and wind direction data; see [5, 34, 53] or [35] for references. Finally, shape matrices provide robust alternatives to traditional covariance matrices; as such, they are obvious candidates for serving as the basic tools in a host of multivariate analysis techniques. Null hypotheses of the form $\mathbf{V} = \mathbf{V}_0$, in their various guises (reducing to sphericity, isotropy or unit shape) are thus of interest in their own right.

Because of its importance for applications, the problem of testing the hypothesis of sphericity has a long history and has generated a considerable body of literature which we only very briefly summarize here. For normal populations, the asymptotic theory has been thoroughly investigated. The Gaussian likelihood ratio test was introduced by Mauchly [36] and belongs to the classical toolkit of multivariate analysis. The (Gaussian) locally most powerful invariant (under shift, scale and orthogonal transformations) test was obtained by John [24, 25] and by Sugiura [48]. In their original versions, these tests are valid under Gaussian assumptions only; however, with slight modifications, they remain valid under elliptical populations with finite fourth-order moments; see Section 5.3 of [38] for the *adjusted* Mauchly test and Section 3.3 of the present paper for an adjusted version of John's test. Without elliptical symmetry, however, these adjusted tests are no longer valid; therefore, they qualify as tests of sphericity, not as tests of isotropy or unit shape. Moreover, it has been shown (see [23]) that they behave rather badly under heavy tails (a fact that is confirmed by the Monte Carlo study in Section 6). Although they still require elliptical symmetry and finite fourth-order radial moments, a more robust behavior can be expected from the test statistics introduced by Tyler [50], who proposes replacing covariance matrices with more robust estimators of scatter.

Non-Gaussian models have been investigated by Kariya and Eaton [26], where elliptical densities, possibly with infinite variances, are considered. Uniformly most powerful unbiased tests are derived, basically against specified nonspherical shape values. The results of that paper do not allow for more general optimality concepts (such as maximinity or stringency) involving unspecified shape alternatives. Despite their obvious theoretical value, such tests thus have limited practical value.

As a reaction to Gaussian and other strong distributional assumptions, nonparametric tests of sphericity have also been constructed. Their main advantage is that they remain consistent against all possible nonspherical alternatives, including the nonelliptical ones. The drawback is that they are computationally heavy and only achieve slow nonparametric consistency rates. Examples include Beran [6] and Koltchinskii and Sakhanenko [28] for the null hypothesis of ellipticity and Baringhaus [5] for sphericity. Another way of escaping Gaussian or fourth-order moment assumptions involves basing the tests on statistics that are measurable with respect to invariant or distribution-free quantities such as the multivariate concepts of signs and ranks developed, mainly, in the robustness literature; see [39] for a review.

This sign-and/or-rank-based approach has been adopted by Tyler [53], Ghosh and Sengupta [13] and Marden and Gao [33]. Tyler [53] addresses the problem of testing uniformity over the sphere for directional data and proposes a sign test related to his celebrated [52] estimator of shape. In a slightly different context, Ghosh and Sengupta [13] also propose a test entirely based on multivariate signs, that is, on cosines of the form $\mathbf{U}_i'\mathbf{U}_j =$ $(\mathbf{X}_i - \boldsymbol{\theta})'(\mathbf{X}_j - \boldsymbol{\theta})/\|\mathbf{X}_i - \boldsymbol{\theta}\|\|\mathbf{X}_j - \boldsymbol{\theta}\|$, where $\mathbf{X}_i$, $i = 1, \ldots, n$, denote the $k$-dimensional observations. These two multivariate sign tests are of a heuristic nature and do not rely on any clear optimality concerns. Their advantages and disadvantages are those which are usually associated with sign tests: they remain valid under a broad class of densities and are consistent against a broad class of alternatives, none of which requires elliptical symmetry. As a test of sphericity, the Ghosh and Sengupta test is not consistent against nonspherical alternatives with unit shape matrix. Therefore, rather than a test of sphericity, it should be considered a test of isotropy, or even a test of the hypothesis of unit shape with isotropic fourth-order directional moments; see Section 5 for details and an extension to the null hypothesis of unit shape. If ellipticity is assumed (so sphericity becomes the null hypothesis of interest), however, restricting to signs leads to a substantial loss of efficiency since the distances $d_i := \|\mathbf{X}_i - \boldsymbol{\theta}\|$, which are not taken into account, then also carry much relevant information.

In Marden and Gao [33], a variety of structural hypotheses on covariance matrices are considered, including sphericity and unit shape. Appropriate multivariate sign- and rank-based competitors of the Gaussian likelihood procedures are proposed. The ranks used by the authors are the *spatial ranks* introduced by Möttönen and Oja [37] and Chaudhuri [9]; see also [32]. Although Pitman efficiencies (with respect to the Gaussian methods) are obtained, no attempt is made to achieve any optimality and the authors restrict themselves to procedures of the Wilcoxon and sign test types; even so, they show that the sign-and-rank (Wilcoxon) procedures perform much better than those based on the signs alone, a finding that will be confirmed (both qualitatively and quantitatively) by the form of the information matrices we will derive in Section 2.

The approach we are adopting in the present paper is in the same spirit. However, throughout, we combine robustness (distribution-freeness under sphericity, without any moment assumptions) and optimality concerns. Our tests are based on multivariate signs and the ranks of the norms of the observations centered at $\boldsymbol{\theta}$ (or an estimate $\hat{\boldsymbol{\theta}}$), with test statistics that have a structure similar to that of John's. According to whether the center of symmetry $\boldsymbol{\theta}$ is specified or not, these statistics are strictly distribution-free under sphericity, or asymptotically so. With adequate scores, they are asymptotically optimal (in the Le Cam sense) against nonspherical elliptical distributions at chosen radial densities. In the elliptical setup, asymptotic relative efficiencies (AREs) with respect to the adjusted John and Mauchly tests are derived and appear to be surprisingly high (particularly for the van der Waerden version). Actually, Paindaveine [43] shows that the celebrated Chernoff and Savage [10] result concerning AREs of normal score tests for location with respect to Student's extends to the present situation:

the AREs of the normal score versions of our tests with respect to the traditional John–Mauchly–Muirhead–Waternaux tests are uniformly larger than one, irrespective of the underlying radial density.

The optimality properties of our tests are related to local elliptical alternatives; however, it is shown in Section 5.2 that provided nonconstant score functions are used (the "constant-score case" corresponds to the extended sign test proposed in Section 5.1), our tests nevertheless remain consistent against most elliptical as well as nonelliptical alternatives (some nonelliptical simulation results are reported in Section 6). We refer to Section 5 for a discussion of this matter.

Some basic reasons for considering sphericity as an alternative to classical Gaussian assumptions have been discussed above. Another non-Gaussian extension of the assumption of Gaussian sphericity is the assumption of i.i.d.-ness, under which the $k$ components of the observed $\mathbf{X}$ are independent and identically distributed (i.i.d.), with common unspecified symmetric marginal density $f$. Under Gaussian marginals, i.i.d.-ness and sphericity coincide, but not under general densities. In fact, Kac [27] shows that this hypothesis of i.i.d.-ness is rotation-invariant only on the class of multivariate Gaussian distributions. If Gaussian assumptions are abandoned, this hypothesis is no longer rotation-invariant and becomes strongly coordinate-dependent. Therefore, it does not fit into the semiparametric, coordinate-free setting we are adopting here.

However, the same assumption of i.i.d.-ness implies unit shape. The null hypothesis of i.i.d.-ness is thus strictly included in that of unit shape. Therefore, one might consider testing i.i.d.-ness by performing a test of unit shape such as the extended sign test proposed in Section 5.1. Although valid from the point of view of type I risk, such a test, for the null hypothesis of i.i.d.-ness, would be severely biased and inconsistent. Indeed the Maxwell–Hershell theorem (see, e.g., pages 51–52 of [8]) indicates that all non-Gaussian spherical distributions are part of the alternative, while our Proposition 5.1(v) establishes that $\alpha$-level extended sign tests at spherical alternatives have asymptotic power $\alpha$. For all of these reasons, it seems that i.i.d.-ness, in this context, is not the appropriate generalization of traditional Gaussian assumptions.

1.2. *Elliptical densities*: *location*, *scale*, *shape and radial density.* Denote by $\mathbf{X}^{(n)} := (\mathbf{X}_1^{(n)\prime}, \ldots, \mathbf{X}_n^{(n)\prime})'$, $n \in \mathbb{N}$, a triangular array of $k$-dimensional observations. Throughout, $\mathbf{X}_1^{(n)}, \ldots, \mathbf{X}_n^{(n)}$ are assumed to be i.i.d., with elliptical density

$$\underline{f}_{\boldsymbol{\theta},\sigma^2,\mathbf{V};f_1}(\mathbf{x}) := c_{k,f_1} \frac{1}{\sigma^k |\mathbf{V}|^{1/2}} f_1\left(\frac{1}{\sigma}((\mathbf{x}-\boldsymbol{\theta})'\mathbf{V}^{-1}(\mathbf{x}-\boldsymbol{\theta}))^{1/2}\right), \qquad \mathbf{x} \in \mathbb{R}^k, \tag{1.1}$$

where $\boldsymbol{\theta} \in \mathbb{R}^k$ is a *location parameter* $\sigma^2 \in \mathbb{R}_0^+ := (0, \infty)$ a *scale parameter*, and $\mathbf{V} := (V_{ij})$, a symmetric positive definite real $k \times k$ matrix with $V_{11} = 1$, a *shape parameter*. The infinite-dimensional parameter $f_1 : \mathbb{R}_0^+ \longrightarrow \mathbb{R}^+ := [0, \infty)$ is an a.e. strictly positive function, the constant $c_{k,f_1}$ a normalization factor depending on the dimension $k$ and $f_1$.

The function $f_1$ will be called, conveniently but improperly, a *radial density* ($f_1$ does not integrate to one, and is therefore not a probability density). Denote by $d_i^{(n)} = d_i^{(n)}(\boldsymbol{\theta}, \mathbf{V}) := \|\mathbf{Z}_i^{(n)}(\boldsymbol{\theta}, \mathbf{V})\|$ the modulus of the centered and *sphericized* observations $\mathbf{Z}_i^{(n)} = \mathbf{Z}_i^{(n)}(\boldsymbol{\theta}, \mathbf{V}) := \mathbf{V}^{-1/2}(\mathbf{X}_i^{(n)} - \boldsymbol{\theta})$, $i = 1, \dots, n$. If the $\mathbf{X}_i^{(n)}$'s have density (1.1), these moduli are i.i.d., with density and distribution functions

$$r \mapsto \frac{1}{\sigma} \tilde{f}_{1k}\left(\frac{r}{\sigma}\right) := \frac{1}{\sigma \mu_{k-1;f_1}} \left(\frac{r}{\sigma}\right)^{k-1} f_1\left(\frac{r}{\sigma}\right) I_{[r>0]}$$

and

$$r \mapsto \tilde{F}_{1k}(r/\sigma) := \int_0^{r/\sigma} \tilde{f}_{1k}(s)\, ds,$$

respectively, provided, however, that

$$\mu_{k-1;f_1} := \int_0^\infty r^{k-1} f_1(r)\, dr < \infty, \tag{1.2}$$

an assumption we shall henceforth always make on $f_1$. This function $\tilde{f}_{1k}$ is the *actual* radial density and (1.2) thus merely ensures that it will be a probability density function; in particular, it does not imply any moment restriction on $\tilde{f}_{1k}$, the $d_i^{(n)}$'s nor the $\mathbf{X}_i^{(n)}$'s. Any square root $\mathbf{V}^{1/2}$ of $\mathbf{V}$ [satisfying $\mathbf{V}^{1/2}(\mathbf{V}^{1/2})' = \mathbf{V}$] can be used in the above definitions, provided, of course, it is used in a consistent way. For the sake of simplicity, however, $\mathbf{A}^{1/2}$ throughout stands for the symmetric root of any symmetric positive semi-definite matrix $\mathbf{A}$, thus avoiding superfluous "prime" notation.

Now, if $\sigma$ and $f_1$ (or, more precisely, $\sigma$ and $c_{k,f_1} f_1$) are to be identifiable, a scale constraint is required. Still endeavoring to avoid moment restrictions, we impose the condition that the $d_i^{(n)}$'s, under (1.1), have common median $\sigma$, that is,

$$\tilde{F}_{1k}(1) = 1/2 \quad \text{or, equivalently,} \quad (\mu_{k-1;f_1})^{-1} \int_0^1 r^{k-1} f_1(r)\, dr = 1/2. \tag{1.3}$$

When this is to be emphasized, we call $f_1$ a *standardized* radial density. Special cases are:

(a) the $k$-variate multinormal distribution, with radial density $f_1(r) = \phi_1(r) := \exp(-a_k r^2/2)$;

(b) the $k$-variate Student distributions, with radial densities (for $\nu$ degrees of freedom) $f_1(r) = f^t_{1,\nu}(r) := (1 + a_{k,\nu} r^2/\nu)^{-(k+\nu)/2}$;

(c) the $k$-variate power-exponential distributions, with radial densities of the form $f_1(r) = f^e_{1,\eta}(r) := \exp(-b_{k,\eta} r^{2\eta})$.

[The constants $a_k > 0$, $a_{k,\nu} > 0$ and $b_{k,\eta} > 0$ are such that (1.3) is satisfied; note that $a_k = 2b_{k,1} = \lim_{\nu\to\infty} a_{k,\nu}$.]

Writing $\operatorname{vech}\mathbf{M} := (M_{11}, (\overset{\circ}{\operatorname{vech}}\,\mathbf{M})')'$ for the $k(k+1)/2$-dimensional vector obtained by stacking the upper-triangular elements of a $k \times k$ symmetric matrix $\mathbf{M} = (M_{ij})$, we denote by $\mathrm{P}^{(n)}_{\boldsymbol{\vartheta};f_1}$ or $\mathrm{P}^{(n)}_{\boldsymbol{\theta},\sigma^2,\mathbf{V};f_1}$ the distribution of $\mathbf{X}^{(n)}$ under given values of $\boldsymbol{\vartheta} = (\boldsymbol{\theta}', \sigma^2, (\overset{\circ}{\operatorname{vech}}\,\mathbf{V})')'$ and $f_1$ [$f_1$ satisfying (1.2) and (1.3)]. The parameter space is thus $\boldsymbol{\Theta} := \mathbb{R}^k \times \mathbb{R}^+_0 \times \mathcal{V}_k$, where $\mathcal{V}_k$ either stands for the set of all $k \times k$ symmetric positive definite matrices $\mathbf{V}$ such that $V_{11} = 1$ or for the corresponding set (in $\mathbb{R}^{(k(k+1)/2)-1}$) of values of $\overset{\circ}{\operatorname{vech}}\,\mathbf{V}$.

The notation $R^{(n)}_i = R^{(n)}_i(\boldsymbol{\theta}, \mathbf{V})$ will be used for the rank of $d^{(n)}_i = d^{(n)}_i(\boldsymbol{\theta}, \mathbf{V})$ among $d^{(n)}_1, \dots, d^{(n)}_n$; under $\mathrm{P}^{(n)}_{\boldsymbol{\vartheta};f_1}$, the vector $(R^{(n)}_1, \dots, R^{(n)}_n)$ is uniformly distributed over the $n!$ permutations of $(1, \dots, n)$. Let $\mathbf{U}^{(n)}_i = \mathbf{U}^{(n)}_i(\boldsymbol{\theta}, \mathbf{V}) := \mathbf{Z}^{(n)}_i / d^{(n)}_i$. The vectors $\mathbf{U}^{(n)}_i$ under $\mathrm{P}^{(n)}_{\boldsymbol{\vartheta};f_1}$ are i.i.d. and uniformly distributed over the unit sphere. They are independent of the ranks $R^{(n)}_i$ and are usually considered as *multivariate signs* associated with the centered observations $(\mathbf{X}_i - \boldsymbol{\theta})$ since they are totally insensitive to transformations of $(\mathbf{X}_i - \boldsymbol{\theta})$ that preserve half-lines through the origin.

The definition of the shape parameter $\mathbf{V}$ under elliptic symmetry readily follows from the special form of the density (1.1). A more general definition, which remains valid under possibly nonelliptical symmetric distributions, has been given by Tyler [52], where $\mathbf{V}$ is defined as the unique symmetric positive definite matrix such that $\mathbf{V}_{11} = 1$ (Tyler actually uses the normalization $\operatorname{tr}\mathbf{V} = k$) and

$$\mathrm{E}[(\mathbf{X} - \boldsymbol{\theta})(\mathbf{X} - \boldsymbol{\theta})' / (\mathbf{X} - \boldsymbol{\theta})' \mathbf{V}^{-1} (\mathbf{X} - \boldsymbol{\theta})] = \frac{1}{k}\mathbf{V}$$

(where $\boldsymbol{\theta}$ denotes the center of symmetry). The sample Tyler matrix $\mathbf{V}^{(n)}_T$ then provides a universally root-$n$ consistent estimator of $\mathbf{V}$. This ingenious extension may be somewhat misleading, however, as this "shape," in the absence of ellipticity, has a much weaker, and purely directional, interpretation. In particular, it is no longer associated with any family of contours and, under finite second-order moments, it loses its relation to covariance matrices—hence much of its intuitive content.

1.3. *Outline of the paper.* The problem we are considering is that of testing the hypothesis that the shape matrix $\mathbf{V}$ is equal to some given value $\mathbf{V}_0$ (admissible for a shape parameter). The special case $\mathbf{V}_0 = \mathbf{I}_k$, where $\mathbf{I}_k$ stands for the $k$-dimensional identity matrix, yields the problem of testing for sphericity. In the notation of the previous section, the shape matrix $\mathbf{V}$ in this problem is thus the parameter of interest, while $\boldsymbol{\theta}$, $\sigma^2$ and $f_1$ play the role of nuisance parameters. Hence, it is highly desirable that the null distributions of the test statistics to be used remain invariant under variations of $\boldsymbol{\theta}$, $\sigma^2$ and $f_1$.

When $\boldsymbol{\theta}$ is specified, we achieve this objective by basing our tests on the signs $\mathbf{U}_i^{(n)}$ and ranks $R_i^{(n)}$ computed from $\mathbf{Z}_i^{(n)}(\boldsymbol{\theta}, \mathbf{V}_0)$, $i = 1, \ldots, n$. These tests are invariant under monotone radial transformations (including scale transformations), rotations and reflections of the observations (with respect to $\boldsymbol{\theta}$), hence distribution-free with respect to such transformations. When $\boldsymbol{\theta}$ is unspecified, the ranks and signs are to be computed from $\mathbf{Z}_i^{(n)}(\hat{\boldsymbol{\theta}}, \mathbf{V}_0)$, $i = 1, \ldots, n$, where $\hat{\boldsymbol{\theta}} = \hat{\boldsymbol{\theta}}^{(n)}$ is an arbitrary root-$n$ consistent estimator of the location parameter $\boldsymbol{\theta}$; however, for $\hat{\boldsymbol{\theta}}$, we recommend the (rotation-equivariant) spatial median of Möttönen and Oja [37] which is itself "sign-based." This issue is treated in Section 4.4.

The tests $\underset{\sim}{\phi}{}_K^{(n)}$ based on these multivariate signed-rank statistics, whether ranks and signs are computed from $\boldsymbol{\theta}$ or from $\hat{\boldsymbol{\theta}}$, are locally asymptotically optimal (actually, *locally asymptotically maximin-efficient*, as the nonspecification of the scale $\sigma$ induces a strict loss of efficiency) in the Le Cam sense, under adequately chosen score functions. The test statistics take the very simple form (dropping superfluous superscripts, $c$ being some positive constant and $K$ a score function, see Section 4.2 for details)

$$\underset{\sim}{Q}{}_K = c\left(\operatorname{tr}\mathbf{S}_K^2 - \frac{1}{k}\operatorname{tr}^2\mathbf{S}_K\right) \qquad \text{with } \mathbf{S}_K := \frac{1}{n}\sum_{i=1}^n K\left(\frac{R_i}{n+1}\right)\mathbf{U}_i\mathbf{U}_i', \tag{1.4}$$

to be compared with the Gaussian statistic of John [24] ($d$ is some positive constant; see Section 3.3 for details),

$$Q_{\mathcal{N}} = \frac{d}{\operatorname{tr}^2\mathbf{S}}\left(\operatorname{tr}\mathbf{S}^2 - \frac{1}{k}\operatorname{tr}^2\mathbf{S}\right) \qquad \text{with } \mathbf{S} := \frac{1}{n}\sum_{i=1}^n \mathbf{Z}_i\mathbf{Z}_i' = \frac{1}{n}\sum_{i=1}^n d_i^2\mathbf{U}_i\mathbf{U}_i'. \tag{1.5}$$

The special case of a constant score $[K(u) = 1,\ 0 < u < 1]$ yields $\mathbf{S}_S := \frac{1}{n}\sum_{i=1}^n \mathbf{U}_i\mathbf{U}_i'$ and a test $\underset{\sim}{\phi}{}_S^{(n)}$ which is essentially that proposed by Ghosh and Sengupta [13].

The rest of the paper is organized as follows. In Section 2 we establish the local asymptotic normality result (with respect to the location, scale and shape parameters) that provides the main theoretical tool of the paper. This result allows the development of asymptotically optimal parametric procedures for $\mathbf{V}$ under specified values of $f_1$ and $\sigma$ (with possibly unspecified $\boldsymbol{\theta}$). This is explained in detail in Section 3 where we also derive the asymptotically optimal (efficient, at given $f_1$) "$\sigma$-free" tests for hypotheses of the form $\mathbf{V} = \mathbf{V}_0$ (tests for sphericity being a special case) and explicitly compute the loss (in local powers) due to the nonspecification of scale. The Gaussian version of this test is investigated further and its link with some classical tests of sphericity is discussed. In Section 4 we propose nonparametric (signed-rank-based) versions of the optimal procedures defined in Section 3 and study their invariance and asymptotic properties. Asymptotic relative efficiencies with respect to the parametric Gaussian tests are derived in the elliptical case. All of these results are obtained under specified $\boldsymbol{\theta}$ first; then, in Section 4.4 we show that under minimal regularity assumptions on the actual underlying density (essentially, those ensuring ULAN), $\boldsymbol{\theta}$ can safely be replaced by any root-$n$ consistent estimator $\hat{\boldsymbol{\theta}}^{(n)}$. In Section 5, we study the validity (under extensions of the null hypothesis of sphericity) and consistency properties under nonlocal alternatives of our testing procedures. An adjusted version of the sign test is proposed, extending the validity of $\underset{\sim}{\phi}{}_S^{(n)}$ to the null hypothesis of unit shape. Necessary and/or sufficient conditions for consistency are established. For Wilcoxon (i.e., linear) scores, these necessary and sufficient conditions take a very simple form which shows that the corresponding rank tests are consistent against essentially all nonspherical alternatives, including the nonelliptical ones. As for the (adjusted) sign test $\underset{\sim}{\phi}{}_S^{(n)}$, it is shown to be consistent against all non-unit-shape alternatives, confirming its qualification as a fully consistent test for unit shape. Section 6 provides some simulation results which indicate that finite-sample performances reflect the asymptotic powers derived in the previous sections, as well as the nonelliptical consistency property established in Section 5. The Appendix compiles some technical proofs.

1.4. *Notation.* The following notation will be used throughout. Denoting by $\mathbf{e}_\ell$ the $\ell$th vector in the canonical basis of $\mathbb{R}^k$ and by $\mathbf{I}_k$ the $k \times k$ unit matrix, let

$$\mathbf{K}_k := \sum_{i,j=1}^{k} (\mathbf{e}_i \mathbf{e}_j') \otimes (\mathbf{e}_j \mathbf{e}_i') \quad \text{and} \quad \mathbf{J}_k := \sum_{i,j=1}^{k} (\mathbf{e}_i \mathbf{e}_j') \otimes (\mathbf{e}_i \mathbf{e}_j') = (\operatorname{vec} \mathbf{I}_k)(\operatorname{vec} \mathbf{I}_k)';$$

the $k^2 \times k^2$ matrix $\mathbf{K}_k$ is known as the *commutation matrix*. With this notation, $\mathbf{K}_k \operatorname{vec}(\mathbf{A}) = \operatorname{vec}(\mathbf{A}')$ and $\mathbf{J}_k \operatorname{vec}(\mathbf{A}) = (\operatorname{tr} \mathbf{A})(\operatorname{vec} \mathbf{I}_k)$. Note that $(1/k)\mathbf{J}_k$ and $\mathbf{J}_k^{\perp} := \mathbf{I}_{k^2} - (1/k)\mathbf{J}_k$ are the matrices of the projections onto the mutually orthogonal subspaces $\{\lambda(\operatorname{vec}\mathbf{I}_k)|\lambda \in \mathbb{R}\}$ and $\{\operatorname{vec}(\mathbf{A})|\operatorname{tr}\mathbf{A} = 0\}$, respectively. Define $\mathbf{M}_k$ as the $(k(k+1)/2-1) \times k^2$ matrix such that $\mathbf{M}_k'(\overset{\circ}{\operatorname{vech}}(\mathbf{v})) = \operatorname{vec}(\mathbf{v})$ for any symmetric $k \times k$ matrix $\mathbf{v} = (v_{ij})$ with $v_{11} = 0$, and let $\mathbf{N}_k$ be the $(k(k+1)/2-1) \times k^2$ real matrix such that $\mathbf{N}_k(\operatorname{vec}\mathbf{v}) = \overset{\circ}{\operatorname{vech}}\,\mathbf{v}$ for any symmetric $k \times k$ matrix $\mathbf{v}$. Finally, we write $\mathbf{V}^{\otimes 2}$ for the Kronecker product $\mathbf{V} \otimes \mathbf{V}$.

## 2. Uniform local asymptotic normality (ULAN).

Our objective is to perform inference on the shape parameter $\mathbf{V}$ under unspecified location $\boldsymbol{\theta}$, unspecified scale $\sigma$ and unspecified standardized radial density $f_1$: $\mathbf{V}$ is thus the parameter of interest, whereas $\boldsymbol{\theta}$, $\sigma^2$ and $f_1$ play the role of nuisance parameters. The relevant statistical experiment involves the nonparametric family

$$
\begin{aligned}
\mathcal{P}^{(n)} &:= \bigcup_{f_1 \in \mathcal{F}_A} \mathcal{P}_{f_1}^{(n)} := \bigcup_{f_1 \in \mathcal{F}_A} \bigcup_{\sigma > 0} \mathcal{P}_{\sigma^2; f_1}^{(n)} \\
&:= \bigcup_{f_1 \in \mathcal{F}_A} \bigcup_{\sigma > 0} \{\mathrm{P}_{\boldsymbol{\theta}, \sigma^2, \mathbf{V}; f_1}^{(n)} | \boldsymbol{\theta} \in \mathbb{R}^k, \mathbf{V} \in \mathcal{V}_k\}
\end{aligned}
\tag{2.1}
$$

[$f_1$ ranges over the set $\mathcal{F}_A$ of standardized densities satisfying Assumptions (A1) and (A2) below], in which the partition of $\mathcal{P}^{(n)}$ into a collection of parametric subexperiments $\mathcal{P}_{\sigma^2; f_1}^{(n)}$ all indexed by the same parameters $\boldsymbol{\theta}$ and $\mathbf{V}$, induces a semiparametric structure. The main technical tool is the uniform local asymptotic normality (ULAN), with respect to $\boldsymbol{\vartheta} = (\boldsymbol{\theta}', \sigma^2, (\overset{\circ}{\operatorname{vech}}\,\mathbf{V})')'$, of the families $\mathcal{P}_{f_1}^{(n)}$. This LAN (ULAN) issue has been briefly touched by Bickel (Example 4 in [7]). The very particular case of bivariate distributions with finite second-order moments has been treated recently by Falk [11] in his investigation of the inefficiency of empirical correlation coefficients.

In order to describe the extremely mild assumptions under which the families $\mathcal{P}_{f_1}^{(n)}$ are ULAN, we introduce the following definitions. Consider the measure space $(\Omega, \mathbb{B}_\Omega, \lambda)$, where $\lambda$ is some measure on the open subset $\Omega \subset \mathbb{R}$ equipped with its Borel $\sigma$-field $\mathbb{B}_\Omega$. Denote by $L^2(\Omega, \lambda)$ the space of measurable functions $h : \Omega \to \mathbb{R}$ satisfying $\int_\Omega [h(x)]^2 \, d\lambda(x) < \infty$. In particular, consider the space $L^2(\mathbb{R}_0^+, \mu_\ell)$ [resp. $L^2(\mathbb{R}, \nu_\ell)$] of square integrable functions w.r.t. the Lebesgue measure with weight $x^\ell$ on $\mathbb{R}_0^+$ (resp. with weight $e^{\ell x}$ on $\mathbb{R}$), that is, the space of measurable functions $h : \mathbb{R}_0^+ \to \mathbb{R}$ satisfying $\int_0^\infty [h(x)]^2 x^\ell \, dx < \infty$ (resp. $h : \mathbb{R} \to \mathbb{R}$ satisfying $\int_{-\infty}^\infty [h(x)]^2 e^{\ell x} \, dx < \infty$). Recall that $g \in L^2(\Omega, \lambda)$ admits a *weak derivative* $T$ iff $\int_\Omega g(x)\varphi'(x)\, dx = -\int_\Omega T(x)\varphi(x)\, dx$ for all infinitely differentiable (in the classical sense) compactly supported functions $\varphi$ on $\Omega$. The mapping $T$ is also called the *derivative of $g$ in the sense of distributions* in $L^2(\Omega, \lambda)$. If, moreover, $T$ itself is in $L^2(\Omega,\lambda)$, then $g$ belongs to $W^{1,2}(\Omega,\lambda)$, the Sobolev space of order 1 on $L^2(\Omega,\lambda)$. For the sake of simplicity, we will write $L^2(\Omega)$ and $W^{1,2}(\Omega)$ when $\lambda$ is the Lebesgue measure on $\Omega$. The family $\mathcal{P}_{f_1}^{(n)}$ is ULAN under the following assumptions on the radial density $f_1$:

ASSUMPTION (A1). The mapping $x \mapsto f_1^{1/2}(x)$ is in $W^{1,2}(\mathbb{R}_0^+, \mu_{k-1})$.

Letting $\varphi_{f_1}(r) := -2(f_1^{1/2})'(r)/f_1^{1/2}(r)$, where $(f_1^{1/2})'$ stands for the weak derivative of $f_1^{1/2}$ in $L^2(\mathbb{R}_0^+, \mu_{k-1})$, Assumption (A1) ensures the finiteness of *radial Fisher information for location* (expectation is taken under $\mathrm{P}_{\boldsymbol{\vartheta};f_1}^{(n)}$),

$$\mathcal{I}_k(f_1) := \mathrm{E}[\varphi_{f_1}^2(d_i^{(n)}(\boldsymbol{\theta},\mathbf{V})/\sigma)] = \int_0^1 \varphi_{f_1}^2(\tilde{F}_{1k}^{-1}(u))\,du.$$

ASSUMPTION (A2). The mapping $x \mapsto f_{1;\exp}^{1/2}(x) := f_1^{1/2}(e^x)$ is in $W^{1,2}(\mathbb{R}, \nu_k)$.

Letting $\psi_{f_1}(r) := -2r^{-1}(f_{1;\exp}^{1/2})'(\ln r)/f_{1;\exp}^{1/2}(\ln r)$, where $(f_{1;\exp}^{1/2})'$ stands for the weak derivative of $f_{1;\exp}^{1/2}$ in $L^2(\mathbb{R},\nu_k)$ and $K_{f_1}(u) := \psi_{f_1}(\tilde{F}_{1k}^{-1}(u))\tilde{F}_{1k}^{-1}(u)$, Assumption (A2) ensures the finiteness of *radial Fisher information for shape* (and *scale*—expectations are still taken under $\mathrm{P}_{\boldsymbol{\vartheta};f_1}^{(n)}$),

$$\mathcal{J}_k(f_1) := \mathrm{E}[\psi_{f_1}^2(d_i^{(n)}(\boldsymbol{\theta},\mathbf{V})/\sigma)(d_i^{(n)}(\boldsymbol{\theta},\mathbf{V})/\sigma)^2] = \int_0^1 K_{f_1}^2(u)\,du.$$

In principle, the functions $\varphi_{f_1}$ and $\psi_{f_1}$ differ. However, they do coincide (a.e.) under the following Assumption (A1-2), which, though slightly more stringent than Assumptions (A1) and (A2), holds for most densities considered in practice:

ASSUMPTION (A1-2). The radial density $f_1$ is absolutely continuous with a.e.-derivative $\dot{f}_1$ and, letting $\varphi_{f_1} = \psi_{f_1} := -\dot{f}_1/f_1$, the integrals

$$\mathcal{I}_k(f_1) := \int_0^1 \varphi_{f_1}^2(\tilde{F}_{1k}^{-1}(u))\,du \quad \text{and} \quad \mathcal{J}_k(f_1) := \int_0^1 K_{f_1}^2(u)\,du$$

are finite.

It should be stressed that none of these assumptions requires the existence of any moment for the radial density $\tilde{f}_{1k}$. They are satisfied, for instance, for all multivariate Student radial densities, including the Cauchy ones. For the radial density $f^t_{1,\nu}$ of the $k$-variate $t$-distribution with $\nu$ degrees of freedom [$\nu \in (0,\infty)$], it can be checked that

$$\mathcal{I}_k(f^t_{1,\nu}) = a_{k,\nu}\frac{k(k+\nu)}{k+\nu+2} \quad \text{and} \quad \mathcal{J}_k(f^t_{1,\nu}) = \frac{k(k+2)(k+\nu)}{k+\nu+2}. \tag{2.2}$$

The same remark holds for the power-exponential distributions, provided that $k \geq 2$ (which is not a limitation, since the problem under consideration is void for $k=1$), with

$$\mathcal{I}_k(f^e_{1,\eta}) = 4\eta^2 b_{k,\eta}\frac{\Gamma((4\eta+k-2)/2\eta)}{\Gamma(k/2\eta)} \quad \text{and} \quad \mathcal{J}_k(f^e_{1,\eta}) = k(k+2\eta) \tag{2.3}$$

($\Gamma$ denotes the Euler Gamma function). The corresponding values for $k$-variate multinormal distributions can be obtained by taking limits of the information quantities in (2.2) as $\nu \to \infty$ or, equivalently, by evaluating (2.3) at $\eta = 1$:

$$\mathcal{I}_k(\phi_1) = a_k k \quad \text{and} \quad \mathcal{J}_k(\phi_1) = k(k+2).$$

Note that $\lim_{\nu\to 0}\mathcal{J}_k(f^t_{1,\nu}) = \lim_{\eta\to 0}\mathcal{J}_k(f^e_{1,\eta}) = k^2$, which is a sharp lower bound for radial shape/scale information since, by Jensen's inequality and integration by parts,

$$(\mathcal{J}_k(f_1))^{1/2} \geq \int_0^1 K_{f_1}(u)\,du = \int_0^1 \psi_{f_1}(\tilde{F}^{-1}_{1k}(u))\tilde{F}^{-1}_{1k}(u)\,du = k. \tag{2.4}$$

Similarly, assuming that the density in (1.1) has finite second-order moments, the radial information for location $\mathcal{I}_k(f_1)$ satisfies (the Cauchy–Schwarz inequality)

$$\mathcal{I}_k(f_1) \geq k^2\left(\int_0^1 (\tilde{F}^{-1}_{1k}(u))^2\,du\right)^{-1},$$

with equality in the multinormal case only.

PROPOSITION 2.1. *Under Assumptions (A1) and (A2), the family* $\mathcal{P}^{(n)}_{f_1} := \{\mathrm{P}^{(n)}_{\boldsymbol{\vartheta};f_1} | \boldsymbol{\vartheta} \in \boldsymbol{\Theta}\}$ *is ULAN, with [writing* $d_i$ *and* $\mathbf{U}_i$*, resp., for* $d^{(n)}_i(\boldsymbol{\theta},\mathbf{V})$ *and* $\mathbf{U}^{(n)}_i(\boldsymbol{\theta},\mathbf{V})$*] central sequence*

$$\begin{aligned}\boldsymbol{\Delta}^{(n)}_{f_1}(\boldsymbol{\vartheta}) &:= \begin{pmatrix}\boldsymbol{\Delta}^{(n)}_{f_1;1}(\boldsymbol{\vartheta})\\ \Delta^{(n)}_{f_1;2}(\boldsymbol{\vartheta})\\ \boldsymbol{\Delta}^{(n)}_{f_1;3}(\boldsymbol{\vartheta})\end{pmatrix}\\ &:= \begin{pmatrix} n^{-1/2}\dfrac{1}{\sigma}\displaystyle\sum_{i=1}^n \varphi_{f_1}\left(\frac{d_i}{\sigma}\right)\mathbf{V}^{-1/2}\mathbf{U}_i \\ \dfrac{1}{2}n^{-1/2}\begin{pmatrix}\sigma^{-2}(\operatorname{vec}\mathbf{I}_k)'\\ \mathbf{M}_k(\mathbf{V}^{\otimes 2})^{-1/2}\end{pmatrix}\displaystyle\sum_{i=1}^n \operatorname{vec}\left(\psi_{f_1}\left(\frac{d_i}{\sigma}\right)\frac{d_i}{\sigma}\mathbf{U}_i\mathbf{U}_i' - \mathbf{I}_k\right)\end{pmatrix}\end{aligned} \tag{2.5}$$

*and full-rank information matrix*

$$\boldsymbol{\Gamma}_{f_1}(\boldsymbol{\vartheta}) := \begin{pmatrix} \boldsymbol{\Gamma}_{f_1;11}(\boldsymbol{\vartheta}) & \mathbf{0} & \mathbf{0} \\ \mathbf{0} & \Gamma_{f_1;22}(\boldsymbol{\vartheta}) & \boldsymbol{\Gamma}'_{f_1;32}(\boldsymbol{\vartheta}) \\ \mathbf{0} & \boldsymbol{\Gamma}_{f_1;32}(\boldsymbol{\vartheta}) & \boldsymbol{\Gamma}_{f_1;33}(\boldsymbol{\vartheta}) \end{pmatrix}, \tag{2.6}$$

*where*

$$\boldsymbol{\Gamma}_{f_1;11}(\boldsymbol{\vartheta}) := \frac{1}{k\sigma^2}\mathcal{I}_k(f_1)\mathbf{V}^{-1},$$

$$\Gamma_{f_1;22}(\boldsymbol{\vartheta}) := \frac{1}{4\sigma^4}(\mathcal{J}_k(f_1) - k^2),$$

$$\boldsymbol{\Gamma}_{f_1;32}(\boldsymbol{\vartheta}) := \frac{1}{4k\sigma^2}(\mathcal{J}_k(f_1) - k^2)\mathbf{M}_k \operatorname{vec}(\mathbf{V}^{-1})$$

*and*

$$\begin{aligned} \boldsymbol{\Gamma}_{f_1;33}(\boldsymbol{\vartheta}) := {} & \frac{1}{4}\mathbf{M}_k(\mathbf{V}^{\otimes 2})^{-1/2}\left[\frac{\mathcal{J}_k(f_1)}{k(k+2)}(\mathbf{I}_{k^2} + \mathbf{K}_k + \mathbf{J}_k) - \mathbf{J}_k\right] \\ & \times (\mathbf{V}^{\otimes 2})^{-1/2}\mathbf{M}'_k. \end{aligned} \tag{2.7}$$

*More precisely, for any* $\boldsymbol{\vartheta}^{(n)} = (\boldsymbol{\theta}^{(n)\prime}, \sigma^{2(n)}, (\operatorname{ve\mathring{c}h} \mathbf{V}^{(n)})')' = \boldsymbol{\vartheta} + O(n^{-1/2})$ *and any bounded sequence* $\boldsymbol{\tau}^{(n)} := (\mathbf{t}^{(n)\prime}, s^{(n)}, (\operatorname{ve\mathring{c}h} \mathbf{v}^{(n)})')' = (\boldsymbol{\tau}_1^{(n)\prime}, \tau_2^{(n)}, \boldsymbol{\tau}_3^{(n)\prime})'$ *in* $\mathbb{R}^{k+k(k+1)/2}$, *we have*

$$\begin{aligned} \Lambda^{(n)}_{\boldsymbol{\vartheta}^{(n)}+n^{-1/2}\boldsymbol{\tau}^{(n)}/\boldsymbol{\vartheta}^{(n)};f_1} &:= \log(d\mathrm{P}^{(n)}_{\boldsymbol{\vartheta}^{(n)}+n^{-1/2}\boldsymbol{\tau}^{(n)};f_1}/d\mathrm{P}^{(n)}_{\boldsymbol{\vartheta}^{(n)};f_1}) \\ &= (\boldsymbol{\tau}^{(n)})'\boldsymbol{\Delta}^{(n)}_{f_1}(\boldsymbol{\vartheta}^{(n)}) - \tfrac{1}{2}(\boldsymbol{\tau}^{(n)})'\boldsymbol{\Gamma}_{f_1}(\boldsymbol{\vartheta})\boldsymbol{\tau}^{(n)} + o_{\mathrm{P}}(1) \end{aligned}$$

*and*

$$\boldsymbol{\Delta}^{(n)}_{f_1}(\boldsymbol{\vartheta}^{(n)}) \xrightarrow{\mathcal{L}} \mathcal{N}(\mathbf{0}, \boldsymbol{\Gamma}_{f_1}(\boldsymbol{\vartheta}))$$

*under* $\mathrm{P}^{(n)}_{\boldsymbol{\vartheta}^{(n)};f_1}$, *as* $n \to \infty$.

PROOF. See Appendix (Section A.1). □

Note that the structure of the information matrix for shape (2.7) is not unfamiliar, having been previously obtained under much more restrictive assumptions; see, for example, page 219 of [8].

**3. Parametrically efficient tests for shape.**

3.1. *An efficient central sequence for shape.* The block-diagonal structure of the information matrix (2.6) and ULAN imply that substituting a (in principle, discretized—see, e.g., [30], page 125) root-$n$ consistent estimator $\hat{\boldsymbol{\theta}} = \hat{\boldsymbol{\theta}}^{(n)}$ for the unknown $\boldsymbol{\theta}$ has no influence, asymptotically, on the $(\sigma^2, \mathbf{V})$-part of the central sequence $\boldsymbol{\Delta}_{f_1}^{(n)}(\boldsymbol{\vartheta})$. Optimal inference about $(\sigma^2, \mathbf{V})$ can thus be based, without any loss of (asymptotic) efficiency, on $(\Delta_{f_1;2}^{(n)}(\hat{\boldsymbol{\theta}}, \sigma^2, \mathbf{V}), \boldsymbol{\Delta}_{f_1;3}^{(n)\prime}(\hat{\boldsymbol{\theta}}, \sigma^2, \mathbf{V}))'$ as if $\hat{\boldsymbol{\theta}}$ were the actual location parameter; see Section 4.4 for details. Therefore, in this section, we assume throughout that $\boldsymbol{\theta}$ is known. Similarly, replacing $\sigma^2$ and $\mathbf{V}$ with root-$n$ consistent estimators $\hat{\sigma}^{2(n)}$ and $\hat{\mathbf{V}}^{(n)}$ in the $\boldsymbol{\theta}$-part of the central sequence $\boldsymbol{\Delta}_{f_1}^{(n)}(\boldsymbol{\vartheta})$ has no impact, asymptotically, on inference about $\boldsymbol{\theta}$.

Unlike the asymptotic covariances between the location and scatter components of the central sequence $\boldsymbol{\Delta}_{f_1}^{(n)}(\boldsymbol{\vartheta})$, the asymptotic covariances between the $\sigma^2$-part $\Delta_{f_1;2}^{(n)}(\boldsymbol{\vartheta})$ and the $\mathbf{V}$-part $\boldsymbol{\Delta}_{f_1;3}^{(n)}(\boldsymbol{\vartheta})$ are not zero. This means that a local perturbation of scale has the same asymptotic impact on $\boldsymbol{\Delta}_{f_1;3}^{(n)}(\boldsymbol{\vartheta})$ as a local perturbation of $\mathbf{V}$. It follows that the cost of not knowing the actual value of $\sigma^2$ is strictly positive when performing inference on $\mathbf{V}$. Since it is hard to think of any practical problem where the scale (but not the shape) is specified, we concentrate on optimality under unspecified scale $\sigma^2$ and explicitly compute the information loss due to the presence of this nuisance.

LAN and the convergence of local experiments to the Gaussian shift experiment

$$\begin{pmatrix} \Delta_2 \\ \boldsymbol{\Delta}_3 \end{pmatrix} \sim \mathcal{N}\left( \begin{pmatrix} \Gamma_{f_1;22}(\boldsymbol{\vartheta}) & \boldsymbol{\Gamma}'_{f_1;32}(\boldsymbol{\vartheta}) \\ \boldsymbol{\Gamma}_{f_1;32}(\boldsymbol{\vartheta}) & \boldsymbol{\Gamma}_{f_1;33}(\boldsymbol{\vartheta}) \end{pmatrix} \begin{pmatrix} \tau_2 \\ \boldsymbol{\tau}_3 \end{pmatrix}, \begin{pmatrix} \Gamma_{f_1;22}(\boldsymbol{\vartheta}) & \boldsymbol{\Gamma}'_{f_1;32}(\boldsymbol{\vartheta}) \\ \boldsymbol{\Gamma}_{f_1;32}(\boldsymbol{\vartheta}) & \boldsymbol{\Gamma}_{f_1;33}(\boldsymbol{\vartheta}) \end{pmatrix} \right), \quad (\tau_2, \boldsymbol{\tau}'_3)' \in \mathbb{R}^{k(k+1)/2}, \tag{3.1}$$

imply that locally optimal inference on shape, in the presence of an unspecified scale parameter, should be based on the residual of the regression [in (3.1)] of $\boldsymbol{\Delta}_3$ with respect to $\Delta_2$, computed at $\boldsymbol{\Delta}_{f_1;3}^{(n)}(\boldsymbol{\vartheta})$ (the shape part of the central sequence) and $\Delta_{f_1;2}^{(n)}(\boldsymbol{\vartheta})$ (the scale part of the same). This residual takes the form $\boldsymbol{\Delta}_3 - \boldsymbol{\Gamma}_{f_1;32}(\boldsymbol{\vartheta})\Gamma_{f_1;22}^{-1}(\boldsymbol{\vartheta})\Delta_2(\boldsymbol{\vartheta})$; the resulting $f_1$-*efficient central sequence for shape* is thus

$$\boldsymbol{\Delta}_{f_1}^{\star(n)}(\boldsymbol{\vartheta}) = \boldsymbol{\Delta}_{f_1;3}^{(n)}(\boldsymbol{\vartheta}) - \boldsymbol{\Gamma}_{f_1;32}(\boldsymbol{\vartheta})\Gamma_{f_1;22}^{-1}(\boldsymbol{\vartheta})\Delta_{f_1;2}^{(n)}(\boldsymbol{\vartheta}),$$

which, after some elementary algebra, reduces to

$$\boldsymbol{\Delta}_{f_1}^{\star(n)}(\boldsymbol{\vartheta}) = \frac{1}{2} n^{-1/2} \mathbf{M}_k (\mathbf{V}^{\otimes 2})^{-1/2} \mathbf{J}_k^{\perp} \sum_{i=1}^{n} \psi_{f_1}\left(\frac{d_i}{\sigma}\right) \frac{d_i}{\sigma} \operatorname{vec}(\mathbf{U}_i \mathbf{U}_i').$$

This efficient central sequence under $\mathrm{P}^{(n)}_{\boldsymbol{\vartheta},f_1}$ is asymptotically normal, with mean zero and covariance (the *efficient Fisher information for shape* under radial density $f_1$) given by

$$\boldsymbol{\Gamma}^{\star}_{f_1}(\boldsymbol{\vartheta}) = \boldsymbol{\Gamma}_{f_1;33}(\boldsymbol{\vartheta}) - \boldsymbol{\Gamma}_{f_1;32}(\boldsymbol{\vartheta})\Gamma^{-1}_{f_1;22}(\boldsymbol{\vartheta})\boldsymbol{\Gamma}'_{f_1;32}(\boldsymbol{\vartheta}).$$

After some routine computation, this efficient information takes the form

$$\begin{aligned}\boldsymbol{\Gamma}^{\star}_{f_1}(\boldsymbol{\vartheta}) &= \frac{1}{4}\mathbf{M}_k(\mathbf{V}^{\otimes 2})^{-1/2}\mathbf{J}_k^{\perp} \\ &\quad \times \left[\frac{\mathcal{J}_k(f_1)}{k(k+2)}(\mathbf{I}_{k^2}+\mathbf{K}_k+\mathbf{J}_k)-\mathbf{J}_k\right]\mathbf{J}_k^{\perp}(\mathbf{V}^{\otimes 2})^{-1/2}\mathbf{M}_k' \\ &= \frac{\mathcal{J}_k(f_1)}{4k(k+2)}\mathbf{M}_k(\mathbf{V}^{\otimes 2})^{-1/2}\left[\mathbf{I}_{k^2}+\mathbf{K}_k-\frac{2}{k}\mathbf{J}_k\right](\mathbf{V}^{\otimes 2})^{-1/2}\mathbf{M}_k' \\ &=: \mathcal{J}_k(f_1)\boldsymbol{\Upsilon}_k^{-1}(\mathbf{V}),\end{aligned} \tag{3.2}$$

a form that is not unfamiliar in the area of robust estimation of covariance matrices; see, for instance, the asymptotic covariances in [40, 42, 50, 51] for the covariances of scatter estimates [as in (2.6), (2.7)], [41, 52] for covariances of shape estimates [as in (3.2)].

In the sequel, *optimality* (in the local and asymptotic sense, at radial density $f_1$) is to be understood in the context of the Gaussian shift experiments associated with the efficient central sequences $\boldsymbol{\Delta}^{\star(n)}_{f_1}(\boldsymbol{\vartheta})$. In particular, a sequence of tests will be called locally and asymptotically *maximin-efficient* (at asymptotic level $\alpha$) if it is asymptotically maximin in the sequence of experiments associated with $\boldsymbol{\Delta}^{\star(n)}_{f_1}(\boldsymbol{\vartheta})$.

3.2. *Optimal parametric tests for shape.* Consider the problem of testing a null hypothesis of the form $\mathcal{H}_0 : \mathbf{V} = \mathbf{V}_0$ in the parametric model where $f_1$ is known and the scale $\sigma^2$ is unspecified. Optimality (in a local and asymptotic sense—see Proposition 3.1 for a precise statement) is reached by tests based on quadratic forms in the $f_1$-efficient central sequence for shape. More precisely, the optimal test statistics take the form

$$Q_{f_1} = Q^{(n)}_{f_1} := (\boldsymbol{\Delta}^{\star(n)}_{f_1}(\hat{\boldsymbol{\vartheta}}_0))'(\boldsymbol{\Gamma}^{\star}_{f_1}(\hat{\boldsymbol{\vartheta}}_0))^{-1}\boldsymbol{\Delta}^{\star(n)}_{f_1}(\hat{\boldsymbol{\vartheta}}_0),$$

where, denoting by $\hat{\sigma}$ a root-$n$ consistent estimator for $\sigma$, we let $\hat{\boldsymbol{\vartheta}}_0 := (\boldsymbol{\theta}', \hat{\sigma}^2, (\overset{\circ}{\mathrm{vech}}\,\mathbf{V}_0)')'$. Note that consistent estimation of $\sigma$ under the family $\bigcup_{f_1}\bigcup_{\sigma>0}\{\mathrm{P}^{(n)}_{\boldsymbol{\theta},\sigma^2,\mathbf{V}_0;f_1}\}$ is easily achieved since $\sigma$ is then the common median of the distances $d_i(\boldsymbol{\theta}, \mathbf{V}_0)$. As we shall see in Section 3.3, the Gaussian version of these optimal parametric tests allows the bypassing of this estimation of $\sigma$.

Lemma 3.1 below leads to the more explicit form

$$Q_{f_1} = \frac{k(k+2)}{2n\mathcal{J}_k(f_1)} \sum_{i,j=1}^{n} \frac{d_i d_j}{\hat{\sigma}^2} \psi_{f_1}\left(\frac{d_i}{\hat{\sigma}}\right) \psi_{f_1}\left(\frac{d_j}{\hat{\sigma}}\right) \left( (\mathbf{U}_i' \mathbf{U}_j)^2 - \frac{1}{k} \right),$$

with $d_i := d_i^{(n)}(\boldsymbol{\theta}, \mathbf{V}_0)$ and $\mathbf{U}_i := \mathbf{U}_i^{(n)}(\boldsymbol{\theta}, \mathbf{V}_0)$.

LEMMA 3.1. *Denote by $\mathbf{e}_{k^2,1}$ the first vector of the canonical basis of $\mathbb{R}^{k^2}$. Then if $\mathbf{V} = (V_{ij})$ is symmetric with $V_{11} = 1$, we have*

$$\begin{aligned} &\frac{1}{k(k+2)} \mathbf{M}_k' \boldsymbol{\Upsilon}_k(\mathbf{V}) \mathbf{M}_k \\ &\quad = [\mathbf{I}_{k^2} + \mathbf{K}_k](\mathbf{V}^{\otimes 2}) - 2(\mathbf{V}^{\otimes 2}) \mathbf{e}_{k^2,1} (\operatorname{vec} \mathbf{V})' \\ &\qquad - 2(\operatorname{vec} \mathbf{V})(\mathbf{e}_{k^2,1})'(\mathbf{V}^{\otimes 2}) + 2(\operatorname{vec} \mathbf{V})(\operatorname{vec} \mathbf{V})'. \end{aligned} \tag{3.3}$$

PROOF. See Appendix (Section A.2). □

PROPOSITION 3.1. *Let $f_1$ satisfy Assumptions (A1) and (A2). Then, denoting by $\|\mathbf{A}\| := [\operatorname{tr}(\mathbf{A}\mathbf{A}')]^{1/2}$ the Frobenius norm of $\mathbf{A}$,*

(i) *$Q_{f_1}^{(n)}$ is asymptotically chi-square with $k(k+1)/2 - 1$ degrees of freedom under $\bigcup_{\sigma^2} \{\mathrm{P}^{(n)}_{\boldsymbol{\theta},\sigma^2,\mathbf{V}_0;f_1}\}$ and asymptotically noncentral chi-square, still with $k(k+1)/2 - 1$ degrees of freedom but with noncentrality parameter*

$$\begin{aligned} &\frac{\mathcal{J}_k(f_1)}{2k(k+2)} \left[ \operatorname{tr}((\mathbf{V}_0^{-1}\mathbf{v})^2) - \frac{1}{k} (\operatorname{tr} \mathbf{V}_0^{-1}\mathbf{v})^2 \right] \\ &\quad = \frac{\mathcal{J}_k(f_1)}{2k(k+2)} (\operatorname{tr} \mathbf{V}_0^{-1}\mathbf{v})^2 \left\| \frac{\mathbf{V}_0^{-1}\mathbf{v}}{\operatorname{tr} \mathbf{V}_0^{-1}\mathbf{v}} - \frac{1}{k}\mathbf{I}_k \right\|^2, \end{aligned}$$

*under $\bigcup_{\sigma^2} \{\mathrm{P}^{(n)}_{\boldsymbol{\theta},\sigma^2,\mathbf{V}_0 + n^{-1/2}\mathbf{v};f_1}\}$;*

(ii) *the sequence of tests $\phi_{f_1}^{(n)}$ which consists of rejecting $\mathcal{H}_0 : \mathbf{V} = \mathbf{V}_0$ as soon as $Q_{f_1}^{(n)}$ exceeds the $\alpha$ upper-quantile of a chi-square variable with $k(k+1)/2 - 1$ degrees of freedom, has asymptotic level $\alpha$ under $\bigcup_{\sigma^2} \{\mathrm{P}^{(n)}_{\boldsymbol{\theta},\sigma^2,\mathbf{V}_0;f_1}\}$ and is locally and asymptotically maximin-efficient, still at asymptotic level $\alpha$, for $\bigcup_{\sigma^2} \{\mathrm{P}^{(n)}_{\boldsymbol{\theta},\sigma^2,\mathbf{V}_0;f_1}\}$ against alternatives of the form $\bigcup_{\sigma^2} \bigcup_{\mathbf{V} \neq \mathbf{V}_0} \{\mathrm{P}^{(n)}_{\boldsymbol{\theta},\sigma^2,\mathbf{V};f_1}\}$.*

PROOF. See Appendix (Section A.2). □

In contrast with this unspecified-$\sigma^2$ test, the locally and asymptotically optimal procedure for testing $\mathcal{H}_0 : \mathbf{V} = \mathbf{V}_0$ under specified radial density $f_1$, specified $\boldsymbol{\theta}$ *and specified scale* $\sigma^2$ rejects $\mathcal{H}_0$ (at asymptotic level $\alpha$) whenever

$$\begin{aligned} Q_{\sigma^2,f_1} &= Q^{(n)}_{\sigma^2,f_1} \\ &:= (\boldsymbol{\Delta}^{(n)}_{f_1;3}(\boldsymbol{\theta},\sigma^2,\mathbf{V}_0))'(\boldsymbol{\Gamma}_{f_1;33}(\boldsymbol{\theta},\sigma^2,\mathbf{V}_0))^{-1}\boldsymbol{\Delta}^{(n)}_{f_1;3}(\boldsymbol{\theta},\sigma^2,\mathbf{V}_0) \end{aligned} \tag{3.4}$$

exceeds the $\alpha$ upper-quantile of a chi-square with $k(k+1)/2-1$ degrees of freedom. The efficiency loss due to an unspecified $\sigma^2$ can thus be measured by the difference between the noncentrality parameters in the asymptotic chi-square distributions of $Q_{\sigma^2,f_1}$ and $Q_{f_1}$ under local alternatives. Along the same lines as the proof of Proposition 3.1, one can show that this difference, under $\mathrm{P}^{(n)}_{\boldsymbol{\theta},\sigma^2,\mathbf{V}_0+n^{-1/2}\mathbf{v};f_1}$, is

$$\frac{1}{4k^2}(\mathcal{J}_k(f_1)-k^2)(\operatorname{tr}\mathbf{V}_0^{-1}\mathbf{v})^2. \tag{3.5}$$

Inequality (2.4) confirms the unsurprising fact that this loss is nonnegative and an increasing function of the information for shape (or scale) $\mathcal{J}_k(f_1)$. Quite remarkably, it does not depend on the scale $\sigma^2$ itself. Also, note that the loss is nil against local alternatives such that $\operatorname{tr}\mathbf{V}_0^{-1}\mathbf{v}=0$. When testing for sphericity ($\mathbf{V}_0=\mathbf{I}_k$), this reduces to $\operatorname{tr}\mathbf{v}=0$; in particular, there is no loss in the case $v_{ii}=0$ for all $i=2,\dots,k$.

Further investigation of (3.5) reveals some interesting facts concerning the relation between this loss and the tails of underlying radial densities. Assume, for the sake of simplicity, that $\mathbf{V}_0=\mathbf{I}_k$ and consider the "elementary diagonal deviations from sphericity" associated with $\mathbf{v}=\lambda\mathbf{e}_i\mathbf{e}_i'$ for some $i=2,\dots,k$. The relative loss in local powers (strictly speaking, the relative loss in the corresponding noncentrality parameters) can be evaluated as the ratio of (3.5) and the noncentrality parameter one would obtain for the specified-$\sigma$ test statistic (3.4)—namely, the sum of (3.5) and the noncentrality parameter in Proposition 3.1(i). This relative loss no longer depends on $\lambda$ and takes the form

$$\frac{(k+2)(\mathcal{J}_k(f_1)-k^2)}{3k(\mathcal{J}_k(f_1)-k^2)+2k^2(k-1)}, \tag{3.6}$$

an increasing function of $\mathcal{J}_k(f_1)$, with lower and upper bounds 0 and $(k+2)/3k$, corresponding to arbitrarily heavy- and arbitrary light-tailed distributions, respectively. Indeed, these bounds can be obtained, for example, by letting $\eta\to 0$ and $\eta\to\infty$, respectively, in the power-exponential family of distributions considered in Section 1.2. We refer to [18] for more general results on efficiency losses in the related problem of estimating the shape parameter.

Some numerical values of those relative losses (3.6) are provided in Table 1 where we consider:

TABLE 1

*Numerical values of the relative power losses in (3.6) under k-variate power-exponential densities (with $\eta = 0.1$, 0.5, 1, 2, 5, along with the limiting values obtained for $\eta \to 0$ and $\eta \to \infty$), and under k-variate Student densities (with $\nu$ degrees of freedom, $\nu = 1$, 3, 5, 8, 15, along with the limiting values obtained for $\nu \to 0$ and $\nu \to \infty$), for $k = 2$, 3, 4, 6, 10 and for $k \to \infty$*

| | **Parameter $\eta$ of the power-exponential density** | | | | | | |
|---|---|---|---|---|---|---|---|
| $k$ | $\to 0$ | **0.1** | **0.5** | **1** | **2** | **5** | $\to \infty$ |
| 2 | 0.000 | 0.154 | 0.400 | 0.500 | 0.571 | 0.625 | 0.667 |
| 3 | 0.000 | 0.072 | 0.238 | 0.333 | 0.417 | 0.490 | 0.556 |
| 4 | 0.000 | 0.045 | 0.167 | 0.250 | 0.333 | 0.417 | 0.500 |
| 6 | 0.000 | 0.025 | 0.103 | 0.167 | 0.242 | 0.333 | 0.444 |
| 10 | 0.000 | 0.013 | 0.057 | 0.100 | 0.160 | 0.250 | 0.400 |
| $\infty$ | 0.000 | 0.000 | 0.000 | 0.000 | 0.000 | 0.000 | 0.333 |
| | **Degrees of freedom of the underlying $t$ density** | | | | | | |
| | $\to 0$ | **1** | **3** | **5** | **8** | **15** | $\to \infty$ |
| 2 | 0.000 | 0.250 | 0.375 | 0.417 | 0.444 | 0.469 | 0.500 |
| 3 | 0.000 | 0.111 | 0.200 | 0.238 | 0.267 | 0.294 | 0.333 |
| 4 | 0.000 | 0.063 | 0.125 | 0.156 | 0.182 | 0.208 | 0.250 |
| 6 | 0.000 | 0.028 | 0.063 | 0.083 | 0.103 | 0.125 | 0.167 |
| 10 | 0.000 | 0.010 | 0.025 | 0.036 | 0.047 | 0.063 | 0.100 |
| $\infty$ | 0.000 | 0.000 | 0.000 | 0.000 | 0.000 | 0.000 | 0.000 |

(a) the family of power-exponential densities (providing a full range of tail behaviors), with relative loss $(k+2)\eta/(k(k+3\eta-1))$ and

(b) the more familiar heavy-tailed Student densities with $\nu$ degrees of freedom (including the Gaussian as $\nu \to \infty$), with relative loss $\nu/(k(k+\nu-1))$, respectively, for several values of the space dimension $k$. Limits as $k \to \infty$ are taken for fixed $\nu$ or $\eta$; note that the $\eta = 1$ power-exponential and $\nu = \infty$ Student columns both correspond to the Gaussian case.

3.3. *Optimal Gaussian tests for shape.* The parametric tests $\phi_{f_1}^{(n)}$ described in part (ii) of Proposition 3.1 achieve local and asymptotic optimality at radial density $f_1$ but are generally not valid when the underlying radial density is $g_1 \neq f_1$. If correctly formulated, the Gaussian version of these tests (obtained for $f_1 = \phi_1$, where $\phi_1$ was defined in Section 1.2) is an interesting exception to this rule and can easily be written in a form that remains valid under the class of all radial densities $g_1$ such that $\tilde{g}_{1k}$ has finite fourth-order moments.

Denote by $D_k(g_1) := \mathrm{E}[(\tilde{G}_{1k}^{-1}(U))^2]$ and $E_k(g_1) := \mathrm{E}[(\tilde{G}_{1k}^{-1}(U))^4] < \infty$, where $U$ stands for a random variable with uniform distribution over $(0,1)$, the second- and fourth-order moments of $\tilde{g}_{1k}$, respectively, and assume that $E_k(g_1) < \infty$ [hence also that $D_k(g_1) < \infty$]. These two quantities are closely related to the *kurtosis* of the elliptical distribution under consideration. To be precise, the kurtosis $3\kappa_k(g_1)$ of an elliptically symmetric random $k$-vector $\mathbf{X} = (X_i)$ with location center $\boldsymbol{\theta} = (\theta_1, \ldots, \theta_k)'$, scale $\sigma^2$, shape matrix $\mathbf{V}$, and radial density $g_1$ is defined to be

$$3\kappa_k(g_1) := \frac{\mathrm{E}[(X_i - \theta_i)^4]}{\mathrm{E}^2[(X_i - \theta_i)^2]} - 3;$$

see, for example, [1], page 54, [38] or [50]. This quantity depends only on the dimension $k$ and the radial density $g_1$, not on $i$ or on the other parameters characterizing the elliptical distribution (which of course justifies the notation); it is related to $D_k(g_1)$ and $E_k(g_1)$ by the simple relation

$$\kappa_k(g_1) = \frac{k}{k+2} \frac{E_k(g_1)}{D_k^2(g_1)} - 1.$$

At the $k$-variate Gaussian distribution and $t$-distribution with $\nu$ degrees of freedom ($\nu > 4$), this kurtosis parameter takes values $\kappa_k(\phi_1) = 0$ and $\kappa_k(f_{1,\nu}^t) = 2/(\nu - 4)$, respectively.

The Gaussian version of the efficient central sequence for shape $\boldsymbol{\Delta}_{f_1}^{\star(n)}(\boldsymbol{\vartheta})$ can be written as $\boldsymbol{\Delta}_{\phi_1}^{\star(n)}(\boldsymbol{\vartheta}) = a_k \sigma^{-2} \mathbf{T}_{\boldsymbol{\theta},\mathbf{V}}$, where

$$\mathbf{T}_{\boldsymbol{\theta},\mathbf{V}} = \mathbf{T}_{\boldsymbol{\theta},\mathbf{V}}^{(n)} := \tfrac{1}{2} n^{-1/2} \mathbf{M}_k (\mathbf{V}^{\otimes 2})^{-1/2} \mathbf{J}_k^{\perp} (\mathbf{V}^{\otimes 2})^{-1/2} \sum_{i=1}^{n} \operatorname{vec}((\mathbf{X}_i - \boldsymbol{\theta})(\mathbf{X}_i - \boldsymbol{\theta})').$$

It is convenient to work with $\mathbf{T}_{\boldsymbol{\theta},\mathbf{V}}$ and an estimate $\hat{\boldsymbol{\Gamma}}^{(n)}$ of its asymptotic covariance rather than with $\boldsymbol{\Delta}_{\phi_1}^{\star(n)}(\boldsymbol{\vartheta})$ and an estimate of the corresponding information matrix since the scalar factor $a_k \sigma^{-2}$ in the quadratic form in $\boldsymbol{\Delta}_{\phi_1}^{\star(n)}(\boldsymbol{\vartheta})$ cancels out. For optimality (at Gaussian radial densities), it is sufficient for $\hat{\boldsymbol{\Gamma}}^{(n)}$ to consistently estimate the asymptotic covariance of $\mathbf{T}_{\boldsymbol{\theta},\mathbf{V}_0}$ under $\bigcup_{\sigma^2} \{\mathrm{P}_{\boldsymbol{\theta},\sigma^2,\mathbf{V}_0;\phi_1}^{(n)}\}$.

Letting

$$\hat{\boldsymbol{\Gamma}}^{(n)} := \left( \frac{1}{n} \sum_{i=1}^{n} d_i^4 \right) \boldsymbol{\Upsilon}_k^{-1}(\mathbf{V}),$$

with the same $d_i = d_i^{(n)}(\boldsymbol{\theta}, \mathbf{V}_0)$'s as in Section 3.2, it is easy to check that $\hat{\boldsymbol{\Gamma}}^{(n)}$ provides, for all $\boldsymbol{\theta}$, a consistent estimate of the asymptotic variance of $\mathbf{T}_{\boldsymbol{\theta},\mathbf{V}_0}$, not only under $\bigcup_{\sigma^2} \{\mathrm{P}_{\boldsymbol{\theta},\sigma^2,\mathbf{V}_0;\phi_1}^{(n)}\}$, but also under $\bigcup_{\sigma^2} \bigcup_{g_1} \{\mathrm{P}_{\boldsymbol{\theta},\sigma^2,\mathbf{V}_0;g_1}^{(n)}\}$, where the union is taken over the set of all densities $g_1$ such that $E_k(g_1) < \infty$. The Gaussian test statistic then takes the form $Q_{\mathcal{N}} = Q_{\mathcal{N}}^{(n)} := \mathbf{T}_{\boldsymbol{\theta},\mathbf{V}_0}^{(n)\prime}(\hat{\boldsymbol{\Gamma}}^{(n)})^{-1}\mathbf{T}_{\boldsymbol{\theta},\mathbf{V}_0}^{(n)}$. Lemma 3.1 and standard algebra yield

$$Q_{\mathcal{N}} = \frac{k(k+2)}{2(\sum_{i=1}^n d_i^4)} \sum_{i,j=1}^n d_i^2 d_j^2 \left( (\mathbf{U}_i'\mathbf{U}_j)^2 - \frac{1}{k} \right), \tag{3.7}$$

with the same $\mathbf{U}_i = \mathbf{U}_i^{(n)}(\boldsymbol{\theta}, \mathbf{V}_0)$ as in Section 3.2. Now, defining

$$\mathbf{S} = \mathbf{S}^{(n)} := \frac{1}{n} \sum_{i=1}^n [\mathbf{V}_0^{-1/2}(\mathbf{X}_i - \boldsymbol{\theta})][\mathbf{V}_0^{-1/2}(\mathbf{X}_i - \boldsymbol{\theta})]'$$

and letting $\hat{\kappa}^{(n)} := [k(n^{-1}\sum_{i=1}^n d_i^4)]/[(k+2)(n^{-1}\sum_{i=1}^n d_i^2)^2] - 1$ be a consistent estimate of the kurtosis parameter $\kappa_k(g_1)$, (3.7) takes the form

$$Q_{\mathcal{N}} = \frac{n^2 k(k+2)}{2(\sum_{i=1}^n d_i^4)} \left( \operatorname{tr}\mathbf{S}^2 - \frac{1}{k}\operatorname{tr}^2\mathbf{S} \right) = \frac{1}{1+\hat{\kappa}^{(n)}} \frac{nk^2}{2} \left\| \frac{\mathbf{S}}{\operatorname{tr}\mathbf{S}} - \frac{1}{k}\mathbf{I}_k \right\|^2. \tag{3.8}$$

It is straightforward to check that $Q_{\mathcal{N}}$ is invariant under rotations, scale transformations and reflections (with respect to $\boldsymbol{\theta}$, in the metric associated with $\mathbf{V}_0$), but that it is not (even asymptotically) invariant under the group of monotone continuous radial transformations (see Section 4.1 below). The following proposition summarizes the asymptotic properties of the Gaussian procedure based on $Q_{\mathcal{N}}$:

PROPOSITION 3.2. *Denote by $\phi_{\mathcal{N}}^{(n)}$ the parametric Gaussian test rejecting the null hypothesis $\mathcal{H}_0 : \mathbf{V} = \mathbf{V}_0$ whenever $Q_{\mathcal{N}}^{(n)}$ exceeds the $\alpha$ upper-quantile of a chi-square distribution with $k(k+1)/2 - 1$ degrees of freedom. Then (unions over $g_1$ are taken over all densities such that $\tilde{g}_{1k}$ has finite fourth-order moments):*

(i) *$Q_{\mathcal{N}}^{(n)}$ is asymptotically chi-square with $k(k+1)/2-1$ degrees of freedom under $\bigcup_{\sigma^2}\bigcup_{g_1}\{\mathrm{P}_{\boldsymbol{\theta},\sigma^2,\mathbf{V}_0;g_1}^{(n)}\}$ and asymptotically noncentral chi-square, still with $k(k+1)/2-1$ degrees of freedom, but with noncentrality parameter*

$$\frac{1}{2(1+\kappa_k(g_1))}\left[\operatorname{tr}((\mathbf{V}_0^{-1}\mathbf{v})^2) - \frac{1}{k}(\operatorname{tr}\mathbf{V}_0^{-1}\mathbf{v})^2\right]$$

*under $\bigcup_{\sigma^2}\{\mathrm{P}_{\boldsymbol{\theta},\sigma^2,\mathbf{V}_0+n^{-1/2}\mathbf{v};g_1}^{(n)}\}$;*

(ii) *the sequence of tests $\phi_{\mathcal{N}}^{(n)}$ under $\bigcup_{\sigma^2}\bigcup_{g_1}\{\mathrm{P}_{\boldsymbol{\theta},\sigma^2,\mathbf{V}_0;g_1}^{(n)}\}$ has asymptotic level $\alpha$ and is locally and asymptotically maximin-efficient, still at asymptotic level $\alpha$, for $\bigcup_{\sigma^2}\bigcup_{g_1}\{\mathrm{P}_{\boldsymbol{\theta},\sigma^2,\mathbf{V}_0;g_1}^{(n)}\}$ against alternatives of the form $\bigcup_{\sigma^2}\bigcup_{\mathbf{V}\neq\mathbf{V}_0}\{\mathrm{P}_{\boldsymbol{\theta},\sigma^2,\mathbf{V};\phi_1}^{(n)}\}$.*

PROOF. See Appendix (Section A.3). □

For $\mathbf{V}_0 = \mathbf{I}_k$, the test statistic $Q_{\mathcal{N}}$ in (3.8) and Proposition 3.2 actually appears as a modification of the test statistic

$$Q_{\text{John}} := \frac{nk^2}{2} \left\| \frac{\mathbf{S}}{\operatorname{tr} \mathbf{S}} - \frac{1}{k} \mathbf{I}_k \right\|^2 = \frac{nk^2}{2} \operatorname{tr} \left[ \left( \frac{\mathbf{S}}{\operatorname{tr} \mathbf{S}} - \frac{1}{k} \mathbf{I}_k \right)^2 \right] \tag{3.9}$$

proposed by John [24, 25]. The only difference is that $Q_{\text{John}}$ relies on the Gaussian value $\kappa = 0$ of the kurtosis parameter, whereas $Q_{\mathcal{N}}$ instead involves an estimation $\hat{\kappa}^{(n)}$ of the same, which makes the asymptotic null distribution of $Q_{\mathcal{N}}$ agree, under any elliptical distribution with finite fourth-order moments, with the limiting distribution of $Q_{\text{John}}$ in the multinormal case.

This adjustment is very much in the spirit of the Muirhead and Waternaux version [38] of Mauchly's Gaussian likelihood ratio test [36]—probably the most widely used test of sphericity. Muirhead and Waternaux [38] actually show that the limiting distribution of $(-2 \log \Lambda^{(n)})/(1 + \kappa_k(g_1))$, where $-2 \log \Lambda^{(n)}$ is the Gaussian likelihood ratio test statistic, is asymptotically chi-square, with $k(k+1)/2 - 1$ degrees of freedom, under $\bigcup_{\sigma^2} \bigcup_{g_1} \{\mathrm{P}^{(n)}_{\boldsymbol{\theta}, \sigma^2, \mathbf{I}_k; g_1}\}$ (the union is taken over all $g_1$ such that $\tilde{g}_{1k}$ has finite fourth-order moments); the population kurtosis parameter $\kappa_k(g_1)$ can of course be replaced by its sample counterpart $\hat{\kappa}^{(n)}$ without modifying the asymptotic chi-square distribution. These results straightforwardly extend to the problem of testing for a specified shape $\mathbf{V}_0$ rather than for sphericity. It also follows from [38] that the adjusted version of John's test statistic, namely our Gaussian test statistic $Q_{\mathcal{N}}$, is asymptotically equivalent to their adjusted version of the Mauchly test. In the sequel, the expression "optimal parametric Gaussian test" will refer to any of these tests. Note, however, that optimality here follows from Proposition 3.2 and is therefore of an asymptotic nature. Actually, only John's original (nonadjusted) test [24] enjoys some finite-sample optimality properties (restricted to the Gaussian case), being *locally most powerful invariant* at the multinormal distribution. Our adjusted tests inherit, under weaker asymptotic form, this optimality property from John's test; on the other hand, they remain valid under non-Gaussian densities, which is not the case for John's.

## 4. Rank-based tests for shape.

4.1. *Rank-based versions of efficient central sequences for shape.* As already mentioned, the problem with tests based on efficient central sequences is that (with the exception of the adjusted Gaussian tests described in Section 3.3) they are only valid under correctly specified radial densities. In practice, a correct specification $f_1$ of the actual radial density $g_1$ is rather unrealistic and thus the problem has to be treated from a semiparametric point of view, where $g_1$ plays the role of a nuisance.

Within the family of distributions $\bigcup_{\sigma^2}\bigcup_{\mathbf{V}}\bigcup_{g_1}\{\mathrm{P}^{(n)}_{\boldsymbol{\theta},\sigma^2,\mathbf{V};g_1}\}$, where $\boldsymbol{\theta}$ is fixed, consider the null hypothesis $\mathcal{H}_0(\boldsymbol{\theta},\mathbf{V}_0)$ under which $\mathbf{V}=\mathbf{V}_0$. Throughout, therefore, $\boldsymbol{\theta}$ and $\mathbf{V}=\mathbf{V}_0$ are fixed, and $\sigma^2$ and the radial density $g_1$ remain unspecified (no moment assumptions are being made here). As we have seen, the scalar nuisance $\sigma^2$ can be taken care of by means of a simple projection, yielding the efficient central sequence. In principle, the infinite-dimensional nuisance $g_1$ can be treated similarly, by projecting central sequences along adequate *tangent spaces*; see Example 4 of [7]. This approach is rather technical, however. Hallin and Werker [20] showed that appropriate group invariance structures allow for the same result by conditioning central sequences with respect to maximal invariants such as ranks or signs. This is the approach we also adopt here.

Clearly, the null hypothesis $\mathcal{H}_0(\boldsymbol{\theta},\mathbf{V}_0)$ is invariant under the following two groups of transformations, acting on the observations $\mathbf{X}_1,\ldots,\mathbf{X}_n$:

(i) the group $\mathcal{G}^{\mathrm{orth}(n)},\circ := \mathcal{G}^{\mathrm{orth}(n)}_{\boldsymbol{\theta},\mathbf{V}_0},\circ$ of $\mathbf{V}_0$-*orthogonal transformations* (centered at $\boldsymbol{\theta}$) consisting of all transformations of the form

$$\begin{aligned}\mathbf{X}\mapsto{}& g_{\mathbf{O}}(\mathbf{X}_1,\ldots,\mathbf{X}_n)\\ &= g_{\mathbf{O}}(\boldsymbol{\theta}+d_1(\boldsymbol{\theta},\mathbf{V}_0)\mathbf{V}_0^{1/2}\mathbf{U}_1(\boldsymbol{\theta},\mathbf{V}_0),\ldots,\boldsymbol{\theta}+d_n(\boldsymbol{\theta},\mathbf{V}_0)\mathbf{V}_0^{1/2}\mathbf{U}_n(\boldsymbol{\theta},\mathbf{V}_0))\\ &:= (\boldsymbol{\theta}+d_1(\boldsymbol{\theta},\mathbf{V}_0)\mathbf{V}_0^{1/2}\mathbf{O}\mathbf{U}_1(\boldsymbol{\theta},\mathbf{V}_0),\ldots,\boldsymbol{\theta}+d_n(\boldsymbol{\theta},\mathbf{V}_0)\mathbf{V}_0^{1/2}\mathbf{O}\mathbf{U}_n(\boldsymbol{\theta},\mathbf{V}_0)),\end{aligned}$$

where $\mathbf{O}$ is an arbitrary $k\times k$ orthogonal matrix. In particular, this group contains "rotations" (in the metric associated with $\mathbf{V}_0$) around $\boldsymbol{\theta}$, as well as the *reflection* with respect to $\boldsymbol{\theta}$, that is, the mapping $(\mathbf{X}_1,\ldots,\mathbf{X}_n)\mapsto(\boldsymbol{\theta}-(\mathbf{X}_1-\boldsymbol{\theta}),\ldots,\boldsymbol{\theta}-(\mathbf{X}_n-\boldsymbol{\theta}))$;

(ii) the group $\mathcal{G}^{(n)},\circ := \mathcal{G}^{(n)}_{\boldsymbol{\theta},\mathbf{V}_0},\circ$ of *continuous monotone radial transformations*, of the form

$$\begin{aligned}\mathbf{X}\mapsto{}& g_h(\mathbf{X}_1,\ldots,\mathbf{X}_n)\\ &= g_h(\boldsymbol{\theta}+d_1(\boldsymbol{\theta},\mathbf{V}_0)\mathbf{V}_0^{1/2}\mathbf{U}_1(\boldsymbol{\theta},\mathbf{V}_0),\ldots,\boldsymbol{\theta}+d_n(\boldsymbol{\theta},\mathbf{V}_0)\mathbf{V}_0^{1/2}\mathbf{U}_n(\boldsymbol{\theta},\mathbf{V}_0))\\ &:= (\boldsymbol{\theta}+h(d_1(\boldsymbol{\theta},\mathbf{V}_0))\mathbf{V}_0^{1/2}\mathbf{U}_1(\boldsymbol{\theta},\mathbf{V}_0),\ldots,\boldsymbol{\theta}+h(d_n(\boldsymbol{\theta},\mathbf{V}_0))\mathbf{V}_0^{1/2}\mathbf{U}_n(\boldsymbol{\theta},\mathbf{V}_0)),\end{aligned}$$

where $h:\mathbb{R}^+\to\mathbb{R}^+$ is continuous, monotone increasing and such that $h(0)=0$ and $\lim_{r\to\infty}h(r)=\infty$. In particular, this group includes the subgroup of *scale transformations* $(\mathbf{X}_1,\ldots,\mathbf{X}_n)\mapsto(\boldsymbol{\theta}+a(\mathbf{X}_1-\boldsymbol{\theta}),\ldots,\ \boldsymbol{\theta}+a(\mathbf{X}_n-\boldsymbol{\theta}))$, $a>0$.

Clearly, the group $\mathcal{G}^{(n)},\circ$ of continuous monotone radial transformations is a generating group for the family of distributions $\bigcup_{\sigma^2}\bigcup_{f_1}\{\mathrm{P}^{(n)}_{\boldsymbol{\theta},\sigma^2,\mathbf{V}_0;f_1}\}$, that is, a generating group for the null hypothesis $\mathcal{H}_0(\boldsymbol{\theta}, \mathbf{V}_0)$ under consideration. The invariance principle therefore leads to the consideration of test statistics that are measurable with respect to the corresponding maximal invariant, namely the vector $(R_1(\boldsymbol{\theta}, \mathbf{V}_0), \ldots, R_n(\boldsymbol{\theta}, \mathbf{V}_0), \mathbf{U}_1(\boldsymbol{\theta}, \mathbf{V}_0), \ldots, \mathbf{U}_n(\boldsymbol{\theta}, \mathbf{V}_0))$, where $R_i(\boldsymbol{\theta}, \mathbf{V}_0)$ denotes the rank of $d_i(\boldsymbol{\theta}, \mathbf{V}_0)$ among $d_1(\boldsymbol{\theta}, \mathbf{V}_0), \ldots, d_n(\boldsymbol{\theta}, \mathbf{V}_0)$. The resulting signed rank test statistics are (strictly) invariant under $\mathcal{G}^{(n)}, \circ$, hence distribution-free under $\mathcal{H}_0(\boldsymbol{\theta}, \mathbf{V}_0)$.

Now, in the construction of the proposed tests for the null hypothesis $\mathcal{H}_0(\boldsymbol{\theta}, \mathbf{V}_0)$, we intend to combine invariance and optimality arguments by considering a (signed-)rank-based version of the $f_1$-efficient central sequences for shape [recall that central sequences are always defined up to $o_{\mathrm{P}}(1)$—under $\mathrm{P}^{(n)}_{\boldsymbol{\vartheta};f_1}$, as $n \to \infty$—terms]. The signed-rank version $\underset{\sim}{\boldsymbol{\Delta}}{}^{(n)}_{f_1}(\boldsymbol{\vartheta})$ of the shape-efficient central sequence $\boldsymbol{\Delta}^{\star(n)}_{f_1}(\boldsymbol{\vartheta})$ we plan to use in our nonparametric tests is the $f_1$-score version (based on the scores $K = K_{f_1}$) of the statistic

$$
\begin{aligned}
\underset{\sim}{\boldsymbol{\Delta}}{}^{(n)}_{K}(\boldsymbol{\vartheta}) &:= \frac{1}{2} n^{-1/2} \mathbf{M}_k (\mathbf{V}^{\otimes 2})^{-1/2} \mathbf{J}_k^{\perp} \sum_{i=1}^{n} K\left(\frac{R_i}{n+1}\right) \operatorname{vec}(\mathbf{U}_i \mathbf{U}_i') \\
&= \frac{1}{2} n^{-1/2} \mathbf{M}_k (\mathbf{V}^{\otimes 2})^{-1/2} \sum_{i=1}^{n} K\left(\frac{R_i}{n+1}\right) \operatorname{vec}\left(\mathbf{U}_i \mathbf{U}_i' - \frac{1}{k} \mathbf{I}_k\right) \\
&= \frac{1}{2} n^{-1/2} \mathbf{M}_k (\mathbf{V}^{\otimes 2})^{-1/2} \\
&\quad \times \sum_{i=1}^{n} \left( K\left(\frac{R_i}{n+1}\right) \operatorname{vec}(\mathbf{U}_i \mathbf{U}_i') - \frac{m_K^{(n)}}{k} \operatorname{vec}(\mathbf{I}_k) \right),
\end{aligned}
\tag{4.1}
$$

where $R_i = R_i^{(n)}(\boldsymbol{\theta}, \mathbf{V})$ denotes the rank of $d_i = d_i^{(n)}(\boldsymbol{\theta}, \mathbf{V})$ among $d_1, \ldots, d_n$, $\mathbf{U}_i = \mathbf{U}_i^{(n)}(\boldsymbol{\theta}, \mathbf{V})$ and $m_K^{(n)} := n^{-1} \sum_{i=1}^{n} K(i/(n+1))$.

Beyond its role in the derivation of the asymptotic distribution of the rank-based random vector (4.1), the following *asymptotic representation* result shows that $\underset{\sim}{\boldsymbol{\Delta}}{}^{(n)}_{f_1}(\boldsymbol{\vartheta})$ is indeed another version of the efficient central sequence $\boldsymbol{\Delta}^{\star(n)}_{f_1}(\boldsymbol{\vartheta})$.

LEMMA 4.1. *Assume that the score function $K : (0,1) \to \mathbb{R}$ is continuous, square integrable and that it can be expressed as the difference of two monotone increasing functions. Then, defining*

$$
\boldsymbol{\Delta}^{\star(n)}_{K;g_1}(\boldsymbol{\vartheta}) := \frac{1}{2} n^{-1/2} \mathbf{M}_k (\mathbf{V}^{\otimes 2})^{-1/2} \mathbf{J}_k^{\perp} \sum_{i=1}^{n} K\left(\tilde{G}_{1k}\left(\frac{d_i}{\sigma}\right)\right) \operatorname{vec}(\mathbf{U}_i \mathbf{U}_i'),
\tag{4.2}
$$

*we have* $\underset{\sim}{\mathbf{\Delta}}{}_K^{(n)}(\boldsymbol{\vartheta}) = \mathbf{\Delta}_{K;g_1}^{\star(n)}(\boldsymbol{\vartheta}) + o_{L^2}(1)$ *as* $n$ *goes to infinity, under* $\mathrm{P}_{\boldsymbol{\vartheta};g_1}^{(n)}$.

PROOF. See Appendix (Section A.3). □

4.2. *The proposed class of tests.* Let $K:(0,1)\to\mathbb{R}$ be some *score* function as in Lemma 4.1. Writing $\mathrm{E}[K(U)]$ and $\mathrm{E}[K^2(U)]$ for $\int_0^1 K(u)\,du$ and $\int_0^1 K^2(u)\,du$, respectively, the $K$-score version of the statistics we propose for testing $\mathcal{H}_0:\mathbf{V}=\mathbf{V}_0$ is

$$\underset{\sim}{Q}{}_K = \underset{\sim}{Q}{}_K^{(n)}$$

$$:= \frac{k(k+2)}{2n\mathrm{E}[K^2(U)]}\sum_{i,j=1}^n K\left(\frac{R_i}{n+1}\right)K\left(\frac{R_j}{n+1}\right)\left((\mathbf{U}_i'\mathbf{U}_j)^2-\frac{1}{k}\right), \tag{4.3}$$

where $R_i = R_i^{(n)}(\boldsymbol{\theta},\mathbf{V}_0)$ and $\mathbf{U}_i = \mathbf{U}_i^{(n)}(\boldsymbol{\theta},\mathbf{V}_0)$. Letting

$$\mathbf{S}_K = \mathbf{S}_K^{(n)} := \frac{1}{n}\sum_{i=1}^n K\left(\frac{R_i}{n+1}\right)\mathbf{U}_i\mathbf{U}_i',$$

these test statistics can be rewritten as

$$\underset{\sim}{Q}{}_K = \frac{nk(k+2)}{2\mathrm{E}[K^2(U)]}\left(\operatorname{tr}\mathbf{S}_K^2 - \frac{1}{k}\operatorname{tr}^2\mathbf{S}_K\right) \tag{4.4}$$

$$= \frac{k(k+2)\mathrm{E}^2[K(U)]}{k^2\mathrm{E}[K^2(U)]}\frac{nk^2}{2}\left\|\frac{\mathbf{S}_K}{\operatorname{tr}\mathbf{S}_K}-\frac{1}{k}\mathbf{I}_k\right\|^2 + o_\mathrm{P}(1)$$

$$= \frac{k(k+2)}{2\mathrm{E}[K^2(U)]}\left\|n^{-1/2}\sum_{i=1}^n K\left(\frac{R_i}{n+1}\right)\left(\mathbf{U}_i\mathbf{U}_i'-\frac{1}{k}\mathbf{I}_k\right)\right\|^2 + o_\mathrm{P}(1) \tag{4.5}$$

as $n$ goes to infinity, under any distribution [cf. (3.8)]. These test statistics are strictly invariant under $\mathcal{G}^{\mathrm{orth}(n)},\circ$ as well as under $\mathcal{G}^{(n)},\circ$. They admit (up to a multiplicative constant) an interesting interpretation as the sum of squared deviations of the eigenvalues of $\mathbf{S}_K$ from their arithmetic mean.

The power functions $K_a(u)=u^a$, $a\geq 0$, provide some traditional score functions. The corresponding test statistics are

$$\underset{\sim}{Q}{}_{K_a} := \frac{(2a+1)k(k+2)}{2n(n+1)^{2a}}\sum_{i,j=1}^n R_i^a R_j^a\left((\mathbf{U}_i'\mathbf{U}_j)^2-\frac{1}{k}\right). \tag{4.6}$$

Important particular cases are the sign-, Wilcoxon- and Spearman-type test statistics, defined by $\underset{\sim}{Q}{}_\mathrm{S} := \underset{\sim}{Q}{}_{K_0}$, $\underset{\sim}{Q}{}_\mathrm{W} := \underset{\sim}{Q}{}_{K_1}$ and $\underset{\sim}{Q}{}_\mathrm{SP} := \underset{\sim}{Q}{}_{K_2}$, respectively. In general, the resulting tests are not optimal at any density (they sometimes are, though—for instance, the Wilcoxon test $\underset{\sim}{Q}_W$ is optimal in dimension $k = 2$ at Student densities with two degrees of freedom; see Section 4.3), but they nevertheless yield good overall performances and are simple to compute. The sign test statistic $\underset{\sim}{Q}_S$ essentially coincides with that proposed by Ghosh and Sengupta [13] where, however, the $U_i'U_j$ are compared from Randles' interdirections (see [46]).

Local asymptotic optimality under radial density $f_1$ is achieved by the test based on $\underset{\sim}{Q}_{f_1} := \underset{\sim}{Q}_{K_{f_1}}$. This test statistic takes the form

$$\underset{\sim}{Q}_{f_1} = \frac{k(k+2)}{2n\mathcal{J}_k(f_1)} \sum_{i,j=1}^{n} K_{f_1}\left(\frac{R_i}{n+1}\right) K_{f_1}\left(\frac{R_j}{n+1}\right)\left((\mathbf{U}_i'\mathbf{U}_j)^2 - \frac{1}{k}\right) \tag{4.7}$$

which, letting $\mathbf{S}_{f_1} = \mathbf{S}_{f_1}^{(n)} := (1/n)\sum_{i=1}^{n} K_{f_1}(R_i/(n+1))\mathbf{U}_i\mathbf{U}_i'$, simplifies to

$$\begin{aligned}\underset{\sim}{Q}_{f_1} &= \frac{nk(k+2)}{2\mathcal{J}_k(f_1)}\left(\operatorname{tr}\mathbf{S}_{f_1}^2 - \frac{1}{k}\operatorname{tr}^2\mathbf{S}_{f_1}\right)\\ &= \frac{k(k+2)}{\mathcal{J}_k(f_1)}\frac{nk^2}{2}\left\|\frac{\mathbf{S}_{f_1}}{\operatorname{tr}\mathbf{S}_{f_1}} - \frac{1}{k}\mathbf{I}_k\right\|^2 + o_{\mathrm{P}}(1)\end{aligned} \tag{4.8}$$

as $n$ goes to infinity, still under any distribution. The van der Waerden (Gaussian scores $f_1 = \phi_1$) test, for instance, is based on the statistic

$$\underset{\sim}{Q}_{\mathrm{vdW}} := \frac{1}{2n}\sum_{i,j=1}^{n} \Psi_k^{-1}\left(\frac{R_i}{n+1}\right)\Psi_k^{-1}\left(\frac{R_j}{n+1}\right)\left((\mathbf{U}_i'\mathbf{U}_j)^2 - \frac{1}{k}\right), \tag{4.9}$$

where $\Psi_k$ stands for the chi-square distribution function with $k$ degrees of freedom. See (4.10) for the rank-based test statistics based on Student scores.

In order to describe the asymptotic behavior of $\underset{\sim}{Q}_K$ and $\underset{\sim}{Q}_{f_1}$, we will need the quantities

$$\mathcal{J}_k(K; g_1) := \int_0^1 K(u)K_{g_1}(u)\,du \quad \text{and} \quad \mathcal{J}_k(f_1, g_1) := \int_0^1 K_{f_1}(u)K_{g_1}(u)\,du$$

[$\mathcal{J}_k(f_1, g_1)$ can be interpreted as a measure of cross-information].

Denote by $\underset{\sim}{\phi}_K^{(n)}$ (resp. by $\underset{\sim}{\phi}_{f_1}^{(n)}$) the rank-based test which consists in rejecting $\mathcal{H}_0 : \mathbf{V} = \mathbf{V}_0$ as soon as $\underset{\sim}{Q}_K^{(n)}$, defined in (4.3) [resp. $\underset{\sim}{Q}_{f_1}^{(n)}$, defined in (4.7)] exceeds the $\alpha$-upper-quantile of a chi-square distribution with $k(k+1)/2 - 1$ degrees of freedom. We can now state the main result of this paper. Note that here the unions over $g_1$ extend over *all* possible standardized radial densities: contrary to the Gaussian tests described in Section 3.3, where finite fourth-order moments are required, the tests $\underset{\sim}{\phi}_K^{(n)}$ and $\underset{\sim}{\phi}_{f_1}^{(n)}$ are valid without any moment restrictions.

PROPOSITION 4.1. *Let $K$ be a continuous, square integrable score function defined on $(0,1)$ that can be expressed as the difference of two monotone increasing functions. Similarly, assume that $f_1$ [satisfying Assumptions (A1) and (A2)] is such that $K_{f_1}$ is continuous and can be expressed as the difference of two monotone increasing functions. Then:*

(i) *$\underset{\sim}{Q}_K^{(n)}$ and $\underset{\sim}{Q}_{f_1}^{(n)}$ are asymptotically chi-square with $k(k+1)/2-1$ degrees of freedom under $\bigcup_{\sigma^2}\bigcup_{g_1}\{\mathrm{P}^{(n)}_{\boldsymbol{\theta},\sigma^2,\mathbf{V}_0;g_1}\}$ and asymptotically noncentral chi-square, still with $k(k+1)/2-1$ degrees of freedom but with noncentrality parameters*

$$\frac{\mathcal{J}_k^2(K;g_1)}{2k(k+2)\mathrm{E}[K^2(U)]}\left[\mathrm{tr}((\mathbf{V}_0^{-1}\mathbf{v})^2)-\frac{1}{k}(\mathrm{tr}\,\mathbf{V}_0^{-1}\mathbf{v})^2\right]$$

*and*

$$\frac{\mathcal{J}_k^2(f_1,g_1)}{2k(k+2)\mathcal{J}_k(f_1)}\left[\mathrm{tr}((\mathbf{V}_0^{-1}\mathbf{v})^2)-\frac{1}{k}(\mathrm{tr}\,\mathbf{V}_0^{-1}\mathbf{v})^2\right],$$

*respectively, under $\bigcup_{\sigma^2}\{\mathrm{P}^{(n)}_{\boldsymbol{\theta},\sigma^2,\mathbf{V}_0+n^{-1/2}\mathbf{v};g_1}\}$;*

(ii) *the sequences of tests $\underset{\sim}{\phi}_K^{(n)}$ and $\underset{\sim}{\phi}_{f_1}^{(n)}$ have asymptotic level $\alpha$ under $\bigcup_{\sigma^2}\bigcup_{g_1}\{\mathrm{P}^{(n)}_{\boldsymbol{\theta},\sigma^2,\mathbf{V}_0;g_1}\}$;*

(iii) *the sequence of tests $\underset{\sim}{\phi}_{f_1}^{(n)}$ is locally and asymptotically maximin-efficient, still at asymptotic level $\alpha$, for $\bigcup_{\sigma^2}\bigcup_{g_1}\{\mathrm{P}^{(n)}_{\boldsymbol{\theta},\sigma^2,\mathbf{V}_0;g_1}\}$ against alternatives of the form $\bigcup_{\sigma^2}\bigcup_{\mathbf{V}\neq\mathbf{V}_0}\{\mathrm{P}^{(n)}_{\boldsymbol{\theta},\sigma^2,\mathbf{V};f_1}\}$.*

PROOF. See Appendix (Section A.3). □

Throughout the paper, our rank-based tests are described in terms of approximate critical values based on asymptotic chi-square null distributions. Of course, exact critical values could also be considered. These exact values can easily be simulated by sampling the $n!$ possible values of the vector of ranks and by independently generating uniformly distributed (over the unit sphere) signs.

4.3. *Asymptotic relative efficiencies.* Propositions 3.2 and 4.1 allow the computation of ARE values for $\underset{\sim}{\phi}_K^{(n)}$ (hence, for $\underset{\sim}{\phi}_{f_1}^{(n)}$) with respect to the adjusted John test $\phi_{\mathcal{N}}^{(n)}$ (therefore, also with respect to the adjusted Mauchly test) as ratios of the noncentrality parameters in the asymptotic distributions of their respective test statistics under local alternatives, for various radial densities $g_1$. These adjusted tests are still not valid unless $\kappa_k(g_1) < \infty$ and, therefore, our ARE values also require finite fourth-order moments. Recall, however, that the signed rank tests $\underset{\sim}{\phi}_K^{(n)}$ remain valid without such moment assumption so that, when $g_1$ is such that $\kappa_k(g_1) = \infty$, the asymptotic relative efficiency of any $\underset{\sim}{\phi}_K^{(n)}$ with respect to $\phi_{\mathcal{N}}^{(n)}$ can actually be considered as being infinite.

Proposition 4.2. *Let $K$ satisfy the assumptions of Proposition 4.1. Then the asymptotic relative efficiency of $\underset{\sim}{\phi}_K$ with respect to the parametric Gaussian test $\phi_{\mathcal{N}}$, under radial density $g_1$ satisfying Assumptions (A1), (A2) and $\kappa_k(g_1) < \infty$, is*

$$\mathrm{ARE}_{k,g_1}(\underset{\sim}{\phi}_K/\phi_{\mathcal{N}}) = \frac{1}{(k+2)^2}\frac{E_k(g_1)}{D_k^2(g_1)}\frac{\mathcal{J}_k^2(K;g_1)}{\mathrm{E}[K^2(U)]}.$$

*For $K$ of the form $K_{f_1}$, this yields*

$$\mathrm{ARE}_{k,g_1}(\underset{\sim}{\phi}_{f_1}/\phi_{\mathcal{N}}) = \frac{1}{(k+2)^2}\frac{E_k(g_1)}{D_k^2(g_1)}\frac{\mathcal{J}_k^2(f_1,g_1)}{\mathcal{J}_k(f_1)}.$$

In order to investigate the numerical values of these AREs, we consider the tests $\phi_{f_{1,\nu}^t}$ based on $t_\nu$-scores, that is, the scores associated with the Student radial densities introduced in Section 1.2. One can easily check that $\psi_{f_{1,\nu}^t}(r) = (k+\nu)a_{k,\nu}r/(\nu + a_{k,\nu}r^2)$. Also, since $a_{k,\nu}\|\mathbf{X}_1\|^2/k$, under $\mathrm{P}^{(n)}_{\mathbf{0},1,\mathbf{I}_k;f_{1,\nu}^t}$, is Fisher–Snedecor with $k$ and $\nu$ degrees of freedom, one can show that the test statistic $\underset{\sim}{Q}_{f_{1,\nu}^t}$ takes the form

$$\begin{aligned}\underset{\sim}{Q}_{f_{1,\nu}^t} &= \frac{k^2(k+\nu)(k+\nu+2)}{2n} \\ &\quad\times \sum_{i,j=1}^n \frac{T_i^{(n)}}{\nu + kT_i^{(n)}}\frac{T_j^{(n)}}{\nu + kT_j^{(n)}}\left((\mathbf{U}_i'\mathbf{U}_j)^2 - \frac{1}{k}\right)\end{aligned} \tag{4.10}$$

[see (2.2)], where, denoting by $G_{k,\nu}$ the Fisher–Snedecor distribution function with $k$ and $\nu$ degrees of freedom, we let $T_i^{(n)} := G_{k,\nu}^{-1}(R_i/(n+1))$. Note that the sign test and the van der Waerden test are obtained by letting $\nu \to 0$ and $\nu \to \infty$, respectively. An easy calculation also shows that for $\nu = 2$, $\underset{\sim}{Q}_{t_\nu}$ and $\underset{\sim}{Q}_{K_a}$ coincide for $a = 2/k$, $k = 2, 3, 4, \ldots$. Hence, for $k = 2$, the Wilcoxon test statistic $\underset{\sim}{Q}_W$ is optimal at Student densities with two degrees of freedom.

Numerical values of the AREs of several of the proposed rank-based tests with respect to the Gaussian test, under various $t_\nu$ and normal densities, are given in Table 2. For the sign test $\underset{\sim}{\phi}_S$, closed-form expressions are

$$\mathrm{ARE}_{k,f^t_{1,\nu}}[\underset{\sim}{\phi}_S/\phi_{\mathcal{N}}] = \frac{k(\nu-2)}{(k+2)(\nu-4)} \quad \text{and} \quad \mathrm{ARE}_{k,\phi_1}[\underset{\sim}{\phi}_S/\phi_{\mathcal{N}}] = \frac{k}{k+2}.$$

[recall that $\kappa_k(f^t_{1,\nu}) < \infty$ iff $\nu > 4$, which is the condition for a Student radial density to satisfy $E_k(f^t_{1,\nu}) < \infty$]. Also, the highest ARE with respect to the Gaussian test $\phi_{\mathcal{N}}$ that can be achieved under $t_\nu$ is

$$\mathrm{ARE}_{k,f^t_{1,\nu}}[\underset{\sim}{\phi}_{f^t_{1,\nu}}/\phi_{\mathcal{N}}] = \frac{(k+\nu)(\nu-2)}{(k+\nu+2)(\nu-4)}.$$

The ARE values in Table 2 are all uniformly good, especially for the van der Waerden test $\underset{\sim}{\phi}_{\mathrm{vdW}}$, for which they are not only uniformly larger than 1, but also uniformly larger than the corresponding AREs for location—namely, the AREs of van der Waerden rank tests with respect to the classical Hotelling ones when testing that the center of symmetry $\boldsymbol{\theta}$ of an elliptical distribution is equal to some fixed $\boldsymbol{\theta}_0$, as in [17]. This Pitman dominance of $\underset{\sim}{\phi}_{\mathrm{vdW}}$ over $\phi_{\mathcal{N}}$ also holds under lighter-than-Gaussian radial tails, as can be checked by again considering the power-exponential radial densities defined in Section 1.2; for instance, in the problem of testing for trivariate sphericity, the corresponding AREs are 1.166, 1.014, 1.000, 1.039, 1.108 and 1.183 for $\eta = 0.5$, 0.8, 1, 1.5, 2 and 2.5, respectively. Actually, it can be shown [43] that this is a general property and that $\underset{\sim}{\phi}_{\mathrm{vdW}}$, from the Pitman point of view, uniformly dominates its Gaussian parametric competitors.

4.4. *Unspecified location $\theta$.* In practice, the center of symmetry $\boldsymbol{\theta}$ is seldom specified and must be replaced, in test statistics, with an estimator $\hat{\boldsymbol{\theta}} = \hat{\boldsymbol{\theta}}^{(n)}$. Under very mild conditions, any root-$n$ consistent estimator will be adequate (in principle, after due discretization), but we recommend the (rotation-equivariant) spatial median (see, e.g., Möttönen and Oja [37]), which is itself "sign-based."

The asymptotic impact of this substitution on the validity of the signed rank tests proposed in Section 4.2 could be studied directly (see, e.g., [45]),

Table 2

*AREs of the $t_6$-, van der Waerden-, sign- and Wilcoxon-score rank-based tests for shape and (in parentheses) location, with respect to the corresponding parametric Gaussian tests, under k-dimensional Student (1, 3, 4, 5, 8, 15 and 20 degrees of freedom) and normal densities, respectively, for $k = 2, 3, 4, 6$ and 10*

| | | **Degrees of freedom of the underlying $t$ density** | | | | | | | |
|---|---|---|---|---|---|---|---|---|---|
| $\nu$ | $k$ | **1** | **3** | **4** | **5** | **8** | **15** | **20** | $\infty$ |
| $\underset{\sim}{\phi}_{t_6}$ | 2 | $+\infty$ | $+\infty$ | $+\infty$ | 2.331 | 1.248 | 1.045 | 1.013 | 0.957 |
| | | $(+\infty)$ | (2.067) | (1.484) | (1.294) | (1.107) | (1.009) | (0.986) | (0.927) |
| | 3 | $+\infty$ | $+\infty$ | $+\infty$ | 2.398 | 1.267 | 1.052 | 1.018 | 0.957 |
| | | $(+\infty)$ | (2.174) | (1.540) | (1.331) | (1.124) | (1.014) | (0.988) | (0.919) |
| | 4 | $+\infty$ | $+\infty$ | $+\infty$ | 2.453 | 1.284 | 1.058 | 1.023 | 0.958 |
| | | $(+\infty)$ | (2.258) | (1.584) | (1.361) | (1.139) | (1.019) | (0.990) | (0.913) |
| | 6 | $+\infty$ | $+\infty$ | $+\infty$ | 2.537 | 1.311 | 1.070 | 1.031 | 0.959 |
| | | $(+\infty)$ | (2.382) | (1.652) | (1.408) | (1.163) | (1.028) | (0.995) | (0.905) |
| | 10 | $+\infty$ | $+\infty$ | $+\infty$ | 2.646 | 1.349 | 1.087 | 1.044 | 0.963 |
| | | $(+\infty)$ | (2.534) | (1.736) | (1.468) | (1.196) | (1.043) | (1.005) | (0.896) |
| $\underset{\sim}{\phi}_{\mathrm{vdW}}$ | 2 | $+\infty$ | $+\infty$ | $+\infty$ | 2.204 | 1.215 | 1.047 | 1.025 | 1.000 |
| | | $(+\infty)$ | (1.729) | (1.301) | (1.171) | (1.060) | (1.016) | (1.009) | (1.000) |
| | 3 | $+\infty$ | $+\infty$ | $+\infty$ | 2.270 | 1.233 | 1.052 | 1.028 | 1.000 |
| | | $(+\infty)$ | (1.798) | (1.336) | (1.194) | (1.069) | (1.019) | (1.011) | (1.000) |
| | 4 | $+\infty$ | $+\infty$ | $+\infty$ | 2.326 | 1.249 | 1.057 | 1.031 | 1.000 |
| | | $(+\infty)$ | (1.853) | (1.364) | (1.212) | (1.077) | (1.022) | (1.012) | (1.000) |
| | 6 | $+\infty$ | $+\infty$ | $+\infty$ | 2.413 | 1.275 | 1.066 | 1.036 | 1.000 |
| | | $(+\infty)$ | (1.935) | (1.408) | (1.242) | (1.092) | (1.027) | (1.016) | (1.000) |
| | 10 | $+\infty$ | $+\infty$ | $+\infty$ | 2.531 | 1.312 | 1.080 | 1.045 | 1.000 |
| | | $(+\infty)$ | (2.041) | (1.467) | (1.283) | (1.112) | (1.035) | (1.021) | (1.000) |
| $\underset{\sim}{\phi}_{S}$ | 2 | $+\infty$ | $+\infty$ | $+\infty$ | 1.500 | 0.750 | 0.591 | 0.563 | 0.500 |
| | | $(+\infty)$ | (2.000) | (1.388) | (1.185) | (0.984) | (0.877) | (0.851) | (0.785) |
| | 3 | $+\infty$ | $+\infty$ | $+\infty$ | 1.800 | 0.900 | 0.709 | 0.675 | 0.600 |
| | | $(+\infty)$ | (2.162) | (1.500) | (1.281) | (1.063) | (0.947) | (0.920) | (0.849) |
| | 4 | $+\infty$ | $+\infty$ | $+\infty$ | 2.000 | (1.000 | 0.788 | 0.750 | 0.667 |
| | | $(+\infty)$ | (2.250) | (1.561) | (1.333) | (1.107) | (0.986) | (0.958) | (0.884) |
| | 6 | $+\infty$ | $+\infty$ | $+\infty$ | 2.250 | (1.125 | 0.886 | 0.844 | 0.750 |
| | | $(+\infty)$ | (2.344) | (1.626) | (1.389) | (1.153) | (1.027) | (0.997) | (0.920) |
| | 10 | $+\infty$ | $+\infty$ | $+\infty$ | 2.500 | 1.250 | 0.985 | 0.938 | 0.833 |
| | | $(+\infty)$ | (2.422) | (1.681) | (1.436) | (1.192) | (1.062) | (1.031) | (0.951) |
| $\underset{\sim}{\phi}_{W}$ | 2 | $+\infty$ | $+\infty$ | $+\infty$ | 2.258 | 1.174 | 0.956 | 0.919 | 0.844 |
| | | $(+\infty)$ | (1.748) | (1.317) | (1.185) | (1.066) | (1.015) | (1.005) | (0.985) |
| | 3 | $+\infty$ | $+\infty$ | $+\infty$ | 2.386 | 1.246 | 1.022 | 0.985 | 0.913 |
| | | $(+\infty)$ | (1.621) | (1.233) | (1.117) | (1.019) | (0.983) | (0.978) | (0.975) |
| | 4 | $+\infty$ | $+\infty$ | $+\infty$ | 2.432 | 1.273 | 1.048 | 1.012 | 0.945 |
| | | $(+\infty)$ | (1.533) | (1.171) | (1.064) | (0.979) | (0.954) | (0.952) | (0.961) |
| | 6 | $+\infty$ | $+\infty$ | $+\infty$ | 2.451 | 1.283 | 1.060 | 1.026 | 0.969 |
| | | $(+\infty)$ | (1.422) | (1.090) | (0.994) | (0.921) | (0.908) | (0.911) | (0.938) |
| | 10 | $+\infty$ | $+\infty$ | $+\infty$ | 2.426 | 1.264 | 1.045 | 1.013 | 0.970 |
| | | $(+\infty)$ | (1.315) | (1.007) | (0.919) | (0.855) | (0.851) | (0.857) | (0.907) |

but is more conveniently handled via Le Cam's third lemma, which allows the derivation of the asymptotic distribution under $\mathrm{P}^{(n)}_{\boldsymbol{\theta},\sigma^2,\mathbf{V};g_1}$ of the test statistic $\underset{\sim}{Q}^{(n)}_{K} =: \underset{\sim}{Q}^{(n)}_{K;\boldsymbol{\theta}}$ considered in Section 4.2, but computed at $\hat{\boldsymbol{\theta}}$ instead of $\boldsymbol{\theta}$. This lemma applies in the parametric location experiment $\mathcal{E}^{(n)}_g := \{\mathrm{P}^{(n)}_{\boldsymbol{\theta},\sigma^2,\mathbf{V};g_1} | \boldsymbol{\theta} \in \mathbb{R}^k\}$, provided that it is ULAN, which essentially requires $g_1$ to satisfy Assumption (A1) (see [17]).

The asymptotic distribution, as $n \to \infty$, of $\underset{\sim}{Q}^{(n)}_{K;\boldsymbol{\theta}+n^{-1/2}\boldsymbol{\tau}^{(n)}}$ under $\mathrm{P}^{(n)}_{\boldsymbol{\theta},\sigma^2,\mathbf{V};g_1}$ for any bounded sequence $\boldsymbol{\tau}^{(n)}$ is the same as under $\mathrm{P}^{(n)}_{\boldsymbol{\theta}+n^{-1/2}\boldsymbol{\tau}^{(n)},\sigma^2,\mathbf{V};g_1}$ [viz., in view of part (i) of Proposition 4.1, chi-square with $k(k+1)/2-1$ degrees of freedom], provided that the asymptotic joint distribution, under $\mathrm{P}^{(n)}_{\boldsymbol{\theta},\sigma^2,\mathbf{V};g_1}$, of $\boldsymbol{\Delta}^{(n)}_{K;g_1}(\boldsymbol{\vartheta})$ [defined in (4.2)] and the central sequence for location $\boldsymbol{\Delta}^{(n)}_{g_1;1}(\boldsymbol{\vartheta})$ in $\mathcal{E}^{(n)}_g$ [as defined in (2.5)] is normal with block-diagonal asymptotic covariance. Now, this is automatically satisfied under the assumptions made on $K$: indeed, both $\boldsymbol{\Delta}^{\star(n)}_{K;g_1}(\boldsymbol{\vartheta})$ and $\boldsymbol{\Delta}^{(n)}_{g_1;1}(\boldsymbol{\vartheta})$ are sums of i.i.d. vectors with finite variances and, in view of the independence under $\mathrm{P}^{(n)}_{\boldsymbol{\theta},\sigma^2,\mathbf{V};g_1}$ between $d^{(n)}_i(\boldsymbol{\theta},\mathbf{V})$ and $\mathbf{U}^{(n)}_i(\boldsymbol{\theta},\mathbf{V})$, have a cross-covariance matrix proportional to $\mathrm{E}[\mathrm{vec}(\mathbf{U}_i\mathbf{U}_i')\mathbf{U}_i'] = \mathbf{0}$. Classical reasoning then extends this to random sequences of the form $\boldsymbol{\tau}^{(n)} = n^{1/2}(\hat{\boldsymbol{\theta}} - \boldsymbol{\theta})$, where $n^{1/2}(\hat{\boldsymbol{\theta}} - \boldsymbol{\theta})$ is $O_{\mathrm{P}}(1)$ and $\hat{\boldsymbol{\theta}}$ is *locally discrete*, that is, such that the number, under $\mathrm{P}^{(n)}_{\boldsymbol{\theta},\sigma^2,\mathbf{V};g_1}$, of its possible values in balls of the form $\{\mathbf{z} \in \mathbb{R}^k | \|\mathbf{z} - \boldsymbol{\theta}\|^2 \le b^2\}$ remains bounded as $n \to \infty$. It is well known that this latter assumption has no practical consequences (see, e.g., [30]). The null distribution of $\underset{\sim}{Q}^{(n)}_{K;\hat{\boldsymbol{\theta}}}$ is thus the same, then, as that of $\underset{\sim}{Q}^{(n)}_{K;\boldsymbol{\theta}}$.

However, Le Cam's third lemma only provides asymptotic equivalence in distribution results. Asymptotic equivalence in probability [i.e., a result of the form $\underset{\sim}{Q}^{(n)}_{K;\hat{\boldsymbol{\theta}}} - \underset{\sim}{Q}^{(n)}_{K;\boldsymbol{\theta}} = o_{\mathrm{P}}(1)$] under $\mathrm{P}^{(n)}_{\boldsymbol{\theta},\sigma^2,\mathbf{V};g_1}$ requires more stringent *asymptotic linearity* results, such as those in Proposition A.1 of [16], or more general methods, such as the one recently developed by Andreou and Werker [2].

Note that $\underset{\sim}{Q}^{(n)}_{K;\hat{\boldsymbol{\theta}}}$ is no longer strictly invariant or distribution-free, but remains asymptotically so, in the sense of being asymptotically equivalent to its genuinely invariant and distribution-free counterpart $\underset{\sim}{Q}^{(n)}_{K;\boldsymbol{\theta}}$. This asymptotic equivalence carries over to contiguous alternatives so that local optimal-ity properties are also preserved. Incidentally, note that $\underset{\sim}{Q}{}^{(n)}_{K;\hat{\boldsymbol\theta}}$ is translation-invariant whenever $\hat{\boldsymbol\theta}$ is translation-equivariant.

## 5. Validity and consistency properties.

5.1. *Null hypothesis*: *sphericity or unit shape*? Our rank tests are basically intended for the null hypothesis of sphericity—not for the hypothesis of isotropy, nor for that of unit shape. Indeed the (asymptotic) size of $\underset{\sim}{\phi}_K$ does not, in general, match the nominal $\alpha$-level under nonelliptical densities, even for unit shape matrices $\mathbf{V} = \mathbf{I}_k$.

One important exception to this general rule is the multivariate sign test $\underset{\sim}{\phi}_S$, based on the test statistic [with scores $K(u) = 1$] $\underset{\sim}{Q}_S := \underset{\sim}{Q}_{K_0}$ given in (4.6). This test in [13] is described as a test of sphericity. However, since the ranks are not involved, $\underset{\sim}{\phi}_S$ remains valid under the hypothesis of isotropy and hence (since only the centering and second-order structure of the matrices $\mathbf{U}_i\mathbf{U}_i'$ matter) under the hypothesis of unit shape with *isotropic fourth-order moments*, that is, provided that the moments of the signs $\mathbf{U}_i$ coincide with those of the uniform distribution over the unit sphere in $\mathbb{R}^k$ up to order four, so that

$$\mathrm{E}[\mathbf{U}_i\mathbf{U}_i'] = \frac{1}{k}\mathbf{I}_k$$

and

$$\mathrm{E}[\mathrm{vec}(\mathbf{U}_i\mathbf{U}_i')(\mathrm{vec}(\mathbf{U}_i\mathbf{U}_i'))'] = \frac{1}{k(k+2)}[\mathbf{I}_{k^2} + \mathbf{K}_k + \mathbf{J}_k].$$

The validity of this test can be extended to the whole hypothesis of unit shape if estimated moments of order four are substituted for the isotropic ones, yielding the *adjusted* sign test $\underset{\sim}{\phi}{}^*_S$, based on the statistic

$$\begin{aligned}\underset{\sim}{Q}{}^*_S &:= \left(\sum_{i=1}^n \mathrm{vec}(\mathbf{U}_i\mathbf{U}_i') - \frac{n}{k}\mathrm{vec}(\mathbf{I}_k)\right)'(\mathbf{V}^{\otimes 2})^{-1/2}\mathbf{M}_k' \\ &\quad\times\Bigg[\mathbf{M}_k(\mathbf{V}^{\otimes 2})^{-1/2} \\ &\quad\quad\times\sum_{i=1}^n\left((\mathrm{vec}(\mathbf{U}_i\mathbf{U}_i'))(\mathrm{vec}(\mathbf{U}_i\mathbf{U}_i'))' - \frac{n}{k^2}\mathbf{J}_k\right)(\mathbf{V}^{\otimes 2})^{-1/2}\mathbf{M}_k'\Bigg]^{-1} \\ &\quad\times\mathbf{M}_k(\mathbf{V}^{\otimes 2})^{-1/2}\left(\sum_{i=1}^n \mathrm{vec}(\mathbf{U}_i\mathbf{U}_i') - \frac{n}{k}\mathrm{vec}(\mathbf{I}_k)\right).\end{aligned} \tag{5.1}$$

Unfortunately, the benefits of Lemma 3.1 are lost and the adjusted test statistic $\underset{\sim}{Q}{}^*_S$ does not retain the elegant and simple structure [cf. (1.4) and (1.5)] of John's test.

One is tempted to apply a similar idea to our rank-based tests $\underset{\sim}{\phi}_K$. An estimate of the covariance matrix of $\underset{\sim}{\mathbf{\Delta}}{}^{(n)}_K(\boldsymbol{\vartheta})$ that does not exploit the elliptical independence between the ranks and the signs is indeed quite possible. But the expectation of $\sum_{i=1}^n K(R_i/(n+1))\,\mathrm{vec}(\mathbf{U}_i\mathbf{U}_i')$ [reducing to $k^{-1}\sum_{i=1}^n K(i/(n+1))\,\mathrm{vec}(\mathbf{I}_k)$ under sphericity] is no longer distribution-free if the assumption of ellipticity is abandoned, and replacing this expectation with an empirical centering would induce a noncentrality parameter in the asymptotic null distribution of the test statistic $\underset{\sim}{Q}_K$.

From the point of view of (asymptotic) validity and with the exception of the multivariate sign test [the adjusted version (5.1)] which is a test of unit shape, our rank-based tests $\underset{\sim}{\phi}_K$ thus only qualify for the null hypothesis of sphericity.

5.2. *Nonlocal alternatives*: *consistency issues.* Validity under the null hypothesis not the only requirement for a test $\phi$ to qualify as a test of $\mathcal{H}_0$, say, against $\mathcal{H}_1$, and consistency under $\mathcal{H}_1$ is certainly an equally important issue. In this respect, the larger the overarching model $\mathcal{H} := \mathcal{H}_0 + \mathcal{H}_1$ (with $+$ standing for disjoint union), the better the test. Although optimality results have been derived under an overall hypothesis $\mathcal{H}$ of ellipticity, the most "natural" $\mathcal{H}$ here should consist of the collection of all i.i.d. sample distributions from nonvanishing $k$-dimensional densities $\underline{f}$.

However, the results of the previous sections are entirely local to the null hypothesis of sphericity and do not allow for any conclusions under nonlocal alternatives. Proposition 5.1 below, on the other hand, provides a characterization of consistency under nonlocal alternatives. Denote by $\mathcal{K}^{(n)}(\underline{f})$ the hypothesis under which the observations $\mathbf{X}_i$ are i.i.d. with nonvanishing, possibly nonelliptic density $\underline{f}$. Our rank tests $\underset{\sim}{\phi}{}^{(n)}_K$ are consistent under $\mathcal{K}^{(n)}(\underline{f})$ iff the quadratic test statistic $\underset{\sim}{Q}{}^{(n)}_K$ is unbounded in probability, that is, iff, for any fixed $q$, $\mathrm{P}[\underset{\sim}{Q}{}^{(n)}_K > q] \to 1$ as $n \to \infty$, under $\mathcal{K}^{(n)}(\underline{f})$ or, equivalently, iff, for all $t \geq 0$,

$$\mathrm{P}\left[n^{-1/2}\left\|\sum_{i=1}^n K\left(\frac{R_i}{n+1}\right)\mathrm{vec}\left(\mathbf{U}_i\mathbf{U}_i' - \frac{1}{k}\mathbf{I}_k\right)\right\| > t\right] \longrightarrow 1 \tag{5.2}$$

as $n \to \infty$, under $\mathcal{K}^{(n)}(\underline{f})$ [see (4.5)], which we unambiguously write as $n^{-1/2} \sum_{i=1}^{n} K((R_i/(n+1))) \operatorname{vec}(\mathbf{U}_i\mathbf{U}_i' - \frac{1}{k}\mathbf{I}_k) \xrightarrow{\mathrm{P}} \infty$ as $n \to \infty$. We then have the following necessary and/or sufficient consistency conditions:

Proposition 5.1. *Assume that the score function $K:(0,1) \to \mathbb{R}$ can be expressed as the difference $K_1 - K_2$ of two monotone increasing, absolutely continuous and square integrable functions. Then:*

(i) *$\underset{\sim}{\phi}_K^{(n)}$ is consistent iff, under $\mathcal{K}^{(n)}(\underline{f})$, as $n \to \infty$,*

$$n^{-1/2} \sum_{i=1}^{n} \mathrm{E}\left[K\left(\frac{R_i^{(n)}}{n+1}\right)\Big|\mathbf{U}^{(n)}\right] \operatorname{vec}\left(\mathbf{U}_i\mathbf{U}_i' - \frac{1}{k}\mathbf{I}_k\right) \xrightarrow{\mathrm{P}} \infty, \tag{5.3}$$

*where $R_i^{(n)} = R_i^{(n)}(\boldsymbol{\theta}, \mathbf{I}_{\boldsymbol{k}})$, $\mathbf{U}_i = \mathbf{U}_i(\boldsymbol{\theta}, \mathbf{I}_{\boldsymbol{k}})$ and $\mathbf{U}^{(n)} := (\mathbf{U}_1, \ldots, \mathbf{U}_n)$.*

(ii) *If the square integrability condition on $K_1$ and $K_2$ in* (i) *is reinforced into*

$$J(K_i) := \int_0^1 u^{1/2}(1-u)^{1/2}\, dK_i(u) < \infty, \qquad i = 1, 2 \tag{5.4}$$

*(a classical condition that goes back to Hoeffding* [22]*), then $\underset{\sim}{\phi}_K^{(n)}$ is consistent iff, under $\mathcal{K}^{(n)}(\underline{f})$, as $n \to \infty$,*

$$\begin{aligned} &n^{-1/2} \sum_{i=1}^{n} \mathrm{E}\left[K\left(\frac{1}{n}\sum_{j=1}^{n} \mathrm{P}[d_j \le d_i | d_i, \mathbf{U}_i, \mathbf{U}_j]\right)\Big|\mathbf{U}^{(n)}\right] \\ &\qquad \times \operatorname{vec}\left(\mathbf{U}_i\mathbf{U}_i' - \frac{1}{k}\mathbf{I}_k\right) \xrightarrow{\mathrm{P}} \infty. \end{aligned} \tag{5.5}$$

(iii) *If the Hoeffding condition* (5.4) *is satisfied and, moreover, $K$ is convex, then a sufficient condition for $\underset{\sim}{\phi}_K^{(n)}$ to be consistent is that, for some $\ell$, either*

$$\begin{aligned} &K\left(\frac{\mathrm{E}[I[d_2 \le d_1] \operatorname{vec}(\mathbf{U}_1\mathbf{U}_1' - (1/k)\mathbf{I}_k)_\ell^-]}{\mathrm{E}[\operatorname{vec}(\mathbf{U}_1\mathbf{U}_1' - (1/k)\mathbf{I}_k)_\ell^-]}\right) \\ &\qquad - \frac{\mathrm{E}[K(\mathrm{P}[d_2 \le d_1 | d_1, \mathbf{U}_1, \mathbf{U}_2]) \operatorname{vec}(\mathbf{U}_1\mathbf{U}_1' - (1/k)\mathbf{I}_k)_\ell^+]}{\mathrm{E}[\operatorname{vec}(\mathbf{U}_1\mathbf{U}_1' - (1/k)\mathbf{I}_k)_\ell^-]} > 0 \end{aligned} \tag{5.6}$$

*or*

$$\begin{aligned} &K\left(\frac{\mathrm{E}[I[d_2 \le d_1] \operatorname{vec}(\mathbf{U}_1\mathbf{U}_1' - (1/k)\mathbf{I}_k)_\ell^+]}{\mathrm{E}[\operatorname{vec}(\mathbf{U}_1\mathbf{U}_1' - (1/k)\mathbf{I}_k)_\ell^+]}\right) \\ &\qquad - \frac{\mathrm{E}[K(\mathrm{P}[d_2 \le d_1 | d_1, \mathbf{U}_1, \mathbf{U}_2]) \operatorname{vec}(\mathbf{U}_1\mathbf{U}_1' - (1/k)\mathbf{I}_k)_\ell^-]}{\mathrm{E}[\operatorname{vec}(\mathbf{U}_1\mathbf{U}_1' - (1/k)\mathbf{I}_k)_\ell^+]} > 0 \end{aligned} \tag{5.7}$$

*under* $\mathcal{K}^{(n)}(\underline{f})$, *where* $\mathrm{vec}(\mathbf{U}_1\mathbf{U}_1' - \frac{1}{k}\mathbf{I}_k)_\ell^{\pm}$ *stand for the positive and negative parts of* $\mathrm{vec}(\mathbf{U}_1\mathbf{U}_1' - \frac{1}{k}\mathbf{I}_k)_\ell$, *respectively.*

(iv) *The Wilcoxon test* $\underset{\sim}{\phi}_W$ *based on*

$$\underset{\sim}{Q}_W := \frac{3k(k+2)}{2n(n+1)^2}\sum_{i=1}^{n} R_i^{(n)} R_j^{(n)}\left((\mathbf{U}_i'\mathbf{U}_j)^2 - \frac{1}{k}\right)$$

[*see* (4.6)] *is consistent iff, under* $\mathcal{K}^{(n)}(\underline{f})$,

$$\mathrm{E}\left[I[d_2 \le d_1]\,\mathrm{vec}\left(\mathbf{U}_1\mathbf{U}_1' - \frac{1}{k}\mathbf{I}_k\right)\right] \neq \mathbf{0}. \tag{5.8}$$

(v) *The adjusted sign test* $\underset{\sim}{\phi}_S^*$ *based on* (5.1) *is consistent iff, under* $\mathcal{K}^{(n)}(\underline{f})$,

$$\mathrm{E}\left[\mathrm{vec}\left(\mathbf{U}_1\mathbf{U}_1' - \frac{1}{k}\mathbf{I}_k\right)\right] \neq \mathbf{0}. \tag{5.9}$$

PROOF. See Appendix (Section A.4). □

Note that the Hoeffding condition in (ii) only slightly reinforces the square integrability condition on $K$: Hoeffding [22] shows that (5.4) holds as soon as $\int_0^1 (K(u))^2[\log(1+|K(u)|)]^{1+\delta}\,du$ is finite for some $\delta > 0$, a condition that is satisfied by all particular score functions considered in this paper.

These consistency results imply that our rank-based tests $\underset{\sim}{\phi}_K^{(n)}$ (excluding the sign test), although unrestrictedly valid under the null hypothesis of sphericity, are consistent under most nonspherical alternatives, which include nonspherical elliptic, nonelliptical unit shape and nonunit-shape cases. For Wilcoxon scores, for instance, only the very particular densities $\underline{f}$ for which $I[d_2 \le d_1]$ is orthogonal to the $k(k+1)/2$ variables $U_{1,r}U_{1,s} - (\delta_{rs}/k)$, $r, s = 1, \ldots, k$ ($\delta_{rs}$ standing for the Kronecker symbol) result in an inconsistent $\underset{\sim}{\phi}_W^{(n)}$. This either corresponds to $\mathrm{E}[\mathbf{U}_1\mathbf{U}_1'] \neq \mathbf{I}_k/k$ and the joint distribution of $d_1, \mathbf{U}_1$ and $d_2$ compensating exactly for the deviations of all $U_{1,r}U_{1,s}$'s from $\delta_{rs}/k$ ($r, s = 1, \ldots, k$), or to unit shape densities under which $I[d_2 \le d_1]$ and $\mathbf{U}_1\mathbf{U}_1'$ are uncorrelated. To the best of our knowledge, the only test retaining consistency under the whole nonspherical alternative is Baringhaus' test [5]; but the price that must be paid is that the separation rates are nonparametric, which entails that its ARE with respect to $\underset{\sim}{\phi}_K^{(n)}$ is zero at elliptical alternatives.

The situation is slightly different with the adjusted sign test. As already mentioned, the natural null hypothesis for this test is that of unit shape and consistency is achieved at all nonunit shape alternatives, since the score $[K(u)=1]$ cannot here compensate for deviations from $\mathbf{I}_k/k$ of $\mathbf{U}_i\mathbf{U}_i'$. On the other hand, the price to be paid in terms of efficiency at elliptical alternatives can be quite high: the AREs of sign tests with respect to their van der Waerden counterparts $\underset{\sim}{\phi}{}^{(n)}_{\mathrm{vdW}}$ are only 0.681, 0.500 and 0.279, respectively, at $t_5$, Gaussian and $e_3$ alternatives, in dimension $k=2$.

As the dimension $k$ of the observation space goes to $\infty$, however, it can easily be shown that, for fixed $n$, $\underset{\sim}{\phi}{}^{(n)}_{\mathrm{vdW}} - \underset{\sim}{\phi}{}^{(n)}_{S} = o(1)$ $\mathcal{P}^{(n)}$-a.s.; this justifies the empirical finding that the AREs of the sign test with respect to the van der Waerden test converge to 1, as $k\to\infty$, irrespective of the underlying distribution. Most interestingly, this convergence also implies that the van der Waerden test in some sense inherits, as $k\to\infty$, most of the nice validity/consistency properties of the sign test, whereas the latter, on the other hand, inherits the attractive efficiency properties of van der Waerden procedures.

**6. Simulation results.** The asymptotic relative efficiencies of the tests (of the null hypothesis $\mathbf{V}=\mathbf{V}_0$) described in Sections 3.3 and 4.2 do not depend on the null value $\mathbf{V}_0$ of the shape matrix. Therefore, in this section, we concentrate on the particular case ($\mathbf{V}_0=\mathbf{I}_k$) of testing for sphericity. We generated $N=2{,}500$ independent samples $\boldsymbol{\varepsilon}_1,\dots,\boldsymbol{\varepsilon}_{500}$ of size $n=500$ from various bivariate spherical densities (the bivariate normal and bivariate $t$-distributions with 0.2, 1 and 6 degrees of freedom, resp.), with center of symmetry $\boldsymbol{\theta}=(0,0)'$. From each of these samples, we constructed four series of 500 spherical (for $m=0$) or elliptical (for $m=1,2,3$) observations $\mathbf{X}_1,\dots,\mathbf{X}_{500}$, characterized by

$$\mathbf{X}_i = (\mathbf{I}_k + m\mathbf{v})\boldsymbol{\varepsilon}_i, \qquad m=0,1,2,3, \tag{6.1}$$

with $\overset{\circ}{\mathrm{vec}}\,\mathrm{h}\,\mathbf{v}=(0,.14)'$.

Although designed against elliptical alternatives, our tests also perform quite well under a broad class of nonelliptical alternatives. In order to show this, we considered the following skew populations. Population $\mathcal{SN}$ refers to samples of $n=500$ observations $\mathbf{X}_1,\dots,\mathbf{X}_{500}$ characterized by

$$\mathbf{X}_i = (\mathrm{sign} V_{m;i})\mathbf{W}_{m;i} - \mathrm{E}[(\mathrm{sign} V_{m;i})\mathbf{W}_{m;i}], \qquad m=0,1,2,3, \tag{6.2}$$

where the i.i.d. vectors $(V_{m;i},\mathbf{W}_{m;i}')'$ are drawn from the trivariate standard normal distribution with mean $\mathbf{0}$ and covariance matrix

$$\begin{pmatrix} 1 & \boldsymbol{\delta}' \\ \boldsymbol{\delta} & \mathbf{I}_2 \end{pmatrix}, \qquad \boldsymbol{\delta}=(1+m^2\mathbf{v}'\mathbf{v})^{-1/2}m\mathbf{v},$$

with $\mathbf{v}=(0.15,0)'$. The distribution of the resulting $\mathbf{X}_i$'s is the so-called *bivariate skew normal distribution* with parameters $\mathbf{0}$, $\mathbf{I}_2$ and $m\mathbf{v}$ (see, e.g., [3] or [4]). Population $\mathcal{S}t_2$ is obtained in the same way, but with trivariate $t_2$-distributed vectors $(V_{m;i}, \mathbf{W}'_{m;i})'$ with the same mean and covariance matrix as in the Gaussian case above, but $\mathbf{v} = (0.25, 0)'$ (see [4]).

On each of these samples, we performed the following eleven tests for sphericity (all at asymptotic level $\alpha = 5\%$): John's test [based on (3.9)], the Gaussian test $\phi_{\mathcal{N}}$ [based on (3.7)], the sign, Wilcoxon and Spearman tests [based on $\underset{\sim}{Q}_{K_0}$, $\underset{\sim}{Q}_{K_1}$ and $\underset{\sim}{Q}_{K_2}$ in (4.6), resp.], the van der Waerden test $\underset{\sim}{\phi}_{\mathrm{vdW}}$ [based on (4.9)], and several $t_\nu$-score tests $\underset{\sim}{\phi}_{f^t_{1,\nu}}$ ($\nu = 0.2, 0.5, 1, 2$ and 6) [based on (4.10)]. Rejection frequencies are reported in Table 3. The corresponding individual confidence intervals (for $N = 2{,}500$ replications) at confidence level 0.95 have half-widths 0.0044, 0.0080 and 0.0100, for frequencies of the order of 0.05 (0.95), 0.20 (0.80) and 0.50, respectively.

Inspection of Table 3 reveals that the Gaussian test $\phi_{\mathcal{N}}$ collapses under the heavy-tailed distributions $t_{0.2}$ and $t_1$ (which have infinite fourth-order moments) and confirms the fact that John's test is only valid under normal distributions. All rank-based tests apparently satisfy the 5% probability level constraint. Power rankings are essentially consistent with the corresponding ARE values, which we also report in Table 3. In particular, the asymptotic optimality of $\underset{\sim}{\phi}_{f^t_{1,\nu}}$ under the Student distribution with $\nu$ degrees of freedom is confirmed. The performances under elliptical and nonelliptical alternatives of the various procedures seem to be quite similar.

Finally, in order to investigate the performances of our tests in very small samples, we generated $N = 2{,}500$ independent samples of size $n = 25$ based on (6.1) [but with $\overset{\circ}{\mathrm{vech}}\,\mathbf{v} = (0, 0.2)'$]. Only Gaussian and $t_{0.2}$ densities were considered. The corresponding rejection frequencies are reported in Table 4. Similar conclusions as in the first Monte Carlo study above hold in this small sample simulation. However, note that for such a small sample size, the asymptotic approximation seems to produce strictly conservative critical values for the van der Waerden- and $t_6$-score versions of our tests.

## APPENDIX

**A.1. Proof of Proposition 2.1.** Our proof relies on Lemma 1 from Swensen [49] (more precisely, on its extension by Garel and Hallin [12]). The sufficient conditions for LAN given in Swensen's result follow from standard arguments once it is shown that $(\boldsymbol{\theta}, \sigma^2, \mathbf{V}) \mapsto \underline{f}^{1/2}_{\boldsymbol{\theta},\sigma^2,\mathbf{V};f_1}(\mathbf{x})$ is differentiable in quadratic mean, where $\underline{f}_{\boldsymbol{\theta},\sigma^2,\mathbf{V};f_1}$ is the density in (1.1), and we therefore focus on this. The main step in establishing this quadratic mean differentiability is the following [here and in the sequel, all $o(\|\cdot\|)$ or $O(\|\cdot\|)$ quantities are taken as $\|\cdot\| \to 0$]:

TABLE 3

*Rejection frequencies (out of $N = 2{,}500$ replications), under various null and nonnull distributions [see (6.1) and (6.2) for details], of John's test ($\phi_{\text{John}}$), the Gaussian parametric test ($\phi_{\mathcal{N}}$) and the signed-rank van der Waerden ($\underset{\sim}{\phi}_{\text{vdW}}$), $t_\nu$-score ($\underset{\sim}{\phi}_{f_{1,\nu}}$, $\nu = 0.2, 0.5, 1, 2, 6$), sign($\underset{\sim}{\phi}_S$), Wilcoxon-type ($\underset{\sim}{\phi}_W$) and Spearman-type ($\underset{\sim}{\phi}_{SP}$) tests, respectively; the sample size is 500 ("ND" means "not defined," which occurs as soon as one of the two tests involved is not valid under the distribution being considered; "?" indicates that no theoretical ARE values are available under nonelliptical alternatives)*

| | | *m* | | | | |
|---|---|---|---|---|---|---|
| **Test** | | **0** | **1** | **2** | **3** | **ARE** |
| $\phi_{\text{John}}$ | $\mathcal{N}$ | 0.0504 | 0.2380 | 0.6856 | 0.9492 | 1.000 |
| $\phi_{\mathcal{N}}$ | | 0.0492 | 0.2348 | 0.6824 | 0.9492 | 1.000 |
| $\underset{\sim}{\phi}_{\text{vdW}}$ | | 0.0460 | 0.2208 | 0.6652 | 0.9432 | 1.000 |
| $\underset{\sim}{\phi}_{f_{1,6}}$ | | 0.0468 | 0.2260 | 0.6644 | 0.9404 | 0.957 |
| $\underset{\sim}{\phi}_{f_{1,2}} = \underset{\sim}{\phi}_W$ | | 0.0544 | 0.2052 | 0.6036 | 0.9028 | 0.844 |
| $\underset{\sim}{\phi}_{f_{1,1}}$ | | 0.0544 | 0.1900 | 0.5532 | 0.8600 | 0.741 |
| $\underset{\sim}{\phi}_{f_{1,0.5}}$ | | 0.0560 | 0.1732 | 0.5000 | 0.8024 | 0.648 |
| $\underset{\sim}{\phi}_{f_{1,0.2}}$ | | 0.0560 | 0.1628 | 0.4536 | 0.7476 | 0.568 |
| $\underset{\sim}{\phi}_S$ | | 0.0568 | 0.1484 | 0.4016 | 0.6908 | 0.500 |
| $\underset{\sim}{\phi}_{SP}$ | | 0.0460 | 0.2180 | 0.6576 | 0.9356 | 0.934 |
| $\phi_{\text{John}}$ | $t_6$ | 0.1928 | 0.3712 | 0.7016 | 0.9092 | ND |
| $\phi_{\mathcal{N}}$ | | 0.0480 | 0.1580 | 0.4528 | 0.7608 | 1.000 |
| $\underset{\sim}{\phi}_{\text{vdW}}$ | | 0.0428 | 0.1816 | 0.5708 | 0.8800 | 1.531 |
| $\underset{\sim}{\phi}_{f_{1,6}}$ | | 0.0460 | 0.1956 | 0.5916 | 0.8956 | 1.600 |
| $\underset{\sim}{\phi}_{f_{1,2}} = \underset{\sim}{\phi}_W$ | | 0.0520 | 0.1904 | 0.5832 | 0.8860 | 1.531 |
| $\underset{\sim}{\phi}_{f_{1,1}}$ | | 0.0500 | 0.1836 | 0.5444 | 0.8588 | 1.408 |
| $\underset{\sim}{\phi}_{f_{1,0.5}}$ | | 0.0464 | 0.1708 | 0.4980 | 0.8148 | 1.269 |
| $\underset{\sim}{\phi}_{f_{1,0.2}}$ | | 0.0468 | 0.1480 | 0.4432 | 0.7648 | 1.172 |
| $\underset{\sim}{\phi}_S$ | | 0.0488 | 0.1284 | 0.3884 | 0.7064 | 1.000 |
| $\underset{\sim}{\phi}_{SP}$ | | 0.0480 | 0.1980 | 0.5956 | 0.8888 | 1.579 |
| $\phi_{\text{John}}$ | $t_1$ | 0.9868 | 0.9872 | 0.9848 | 0.9840 | ND |
| $\phi_{\mathcal{N}}$ | | 0.0060 | 0.0052 | 0.0064 | 0.0088 | ND |
| $\underset{\sim}{\phi}_{\text{vdW}}$ | | 0.0432 | 0.1244 | 0.3620 | 0.6508 | ND |
| $\underset{\sim}{\phi}_{f_{1,6}}$ | | 0.0456 | 0.1492 | 0.4256 | 0.7376 | ND |
| $\underset{\sim}{\phi}_{f_{1,2}} = \underset{\sim}{\phi}_W$ | | 0.0480 | 0.1636 | 0.4668 | 0.7936 | ND |
| $\underset{\sim}{\phi}_{f_{1,1}}$ | | 0.0468 | 0.1632 | 0.4724 | 0.8028 | ND |
| $\underset{\sim}{\phi}_{f_{1,0.5}}$ | | 0.0460 | 0.1636 | 0.4700 | 0.7964 | ND |
| $\underset{\sim}{\phi}_{f_{1,0.2}}$ | | 0.0428 | 0.1548 | 0.4404 | 0.7644 | ND |
| $\underset{\sim}{\phi}_S$ | | 0.0452 | 0.1408 | 0.4020 | 0.7064 | ND |
| $\underset{\sim}{\phi}_{SP}$ | | 0.0488 | 0.1444 | 0.4092 | 0.7240 | ND |

*Continued*

Table 3 (*Continued*)

| Test | | $m$ 0 | 1 | 2 | 3 | ARE |
|---|---|---|---|---|---|---|
| $\phi_{\text{John}}$ | $t_{0.2}$ | 0.9468 | 0.9460 | 0.9460 | 0.9500 | ND |
| $\phi_{\mathcal{N}}$ | | 0.0196 | 0.0184 | 0.0252 | 0.0352 | ND |
| $\underset{\sim}{\phi}_{\text{vdW}}$ | | 0.0412 | 0.0924 | 0.2468 | 0.4644 | ND |
| $\underset{\sim}{\phi}_{f_{1,6}}$ | | 0.0452 | 0.1144 | 0.2996 | 0.5572 | ND |
| $\underset{\sim}{\phi}_{f_{1,2}} = \underset{\sim}{\phi}_{W}$ | | 0.0528 | 0.1284 | 0.3460 | 0.6220 | ND |
| $\underset{\sim}{\phi}_{f_{1,1}}$ | | 0.0544 | 0.1348 | 0.3760 | 0.6672 | ND |
| $\underset{\sim}{\phi}_{f_{1,0.5}}$ | | 0.0476 | 0.1356 | 0.3908 | 0.6996 | ND |
| $\underset{\sim}{\phi}_{f_{1,0.2}}$ | | 0.0500 | 0.1372 | 0.3940 | 0.7016 | ND |
| $\underset{\sim}{\phi}_{S}$ | | 0.0468 | 0.1296 | 0.3724 | 0.6764 | ND |
| $\underset{\sim}{\phi}_{SP}$ | | 0.0468 | 0.1056 | 0.2752 | 0.5100 | ND |
| $\phi_{\text{John}}$ | $\mathcal{SN}$ | 0.0520 | 0.0624 | 0.2596 | 0.8000 | ? |
| $\phi_{\mathcal{N}}$ | | 0.0528 | 0.0664 | 0.2600 | 0.8000 | ? |
| $\underset{\sim}{\phi}_{\text{vdW}}$ | | 0.0472 | 0.0608 | 0.2488 | 0.7828 | ? |
| $\underset{\sim}{\phi}_{f_{1,6}}$ | | 0.0508 | 0.0620 | 0.2456 | 0.7808 | ? |
| $\underset{\sim}{\phi}_{f_{1,2}} = \underset{\sim}{\phi}_{W}$ | | 0.0492 | 0.0620 | 0.2304 | 0.7336 | ? |
| $\underset{\sim}{\phi}_{f_{1,1}}$ | | 0.0488 | 0.0608 | 0.2012 | 0.6784 | ? |
| $\underset{\sim}{\phi}_{f_{1,0.5}}$ | | 0.0476 | 0.0620 | 0.1796 | 0.6112 | ? |
| $\underset{\sim}{\phi}_{f_{1,0.2}}$ | | 0.0492 | 0.0568 | 0.1568 | 0.5540 | ? |
| $\underset{\sim}{\phi}_{S}$ | | 0.0512 | 0.0544 | 0.1412 | 0.4972 | ? |
| $\underset{\sim}{\phi}_{SP}$ | | 0.0528 | 0.0652 | 0.2504 | 0.7752 | ? |
| $\phi_{\text{John}}$ | $\mathcal{S}t_2$ | 0.8640 | 0.8616 | 0.9044 | 0.9520 | ? |
| $\phi_{\mathcal{N}}$ | | 0.0196 | 0.0188 | 0.0640 | 0.1896 | ? |
| $\underset{\sim}{\phi}_{\text{vdW}}$ | | 0.0536 | 0.0740 | 0.4144 | 0.8504 | ? |
| $\underset{\sim}{\phi}_{f_{1,6}}$ | | 0.0536 | 0.0724 | 0.4184 | 0.8276 | ? |
| $\underset{\sim}{\phi}_{f_{1,2}} = \underset{\sim}{\phi}_{W}$ | | 0.0512 | 0.0744 | 0.3592 | 0.6964 | ? |
| $\underset{\sim}{\phi}_{f_{1,1}}$ | | 0.0472 | 0.0724 | 0.2964 | 0.5048 | ? |
| $\underset{\sim}{\phi}_{f_{1,0.5}}$ | | 0.0484 | 0.0720 | 0.2324 | 0.3280 | ? |
| $\underset{\sim}{\phi}_{f_{1,0.2}}$ | | 0.0464 | 0.0688 | 0.1744 | 0.2076 | ? |
| $\underset{\sim}{\phi}_{S}$ | | 0.0468 | 0.0604 | 0.1524 | 0.1556 | ? |
| $\underset{\sim}{\phi}_{SP}$ | | 0.0552 | 0.0756 | 0.4592 | 0.8820 | ? |

TABLE 4

*Rejection frequencies (still out of $N = 2{,}500$ replications) under spherical and elliptic Gaussian and $t_{0.2}$ distributions, of the same tests as in Table 3; the sample size is now 25*

| | | *m* | | | | |
|---|---|---|---|---|---|---|
| **Test** | | **0** | **1** | **2** | **3** | **ARE** |
| $\phi_{\text{John}}$ | $\mathcal{N}$ | 0.0412 | 0.6032 | 0.9252 | 0.9860 | 1.000 |
| $\phi_{\mathcal{N}}$ | | 0.0424 | 0.5848 | 0.8924 | 0.9708 | 1.000 |
| $\underset{\sim}{\phi}_{\text{vdW}}$ | | 0.0172 | 0.4136 | 0.8088 | 0.9408 | 1.000 |
| $\underset{\sim}{\phi}_{f_{1,6}}$ | | 0.0356 | 0.5280 | 0.8684 | 0.9628 | 0.957 |
| $\underset{\sim}{\phi}_{f_{1,2}} = \underset{\sim}{\phi}_W$ | | 0.0416 | 0.5400 | 0.8612 | 0.9584 | 0.844 |
| $\underset{\sim}{\phi}_{f_{1,1}}$ | | 0.0468 | 0.5036 | 0.8316 | 0.9432 | 0.741 |
| $\underset{\sim}{\phi}_{f_{1,0.5}}$ | | 0.0496 | 0.4500 | 0.7924 | 0.9132 | 0.648 |
| $\underset{\sim}{\phi}_{f_{1,0.2}}$ | | 0.0484 | 0.4016 | 0.7328 | 0.8724 | 0.568 |
| $\underset{\sim}{\phi}_S$ | | 0.0480 | 0.3580 | 0.6736 | 0.8216 | 0.500 |
| $\underset{\sim}{\phi}_{SP}$ | | 0.0396 | 0.5600 | 0.8856 | 0.9696 | 0.934 |
| $\phi_{\text{John}}$ | $t_{0.2}$ | 0.8652 | 0.9076 | 0.9360 | 0.9484 | ND |
| $\phi_{\mathcal{N}}$ | | 0.0004 | 0.0008 | 0.0016 | 0.0020 | ND |
| $\underset{\sim}{\phi}_{\text{vdW}}$ | | 0.0148 | 0.1476 | 0.3608 | 0.5192 | ND |
| $\underset{\sim}{\phi}_{f_{1,6}}$ | | 0.0308 | 0.2492 | 0.5080 | 0.6844 | ND |
| $\underset{\sim}{\phi}_{f_{1,2}} = \underset{\sim}{\phi}_W$ | | 0.0452 | 0.3288 | 0.6168 | 0.7968 | ND |
| $\underset{\sim}{\phi}_{f_{1,1}}$ | | 0.0496 | 0.3592 | 0.6784 | 0.8376 | ND |
| $\underset{\sim}{\phi}_{f_{1,0.5}}$ | | 0.0488 | 0.3824 | 0.7172 | 0.8584 | ND |
| $\underset{\sim}{\phi}_{f_{1,0.2}}$ | | 0.0508 | 0.3892 | 0.7272 | 0.8692 | ND |
| $\underset{\sim}{\phi}_S$ | | 0.0480 | 0.3752 | 0.7044 | 0.8504 | ND |
| $\underset{\sim}{\phi}_{SP}$ | | 0.0348 | 0.2320 | 0.4620 | 0.6352 | ND |

LEMMA A.1. *Let Assumptions (A1) and (A2) hold. Define*

$$\underline{g}_{\boldsymbol{\theta},\boldsymbol{\Sigma};f_1}(\mathbf{x}) := c_{k,f_1}|\boldsymbol{\Sigma}|^{-1/2} f_1(\|\mathbf{x}-\boldsymbol{\theta}\|_{\boldsymbol{\Sigma}}), \qquad \mathbf{x} \in \mathbb{R}^k,$$

$$D_{\boldsymbol{\theta}}\underline{g}^{1/2}_{\boldsymbol{\theta},\boldsymbol{\Sigma};f_1}(\mathbf{x}) := \tfrac{1}{2}\underline{g}^{1/2}_{\boldsymbol{\theta},\boldsymbol{\Sigma};f_1}(\mathbf{x})\varphi_{f_1}(\|\mathbf{x}-\boldsymbol{\theta}\|_{\boldsymbol{\Sigma}})\boldsymbol{\Sigma}^{-1/2}\mathbf{u}(\boldsymbol{\theta},\boldsymbol{\Sigma})$$

*and*

$$D_{\boldsymbol{\Sigma}}\underline{g}^{1/2}_{\boldsymbol{\theta},\boldsymbol{\Sigma};f_1}(\mathbf{x}) := \tfrac{1}{4}\underline{g}^{1/2}_{\boldsymbol{\theta},\boldsymbol{\Sigma};f_1}(\mathbf{x})\mathbf{P}_k(\boldsymbol{\Sigma}^{\otimes 2})^{-1/2} \times \operatorname{vec}(\psi_{f_1}(\|\mathbf{x}-\boldsymbol{\theta}\|_{\boldsymbol{\Sigma}})\|\mathbf{x}-\boldsymbol{\theta}\|_{\boldsymbol{\Sigma}}\mathbf{u}(\boldsymbol{\theta},\boldsymbol{\Sigma})\mathbf{u}'(\boldsymbol{\theta},\boldsymbol{\Sigma}) - \mathbf{I}_k),$$

*where* $\|\mathbf{z}\|_{\mathbf{\Sigma}} := (\mathbf{z}'\mathbf{\Sigma}^{-1}\mathbf{z})^{1/2}$, $\mathbf{u}(\boldsymbol{\theta},\mathbf{\Sigma}) := \mathbf{\Sigma}^{-1/2}(\mathbf{x}-\boldsymbol{\theta})/\|\mathbf{x}-\boldsymbol{\theta}\|_{\mathbf{\Sigma}}$ *and* $\mathbf{P}_k$ *is such that* $\mathbf{P}_k'(\operatorname{vech}\mathbf{H}) = \operatorname{vec}\mathbf{H}$ *for any symmetric* $k\times k$ *matrix* $\mathbf{H}=(H_{ij})$. *Then:*

(i) $\int \{g^{1/2}_{\boldsymbol{\theta}+\mathbf{t},\mathbf{\Sigma};f_1}(\mathbf{x}) - g^{1/2}_{\boldsymbol{\theta},\mathbf{\Sigma};f_1}(\mathbf{x}) - \mathbf{t}'(D_{\boldsymbol{\theta}} g^{1/2}_{\boldsymbol{\theta},\mathbf{\Sigma};f_1}(\mathbf{x}))\}^2\,d\mathbf{x} = o(\|\mathbf{t}\|^2)$,

(ii) $\int \{g^{1/2}_{\boldsymbol{\theta},\mathbf{\Sigma}+\mathbf{H};f_1}(\mathbf{x}) - g^{1/2}_{\boldsymbol{\theta},\mathbf{\Sigma};f_1}(\mathbf{x}) - (\operatorname{vech}\mathbf{H})'(D_{\mathbf{\Sigma}} g^{1/2}_{\boldsymbol{\theta},\mathbf{\Sigma};f_1}(\mathbf{x}))\}^2\,d\mathbf{x}$
$= o(\|\mathbf{H}\|^2)$ *and*

(iii) $\int \left\{ g^{1/2}_{\boldsymbol{\theta}+\mathbf{t},\mathbf{\Sigma}+\mathbf{H};f_1}(\mathbf{x}) - g^{1/2}_{\boldsymbol{\theta},\mathbf{\Sigma};f_1}(\mathbf{x}) - \begin{pmatrix}\mathbf{t}\\ \operatorname{vech}\mathbf{H}\end{pmatrix}' \begin{pmatrix} D_{\boldsymbol{\theta}} g^{1/2}_{\boldsymbol{\theta},\mathbf{\Sigma};f_1}(\mathbf{x}) \\ D_{\mathbf{\Sigma}} g^{1/2}_{\boldsymbol{\theta},\mathbf{\Sigma};f_1}(\mathbf{x})\end{pmatrix}\right\}^2 d\mathbf{x}$
$= o\left(\left\| \begin{pmatrix}\mathbf{t}\\ \operatorname{vech}\mathbf{H}\end{pmatrix}\right\|^2\right)$.

To prove Lemma A.1, we need the following reformulation of Assumption (A2):

LEMMA A.2. *Assumption* (A2) *holds iff* (i) $f^{1/2}_{1;\exp} \in L^2(\mathbb{R},\nu_k)$ *and* (ii) *there exists* $Df^{1/2}_{1;\exp} \in L^2(\mathbb{R},\nu_k)$ *such that*

$$\int [f^{1/2}_{1;\exp}(x+h) - f^{1/2}_{1;\exp}(x) - h(Df^{1/2}_{1;\exp})(x)]^2 e^{kx}\,dx = o(h^2)$$

*as* $h\to 0$. *In that case,* $Df^{1/2}_{1;\exp}$ *and* $(f^{1/2}_{1;\exp})'$ *are equal in* $L^2(\mathbb{R},\nu_k)$.

The proof of this lemma relies on the following result by Schwartz (see [47], pages 186–188):

LEMMA A.3 (Schwartz). *The real function* $g$ *is in* $W^{1,2}(\mathbb{R})$ (*with weak derivative* $g'$, *say*) *iff* (i) $g\in L^2(\mathbb{R})$ *and* (ii) *there exists* $Dg\in L^2(\mathbb{R})$ *such that* $x\mapsto g(x+h)-g(x)-h(Dg(x))$ *is* $o(h)$ *in* $L^2(\mathbb{R})$ (*as* $h\to 0$), *that is,* $\int [g(x+h)-g(x)-h(Dg(x))]^2\,dx = o(h^2)$ *as* $h\to 0$. *In that case,* $Dg$ *and* $g'$ *are equal in* $L^2(\mathbb{R})$.

PROOF OF LEMMA A.2. Throughout this proof, we write $f$ instead of $f^{1/2}_{1;\exp}$ and all $o(h)$'s are taken as $h\to 0$.

(*Necessity*) It is easy to show that the real function $x\mapsto g(x) := f(x)e^{kx/2}$ admits the weak derivative $x\mapsto g'(x) = f'(x)e^{kx/2} + (k/2)g(x)$, where $f'$ denotes the weak derivative of $f$. In view of the assumptions on $f$, both $g$ and $g'$ are in $L^2(\mathbb{R})$. Lemma A.3 therefore yields that $x\mapsto M_h(x) := g(x+h) - g(x) - hg'(x)$ is $o(h)$ in $L^2(\mathbb{R})$. But $M_h = I_h + J_h + K_h + L_h$, where

$$I_h(x) := (f(x+h) - f(x) - hf'(x))e^{kx/2},$$

$$J_h(x) := f(x+h)e^{k(x+h)/2}e^{-kh/2}(e^{kh/2} - 1 - hk/2),$$
$$K_h(x) := (f(x+h)e^{k(x+h)/2} - f(x)e^{kx/2})hk/2$$

and

$$L_h(x) := f(x+h)e^{k(x+h)/2}(e^{-kh/2} - 1)hk/2.$$

Since $J_h$, $K_h$ and $L_h$ are also $o(h)$ in $L^2(\mathbb{R})$, so is $I_h$.

(*Sufficiency*) Assume now that $f \in L^2(\mathbb{R}, \nu_k)$ is such that $x \mapsto I_h(x) := (f(x+h) - f(x) - hDf(x))e^{kx/2}$ is $o(h)$ in $L^2(\mathbb{R})$ for some $Df \in L^2(\mathbb{R}, \nu_k)$ and again define $x \mapsto g(x) := f(x)e^{kx/2}$ $[g \in L^2(\mathbb{R})]$. With $Dg(x) := Df(x)e^{kx/2} + (k/2)g(x)$ $[Dg \in L^2(\mathbb{R})]$, we have that

$$\begin{aligned} x \mapsto \widetilde{M}_h(x) &:= g(x+h) - g(x) - hDg(x) \\ &= (f(x+h) - f(x) - hDf(x))e^{kx/2} + J_h(x) + K_h(x) + L_h(x) \end{aligned}$$

is $o(h)$ in $L^2(\mathbb{R})$. Lemma A.3 thus yields that $Dg$ is the weak derivative of $g$; this implies that, for all infinitely differentiable compactly supported functions $\varphi$,

$$\begin{aligned} &\int [\varphi(x)e^{-kx/2}][Df(x)e^{kx/2} + (k/2)g(x)]\,dx \\ &\qquad = -\int [\varphi'(x)e^{-kx/2} - (k/2)\varphi(x)e^{-kx/2}][f(x)e^{kx/2}]\,dx, \end{aligned}$$

that is, that $Df$ is the weak derivative of $f$. □

Proof of Lemma A.1. (i) See [17].

(ii) Using the fact that $(\mathbf{C}' \otimes \mathbf{A}) \operatorname{vec} \mathbf{B} = \operatorname{vec}(\mathbf{ABC})$ and letting $\mathbf{y} := \mathbf{\Sigma}^{-1/2}(\mathbf{x} - \boldsymbol{\theta})$, the left-hand side of (ii) takes the form

$$\begin{aligned} &c_{k,f_1} \int \biggl\{ \frac{1}{|\mathbf{I}_k + \mathbf{H}_{\mathbf{\Sigma}}|^{1/4}} f_1^{1/2}(\|\mathbf{y}\|_{\mathbf{I}_k + \mathbf{H}_{\mathbf{\Sigma}}}) - f_1^{1/2}(\|\mathbf{y}\|) - \frac{1}{4} f_1^{1/2}(\|\mathbf{y}\|) \\ &\qquad\qquad \times (\operatorname{vec} \mathbf{H}_{\mathbf{\Sigma}})' \operatorname{vec}\biggl(\psi_{f_1}(\|\mathbf{y}\|)\frac{\mathbf{y}\mathbf{y}'}{\|\mathbf{y}\|} - \mathbf{I}_k\biggr)\biggr\}^2 d\mathbf{y} \\ &\quad \le C(T_1 + T_2 + T_3), \end{aligned}$$

where $\mathbf{H}_{\mathbf{\Sigma}} := \mathbf{\Sigma}^{-1/2}\mathbf{H}\mathbf{\Sigma}^{-1/2}$, $C$ is some positive constant,

$$T_1 := \int \biggl\{ \frac{1}{|\mathbf{I}_k + \mathbf{H}_{\mathbf{\Sigma}}|^{1/4}} - 1 + \frac{1}{4}(\operatorname{vec} \mathbf{H}_{\mathbf{\Sigma}})'(\operatorname{vec} \mathbf{I}_k) \biggr\}^2 f_1(\|\mathbf{y}\|_{\mathbf{I}_k + \mathbf{H}_{\mathbf{\Sigma}}})\,d\mathbf{y},$$

$$T_2 := \int \tfrac{1}{16}[(\operatorname{vec} \mathbf{H}_{\mathbf{\Sigma}})'(\operatorname{vec} \mathbf{I}_k)]^2 \{f_1^{1/2}(\|\mathbf{y}\|_{\mathbf{I}_k + \mathbf{H}_{\mathbf{\Sigma}}}) - f_1^{1/2}(\|\mathbf{y}\|)\}^2\,d\mathbf{y}$$

and

$$T_3 := \int \biggl\{ f_1^{1/2}(\|\mathbf{y}\|_{\mathbf{I}_k+\mathbf{H}_\mathbf{\Sigma}}) - f_1^{1/2}(\|\mathbf{y}\|) \\ - \frac{1}{4} f_1^{1/2}(\|\mathbf{y}\|)(\mathrm{vec}\,\mathbf{H}_\mathbf{\Sigma})' \,\mathrm{vec}\biggl(\psi_{f_1}(\|\mathbf{y}\|)\frac{\mathbf{y}\mathbf{y}'}{\|\mathbf{y}\|}\biggr)\biggr\}^2 \,d\mathbf{y}.$$

Since $(\mathrm{vec}\,\mathbf{A})'(\mathrm{vec}\,\mathbf{B}) = \mathrm{tr}(\mathbf{A}'\mathbf{B})$ and $|\mathbf{A}+\mathbf{B}|^a = |\mathbf{A}|^a + a|\mathbf{A}|^a\,\mathrm{tr}(\mathbf{A}^{-1}\mathbf{B}) + o(\|\mathbf{B}\|)$ for all $a$ (see, e.g., [31], page 149),

$$T_1 = \frac{|\mathbf{I}_k + \mathbf{H}_\mathbf{\Sigma}|^{1/2}}{c_{k,f_1}} \biggl\{ |\mathbf{I}_k + \mathbf{H}_\mathbf{\Sigma}|^{-1/4} - 1 + \frac{1}{4}(\mathrm{tr}\,\mathbf{H}_\mathbf{\Sigma})\biggr\}^2 = o(\|\mathbf{H}\|^2).$$

Now, working in spherical coordinates $(r,\mathbf{u}) := (\|\mathbf{y}\|, \mathbf{y}/\|\mathbf{y}\|)$, we obtain

$$\begin{aligned} T_3 &= C \iint \{ f_1^{1/2}(r\|\mathbf{u}\|_{\mathbf{I}_k+\mathbf{H}_\mathbf{\Sigma}}) - f_1^{1/2}(r) \\ &\qquad - \tfrac{1}{4} f_1^{1/2}(r)\psi_{f_1}(r) r[\mathbf{u}'\mathbf{H}_\mathbf{\Sigma}\mathbf{u}]\}^2 r^{k-1}\,dr\,d\sigma(\mathbf{u}) \\ &= C \iint \{ f_{1;\exp}^{1/2}((\ln r) + (\ln\|\mathbf{u}\|_{\mathbf{I}_k+\mathbf{H}_\mathbf{\Sigma}})) - f_{1;\exp}^{1/2}(\ln r) \\ &\qquad + (f_{1;\exp}^{1/2})'(\ln r)[\tfrac{1}{2}\mathbf{u}'\mathbf{H}_\mathbf{\Sigma}\mathbf{u}]\}^2 r^{k-1}\,dr\,d\sigma(\mathbf{u}) \\ &= C \iint \{ f_{1;\exp}^{1/2}(s + (\ln\|\mathbf{u}\|_{\mathbf{I}_k+\mathbf{H}_\mathbf{\Sigma}})) \\ &\qquad - f_{1;\exp}^{1/2}(s) + (f_{1;\exp}^{1/2})'(s)[\tfrac{1}{2}\mathbf{u}'\mathbf{H}_\mathbf{\Sigma}\mathbf{u}]\}^2 e^{ks}\,ds\,d\sigma(\mathbf{u}) \\ &\le C(T_{3a} + T_{3b}), \end{aligned}$$

where

$$T_{3a} := \iint \{ f_{1;\exp}^{1/2}(s + (\ln\|\mathbf{u}\|_{\mathbf{I}_k+\mathbf{H}_\mathbf{\Sigma}})) \\ - f_{1;\exp}^{1/2}(s) - (f_{1;\exp}^{1/2})'(s)[\ln\|\mathbf{u}\|_{\mathbf{I}_k+\mathbf{H}_\mathbf{\Sigma}}]\}^2 e^{ks}\,ds\,d\sigma(\mathbf{u})$$

and

$$T_{3b} := \iint \{[\ln\|\mathbf{u}\|_{\mathbf{I}_k+\mathbf{H}_\mathbf{\Sigma}}] + [\tfrac{1}{2}\mathbf{u}'\mathbf{H}_\mathbf{\Sigma}\mathbf{u}]\}^2 [(f_{1;\exp}^{1/2})'(s)]^2 e^{ks}\,ds\,d\sigma(\mathbf{u}).$$

By using Lemma A.2 and the fact that $\ln\|\mathbf{u}\|_{\mathbf{I}_k+\mathbf{H}_\mathbf{\Sigma}} = O(\|\mathbf{H}\|)$ for all $\mathbf{u}$, we obtain that

$$\int \{ f_{1;\exp}^{1/2}(s + (\ln\|\mathbf{u}\|_{\mathbf{I}_k+\mathbf{H}_\mathbf{\Sigma}})) - f_{1;\exp}^{1/2}(s) - (f_{1;\exp}^{1/2})'(s)[\ln\|\mathbf{u}\|_{\mathbf{I}_k+\mathbf{H}_\mathbf{\Sigma}}]\}^2 e^{ks}\,ds \\ = o(\|\mathbf{H}\|^2),$$

for all $\mathbf{u}$. Therefore, from Lebesgue's dominated convergence theorem, it follows that $T_{3a} = o(\|\mathbf{H}\|^2)$. As for $T_{3b}$, we have that

$$T_{3b} \leq \sup_{\mathbf{u}\in\mathcal{S}^{k-1}} \{[\ln\|\mathbf{u}\|_{\mathbf{I}_k+\mathbf{H}_{\boldsymbol{\Sigma}}}] + [\tfrac{1}{2}\mathbf{u}'\mathbf{H}_{\boldsymbol{\Sigma}}\mathbf{u}]\}^2 = o(\|\mathbf{H}\|^2)$$

since $[\ln\|\mathbf{u}\|_{\mathbf{I}_k+\mathbf{H}_{\boldsymbol{\Sigma}}}] + [\frac{1}{2}\mathbf{u}'\mathbf{H}_{\boldsymbol{\Sigma}}\mathbf{u}] = o(\|\mathbf{H}\|)$, uniformly for $\mathbf{u}\in\mathcal{S}^{k-1}$ (see, e.g., [31], page 151). Consequently, $T_3 = o(\|\mathbf{H}\|^2)$, so $T_3 = o(1)$ as $\|\mathbf{H}\|$ goes to zero, and hence

$$\begin{aligned} T_2 &\leq C\|\mathbf{H}_{\boldsymbol{\Sigma}}\|^2 \int \{f_1^{1/2}(\|\mathbf{y}\|_{\mathbf{I}_k+\mathbf{H}_{\boldsymbol{\Sigma}}}) - f_1^{1/2}(\|\mathbf{y}\|)\}^2\,d\mathbf{y} \\ &\leq C\|\mathbf{H}_{\boldsymbol{\Sigma}}\|^2 \int \left\{\frac{1}{4} f_1^{1/2}(\|\mathbf{y}\|)(\operatorname{vec}\mathbf{H}_{\boldsymbol{\Sigma}})' \operatorname{vec}\left(\psi_{f_1}(\|\mathbf{y}\|)\frac{\mathbf{y}\mathbf{y}'}{\|\mathbf{y}\|}\right)\right\}^2 d\mathbf{y} + o(\|\mathbf{H}\|^2), \end{aligned}$$

which shows that $T_2 = o(\|\mathbf{H}\|^2)$. This proves (ii).

(iii) The left-hand side in (iii) is bounded by $C(S_1 + S_2 + \|\operatorname{vech}\mathbf{H}\|^2 S_3)$, where

$$S_1 := \int \{g^{1/2}_{\boldsymbol{\theta}+\mathbf{t},\boldsymbol{\Sigma};f_1}(\mathbf{x}) - g^{1/2}_{\boldsymbol{\theta},\boldsymbol{\Sigma};f_1}(\mathbf{x}) - \mathbf{t}'(D_{\boldsymbol{\theta}} g^{1/2}_{\boldsymbol{\theta},\boldsymbol{\Sigma};f_1}(\mathbf{x}))\}^2\,d\mathbf{x},$$

$$S_2 := \int \{g^{1/2}_{\boldsymbol{\theta},\boldsymbol{\Sigma}+\mathbf{H};f_1}(\mathbf{x}) - g^{1/2}_{\boldsymbol{\theta},\boldsymbol{\Sigma};f_1}(\mathbf{x}) - (\operatorname{vech}\mathbf{H})'(D_{\boldsymbol{\Sigma}} g^{1/2}_{\boldsymbol{\theta},\boldsymbol{\Sigma};f_1}(\mathbf{x}))\}^2\,d\mathbf{x}$$

and

$$\begin{aligned} S_3 &:= \int \|D_{\boldsymbol{\Sigma}} g^{1/2}_{\boldsymbol{\theta}+\mathbf{t},\boldsymbol{\Sigma};f_1}(\mathbf{x}) - D_{\Sigma} g^{1/2}_{\boldsymbol{\theta},\boldsymbol{\Sigma};f_1}(\mathbf{x})\|^2\,d\mathbf{x} \\ &= \int \|D_{\boldsymbol{\Sigma}} g^{1/2}_{\boldsymbol{\theta},\boldsymbol{\Sigma};f_1}(\mathbf{x}-\mathbf{t}) - D_{\Sigma} g^{1/2}_{\boldsymbol{\theta},\boldsymbol{\Sigma};f_1}(\mathbf{x})\|^2\,d\mathbf{x}. \end{aligned}$$

Now, from (i) and (ii), respectively, $S_1$ and $S_2$ are $o(\|(\mathbf{t}' \vdots (\operatorname{vech}\mathbf{H})')'\|^2)$. As for $S_3$, the quadratic mean continuity of $\mathbf{x} \to D_{\boldsymbol{\Sigma}} g^{1/2}_{\boldsymbol{\theta},\boldsymbol{\Sigma};f_1}(\mathbf{x}) \in L^2(\mathbb{R}^k)$ implies that it is $o(1)$ as $\mathbf{t}\to\mathbf{0}$. The result follows. □

LEMMA A.4. *Let $\mathbf{x} \mapsto G_{\boldsymbol{\eta}}(\mathbf{x})$ be differentiable in quadratic mean at $\boldsymbol{\eta}_0$, with gradient $\mathbf{x} \mapsto DG_{\boldsymbol{\eta}_0}(\mathbf{x})$, say. Let $h$ be a diffeomorphism in a neighborhood of $\boldsymbol{\xi}_0 := h^{-1}(\boldsymbol{\eta}_0)$. Then $\mathbf{x} \mapsto G_{h(\boldsymbol{\xi})}(\mathbf{x})$ is differentiable in quadratic mean at $\boldsymbol{\xi}_0$, with gradient $\mathbf{x} \mapsto (Dh_{\boldsymbol{\xi}_0})'(DG_{h(\boldsymbol{\xi}_0)}(\mathbf{x}))$, where $Dh_{\boldsymbol{\xi}_0} := (\frac{\partial h_i}{\partial \boldsymbol{\xi}_j}(\boldsymbol{\xi}_0))$ denotes the Jacobian matrix of $h$ at $\boldsymbol{\xi}_0$.*

PROOF OF LEMMA A.4. This is straightforward. □

Applied to Lemma A.1(iii), the latter result implies that $\mathbf{x} \mapsto \underline{f}^{1/2}_{\boldsymbol{\vartheta};f_1}(\mathbf{x}) = \underline{f}^{1/2}_{\boldsymbol{\theta},\sigma^2,\mathbf{V};f_1}(\mathbf{x}) = \underline{g}^{1/2}_{\boldsymbol{\theta},\sigma^2\mathbf{V};f_1}(\mathbf{x})$ is differentiable in quadratic mean, with gradi-

ent

$$D\underline{f}^{1/2}_{\boldsymbol{\vartheta};f_1}(\mathbf{x}) = \begin{pmatrix} D_{\boldsymbol{\theta}} g^{1/2}_{\boldsymbol{\theta},\sigma^2\mathbf{V};f_1}(\mathbf{x}) \\ \begin{pmatrix} 1 & (\operatorname{ve\overset{\circ}{c}h}\mathbf{V})' \\ \mathbf{0} & \sigma^2\mathbf{I} \end{pmatrix} D_{\boldsymbol{\Sigma}} g^{1/2}_{\boldsymbol{\theta},\sigma^2\mathbf{V};f_1}(\mathbf{x}) \end{pmatrix} = \frac{1}{2}\underline{f}^{1/2}_{\boldsymbol{\vartheta};f_1}(\mathbf{x}) W_{\boldsymbol{\vartheta};f_1}(\mathbf{x}),$$

where

$$W_{\boldsymbol{\vartheta};f_1}(\mathbf{x}) := \begin{pmatrix} \frac{1}{\sigma}\varphi_{f_1}\left(\frac{\|\mathbf{x}-\boldsymbol{\theta}\|_{\mathbf{V}}}{\sigma}\right)\mathbf{V}^{-1/2}\mathbf{u}(\boldsymbol{\theta},\mathbf{V}) \\ \frac{1}{2}\begin{pmatrix} \sigma^{-2}(\operatorname{vec}\mathbf{I}_k)' \\ \mathbf{M}_k(\mathbf{V}^{\otimes 2})^{-1/2} \end{pmatrix} \\ \times \operatorname{vec}\left(\psi_{f_1}\left(\frac{\|\mathbf{x}-\boldsymbol{\theta}\|_{\mathbf{V}}}{\sigma}\right)\frac{\|\mathbf{x}-\boldsymbol{\theta}\|_{\mathbf{V}}}{\sigma}\mathbf{u}(\boldsymbol{\theta},\mathbf{V})\mathbf{u}'(\boldsymbol{\theta},\mathbf{V}) - \mathbf{I}_k\right) \end{pmatrix}.$$

Checking Swensen's sufficient conditions for LAN is then a routine task. For example, letting $\nu_i^{(n)} := (\underline{f}^{1/2}_{\boldsymbol{\vartheta}+n^{-1/2}\boldsymbol{\tau}^{(n)};f_1}(\mathbf{X}_i)/\underline{f}^{1/2}_{\boldsymbol{\vartheta};f_1}(\mathbf{X}_i)) - 1$ and $Z_i^{(n)} := (1/2)(\boldsymbol{\tau}^{(n)})' n^{-1/2} W_{\boldsymbol{\vartheta};f_1}(\mathbf{X}_i)$, $i = 1,\ldots,n$, we have

$$\begin{aligned} &\mathrm{E}\left[\sum_{i=1}^n (\nu_i^{(n)} - Z_i^{(n)})^2\right] \\ &\quad = n\int \{\underline{f}^{1/2}_{\boldsymbol{\vartheta}+n^{-1/2}\boldsymbol{\tau}^{(n)};f_1}(\mathbf{x}) \\ &\qquad\quad - \underline{f}^{1/2}_{\boldsymbol{\vartheta};f_1}(\mathbf{x}) - (1/2)(\boldsymbol{\tau}^{(n)})' n^{-1/2}\underline{f}^{1/2}_{\boldsymbol{\vartheta};f_1}(\mathbf{x}) W_{\boldsymbol{\vartheta};f_1}(\mathbf{x})\}^2\,d\mathbf{x} \\ &\quad = n\int \{\underline{f}^{1/2}_{\boldsymbol{\vartheta}+n^{-1/2}\boldsymbol{\tau}^{(n)};f_1}(\mathbf{x}) - \underline{f}^{1/2}_{\boldsymbol{\vartheta};f_1}(\mathbf{x}) - (n^{-1/2}\boldsymbol{\tau}^{(n)})'(D\underline{f}^{1/2}_{\boldsymbol{\vartheta};f_1}(\mathbf{x}))\}^2\,d\mathbf{x}, \end{aligned}$$

which is $o(1)$ as $n \to \infty$. The other conditions easily follow. Now, the linear term in the second-order decomposition of the local log-likelihood ratio is $2\sum_{i=1}^n Z_i^{(n)} = (\boldsymbol{\tau}^{(n)})'\boldsymbol{\Delta}^{(n)}_{f_1}(\boldsymbol{\vartheta})$, where $\boldsymbol{\Delta}^{(n)}_{f_1}(\boldsymbol{\vartheta})$ is the central sequence given in (2.5).

### A.2. Proofs of Lemma 3.1 and Proposition 3.1.

PROOF OF LEMMA 3.1. Denote by $\mathbf{Q}_k(\mathbf{V})$ the matrix in the right-hand side of (3.3). Tedious but routine algebra yields

$$\mathbf{N}_k\mathbf{Q}_k(\mathbf{V})\mathbf{N}_k' = \frac{1}{k(k+2)}\boldsymbol{\Upsilon}_k(\mathbf{V})$$

(where $\mathbf{N}_k$ is defined in Section 1.4). In order to prove the lemma, it is therefore sufficient to show that $\mathbf{M}_k'\mathbf{N}_k\mathbf{Q}_k(\mathbf{V}) = \mathbf{Q}_k(\mathbf{V})$. Now, it is easily seen that

$$\mathbf{Q}_k(\mathbf{V}) = [\mathbf{I}_{k^2} - (\operatorname{vec}\mathbf{V})(\mathbf{e}_{k^2,1})'][\mathbf{I}_{k^2} + \mathbf{K}_k](\mathbf{V}^{\otimes 2})[\mathbf{I}_{k^2} - (\operatorname{vec}\mathbf{V})(\mathbf{e}_{k^2,1})']'.$$

But, letting $\mathbf{E}_{ij} := \mathbf{e}_i\mathbf{e}_j' + \mathbf{e}_j\mathbf{e}_i'$ [where $(\mathbf{e}_1, \ldots, \mathbf{e}_k)$ stands for the canonical basis of $\mathbb{R}^k$], we have

$$\begin{aligned}
&[\mathbf{I}_{k^2} - (\operatorname{vec}\mathbf{V})(\mathbf{e}_{k^2,1})'][\mathbf{I}_{k^2} + \mathbf{K}_k] \\
&\quad = \mathbf{I}_{k^2} + \mathbf{K}_k - 2(\operatorname{vec}\mathbf{V})(\mathbf{e}_{k^2,1})' \\
&\quad = \tfrac{1}{2} \sum_{\substack{i,j=1 \\ (i,j)\neq(1,1)}}^{k} (\operatorname{vec}\mathbf{E}_{ij})(\operatorname{vec}\mathbf{E}_{ij})' + 2(\operatorname{vec}(\mathbf{e}_1\mathbf{e}_1' - \mathbf{V}))(\mathbf{e}_{k^2,1})'.
\end{aligned}$$

The result follows, since $\mathbf{M}_k'\mathbf{N}_k(\operatorname{vec}\mathbf{W}) = (\operatorname{vec}\mathbf{W})$ for any symmetric $k \times k$ matrix $\mathbf{W} = (W_{ij})$ such that $W_{11} = 0$ [recall that it is assumed that $\mathbf{V} = (V_{ij})$ is symmetric with $V_{11} = 1$].

Proof of Proposition 3.1. Under $\mathrm{P}^{(n)}_{\boldsymbol{\vartheta}_0;f_1}$, for any fixed $\boldsymbol{\vartheta}_0' := (\boldsymbol{\theta}', \sigma^2, (\operatorname{vec}^{\circ}\mathrm{h}\,\mathbf{V}_0)')$, we have

$$Q^{(n)}_{f_1} = (\boldsymbol{\Delta}^{\star(n)}_{f_1}(\boldsymbol{\vartheta}_0))'(\boldsymbol{\Gamma}^{\star}_{f_1}(\boldsymbol{\vartheta}_0))^{-1}\boldsymbol{\Delta}^{\star(n)}_{f_1}(\boldsymbol{\vartheta}_0) + o_{\mathrm{P}}(1)$$

as $n \to \infty$. The proof of the first statement in part (i) of Proposition 3.1 follows, since $\boldsymbol{\Delta}^{\star(n)}_{f_1}(\boldsymbol{\vartheta}_0)$ is asymptotically $\mathcal{N}_{k(k+1)/2-1}(\mathbf{0}, \boldsymbol{\Gamma}^{\star}_{f_1}(\boldsymbol{\vartheta}_0))$ under $\mathrm{P}^{(n)}_{\boldsymbol{\vartheta}_0;f_1}$. On the other hand, it is easy to see, still under $\mathrm{P}^{(n)}_{\boldsymbol{\vartheta}_0;f_1}$, that $\boldsymbol{\Delta}^{\star(n)}_{f_1}(\boldsymbol{\vartheta}_0)$ and the local log-likelihood ratio $\Lambda^{(n)}_{\boldsymbol{\vartheta}_0+n^{-1/2}\boldsymbol{\tau}/\boldsymbol{\vartheta}_0;f_1}$, where $\boldsymbol{\tau}' := (\mathbf{t}', s, (\operatorname{vec}^{\circ}\mathrm{h}\,\mathbf{v})')$, are jointly multinormal, with asymptotic covariance $(\boldsymbol{\Gamma}^{\star}_{f_1}(\boldsymbol{\vartheta}_0))(\operatorname{vec}^{\circ}\mathrm{h}\,\mathbf{v})$. Le Cam's third lemma thus implies that $\boldsymbol{\Delta}^{\star(n)}_{f_1}(\boldsymbol{\vartheta}_0)$ is asymptotically

$$\mathcal{N}_{k(k+1)/2-1}((\boldsymbol{\Gamma}^{\star}_{f_1}(\boldsymbol{\vartheta}_0))(\operatorname{vec}^{\circ}\mathrm{h}\,\mathbf{v}), \boldsymbol{\Gamma}^{\star}_{f_1}(\boldsymbol{\vartheta}_0))$$

under $\mathrm{P}^{(n)}_{\boldsymbol{\vartheta}_0+n^{-1/2}\boldsymbol{\tau};f_1}$, which establishes the second statement in part (i) of the proposition.

As for part (ii), the fact that $\phi^{(n)}_{f_1}$ has asymptotic level $\alpha$ follows directly from the asymptotic null distribution given in part (i) and the classical Helly–Bray theorem, while local asymptotic maximinity is a consequence of the weak convergence to Gaussian shifts of local shape experiments (see, e.g., Section 11.9 of [29]). □

### A.3. Proofs of Propositions 3.2, 4.1 and Lemma 4.1.

Proof of Proposition 3.2. Under $\mathrm{P}^{(n)}_{\boldsymbol{\vartheta}_0;\phi_1}$, for any fixed $\boldsymbol{\vartheta}_0' := (\boldsymbol{\theta}', \sigma^2, (\operatorname{vec}^{\circ}\mathrm{h}\,\mathbf{V}_0)')$, we have

$$Q^{(n)}_{\mathcal{N}} = (\boldsymbol{\Delta}^{\star(n)}_{\phi_1}(\boldsymbol{\vartheta}_0))'(a_k^2 E_k(g_1)\boldsymbol{\Upsilon}_k^{-1}(\mathbf{V}_0))^{-1}\boldsymbol{\Delta}^{\star(n)}_{\phi_1}(\boldsymbol{\vartheta}_0) + o_{\mathrm{P}}(1)$$

as $n \to \infty$, where $\mathbf{\Upsilon}_k^{-1}(\mathbf{V}_0)$ was defined in (3.2). The result then follows—as in Proposition 3.1—by proving that, under $\mathrm{P}^{(n)}_{\boldsymbol{\vartheta}_0+n^{-1/2}\boldsymbol{\tau};g_1}$ [with $\boldsymbol{\tau}' := (\mathbf{t}', s, (\overset{\circ}{\mathrm{vech}}\,\mathbf{v})')$], we have

$$\boldsymbol{\Delta}^{\star(n)}_{\phi_1}(\boldsymbol{\vartheta}_0) \xrightarrow{\mathcal{L}} \mathcal{N}(a_k \mathrm{E}[\psi_{g_1}(\tilde{G}_1^{-1}(u))(\tilde{G}_1^{-1}(u))^3]\mathbf{\Upsilon}_k^{-1}(\mathbf{V}_0)(\overset{\circ}{\mathrm{vech}}\,\mathbf{v}), a_k^2 E_k(g_1)\mathbf{\Upsilon}_k^{-1}(\mathbf{V}_0))$$

[also note that integration by parts yields $\mathrm{E}[\psi_{g_1}(\tilde{G}_1^{-1}(u))(\tilde{G}_1^{-1}(u))^3] = (k+2)D_k(g_1)$]. As for the optimality statement in part (ii) of the proposition, it is obtained as in the proof of Proposition 3.1 and by noting that $a_k^2 E_k(\phi_1)\mathbf{\Upsilon}_k^{-1}(\mathbf{V}_0) = \mathbf{\Gamma}^{\star}_{\phi_1}(\boldsymbol{\vartheta}_0)$. □

PROOF OF LEMMA 4.1. Let

$$\underset{\sim}{\mathbf{T}}{}^{(n)}_{\boldsymbol{\vartheta};K} := n^{-1/2}\mathbf{J}_k^{\perp}\sum_{i=1}^n K\left(\frac{R_i}{n+1}\right)\mathrm{vec}(\mathbf{U}_i\mathbf{U}_i')$$

and

$$\mathbf{T}^{(n)}_{\boldsymbol{\vartheta};K;g_1} := n^{-1/2}\mathbf{J}_k^{\perp}\sum_{i=1}^n K\left(\tilde{G}_{1k}\left(\frac{d_i}{\sigma}\right)\right)\mathrm{vec}(\mathbf{U}_i\mathbf{U}_i').$$

Clearly, it is sufficient to prove that $\underset{\sim}{\mathbf{T}}{}^{(n)}_{\boldsymbol{\vartheta};K} - \mathbf{T}^{(n)}_{\boldsymbol{\vartheta};K;g_1}$ goes to zero in quadratic mean, under $\mathrm{P}^{(n)}_{\boldsymbol{\vartheta};g_1}$, as $n \to \infty$. For all $\ell = 1, 2, \ldots, k^2$, we have

$$\mathrm{E}[(\underset{\sim}{\mathbf{T}}{}^{(n)}_{\boldsymbol{\vartheta};K} - \mathbf{T}^{(n)}_{\boldsymbol{\vartheta};K;g_1})^2_\ell] = C_{\ell,k} n^{-1}\sum_{i=1}^n \mathrm{E}\left[\left(K\left(\frac{R_i}{n+1}\right) - K\left(\tilde{G}_{1k}\left(\frac{d_i}{\sigma}\right)\right)\right)^2\right],$$

where, denoting by $U_{i,j}$ the $j$th component of $\mathbf{U}_i$, $C_{\ell,k} = \mathrm{Var}[U_{1,1}^2] = 2(k-1)/(k^2(k+2))$ for $\ell \in \mathcal{L}_k := \{mk+m+1, m = 0, 1, \ldots, k-1\}$ and $C_{\ell,k} = \mathrm{Var}[U_{1,1}U_{1,2}] = 1/k^2$ for $\ell \notin \mathcal{L}_k$. Hájek's classical projection result for linear signed rank statistics ([15]; see also [44], Chapter 3) thus yields the desired result. □

PROOF OF PROPOSITION 4.1. From Lemma 4.1, we easily obtain [for any fixed value $\boldsymbol{\vartheta}_0' := (\boldsymbol{\theta}', \sigma^2, (\overset{\circ}{\mathrm{vech}}\,\mathbf{V}_0)')$ of the parameter]

$$\underset{\sim}{Q}{}^{(n)}_K = (\boldsymbol{\Delta}^{\star(n)}_{K;g_1}(\boldsymbol{\vartheta}_0))'(\mathrm{E}[K^2(U)]\mathbf{\Upsilon}_k^{-1}(\mathbf{V}_0))^{-1}\boldsymbol{\Delta}^{\star(n)}_{K;g_1}(\boldsymbol{\vartheta}_0) + o_{\mathrm{P}}(1)$$

as $n \to \infty$, under $\bigcup_{\sigma^2}\bigcup_{g_1}\{\mathrm{P}^{(n)}_{\boldsymbol{\theta},\sigma^2,\mathbf{V}_0;g_1}\}$. Part (i) of Proposition 4.1 follows, since

$$\boldsymbol{\Delta}^{\star(n)}_{K;g_1}(\boldsymbol{\vartheta}_0) \xrightarrow{\mathcal{L}} \mathcal{N}(\mathcal{J}_k(K;g_1)\mathbf{\Upsilon}_k^{-1}(\mathbf{V}_0)(\overset{\circ}{\mathrm{vech}}\,\mathbf{v}), \mathrm{E}[K^2(U)]\mathbf{\Upsilon}_k^{-1}(\mathbf{V}_0))$$

as $n \to \infty$, under $\bigcup_{\sigma^2} \bigcup_{g_1} \{\mathrm{P}^{(n)}_{\boldsymbol{\theta},\sigma^2,\mathbf{V}_0;g_1}\}$, with $\boldsymbol{\tau}' := (\mathbf{t}', s, (\mathrm{ve}\overset{\circ}{\mathrm{c}}\mathrm{h}\, \mathbf{v})')$. Again, part (ii) follows—as in the proof of Proposition 3.1—by noting that the asymptotic variance of $\boldsymbol{\Delta}^{\star(n)}_{K_{f_1};f_1}(\boldsymbol{\vartheta}_0) = \boldsymbol{\Delta}^{\star(n)}_{f_1}(\boldsymbol{\vartheta}_0)$ under $\bigcup_{\sigma^2} \{\mathrm{P}^{(n)}_{\boldsymbol{\theta},\sigma^2,\mathbf{V}_0;f_1}\}$ is $\mathcal{J}_k(f_1)\boldsymbol{\Upsilon}_k^{-1}(\mathbf{V}) = \boldsymbol{\Gamma}^{\star}_{f_1}(\boldsymbol{\vartheta}_0)$. $\square$

**A.4. Proof of Proposition 5.1.** (i) Letting

$$\underset{\sim}{\mathbf{T}}{}^{(n)}_{K} := n^{-1/2} \sum_{i=1}^{n} K\left(\frac{R_i}{n+1}\right) \operatorname{vec}\left(\mathbf{U}_i \mathbf{U}_i' - \frac{1}{k}\mathbf{I}_k\right),$$

the necessary and sufficient consistency condition (5.2) holds iff $\mathrm{E}[\mathrm{P}[\| \underset{\sim}{\mathbf{T}}{}^{(n)}_{K} \| > t | \mathbf{U}^{(n)}]] \to 1$ under $\mathcal{K}(\underline{f})$ as $n \to \infty$, for any $t \in \mathbb{R}$. Since $\mathrm{P}[\| \underset{\sim}{\mathbf{T}}{}^{(n)}_{K} \| > t | \mathbf{U}^{(n)}]$ is a strictly bounded random variable, this is equivalent to

$$\mathrm{P}[\| \underset{\sim}{\mathbf{T}}{}^{(n)}_{K} \| > t | \mathbf{U}^{(n)}] = 1 + o_{\mathrm{P}}(1), \qquad \text{under } \mathcal{K}(\underline{f}), \text{ as } n \to \infty. \tag{A.1}$$

Now, conditional on $\mathbf{U}^{(n)}$, each component $\underset{\sim}{T}{}^{(n)}_{K,\ell}$ of $\underset{\sim}{\mathbf{T}}{}^{(n)}_{K}$ is a linear rank statistic with approximate scores $K(\frac{i}{n+1})$. Under the assumptions made, the Hájek variance inequality (Theorem 3.1 in [14]) applies (conditional on $\mathbf{U}^{(n)}$), yielding, for all $\ell$ [with appropriate $r$ and $s$, $m^{(n)}_K := \frac{1}{n}\sum_{i=1}^n K(\frac{i}{n+1})$ and $\sigma^2_K := \int_0^1 K^2(u)\, du - (\int_0^1 K(u)\, du)^2$],

$$\begin{aligned} \operatorname{Var}(\underset{\sim}{T}{}^{(n)}_{K,\ell} | \mathbf{U}^{(n)}) &\le 21 \max_{1 \le i \le n} \left(U_{i,r}U_{i,s} - \frac{1}{k}\delta_{rs}\right)^2 \frac{1}{n} \sum_{i=1}^{n} \left(K\left(\frac{i}{n+1}\right) - m^{(n)}_K\right)^2 \\ &< 21 \sigma^2_K, \end{aligned} \tag{A.2}$$

since $\max_{1 \le i \le n} |U_{i,r}U_{i,s} - \frac{1}{k}\delta_{rs}| < 1$ and

$$n^{-1} \sum_{i=1}^{n} \left(K\left(\frac{i}{n+1}\right) - m^{(n)}_K\right)^2 = n^{-1} \sum_{i=1}^{n} K^2\left(\frac{i}{n+1}\right) - (m^{(n)}_K)^2 \to \sigma^2_K < \infty.$$

The bound (A.2) on the conditional variance being uniform, it follows that $\underset{\sim}{\mathbf{T}}{}^{(n)}_{K} = \underset{\sim}{\boldsymbol{\mu}}{}^{(n)}_{\mathbf{T}}(\mathbf{U}^{(n)}) + O_{\mathrm{P}}(1)$, with

$$\begin{aligned} \underset{\sim}{\boldsymbol{\mu}}{}^{(n)}_{\mathbf{T}}(\mathbf{U}^{(n)}) &:= \mathrm{E}[\underset{\sim}{\mathbf{T}}{}^{(n)}_{K} | \mathbf{U}^{(n)}] \\ &= n^{-1/2} \sum_{i=1}^{n} \mathrm{E}\left[K\left(\frac{R^{(n)}_i}{n+1}\right) \Big| \mathbf{U}^{(n)}\right] \operatorname{vec}\left(\mathbf{U}_i \mathbf{U}_i' - \frac{1}{k}\mathbf{I}_k\right). \end{aligned}$$

Consequently, the necessary and sufficient condition (A.1) takes the form $\underset{\sim}{\boldsymbol{\mu}}_{\mathbf{T}}^{(n)}(\mathbf{U}^{(n)}) \xrightarrow{\mathrm{P}} \infty$ [under $\mathcal{K}(\underline{f})$, as $n \to \infty$], which concludes the proof of (i).

(ii) Returning to $\underset{\sim}{\boldsymbol{\mu}}_{\mathbf{T}}^{(n)}(\mathbf{U}^{(n)})$ and denoting by $F^{d_i|\mathbf{U}_i}$ the distribution function, under $\mathcal{K}^{(n)}(\underline{f})$, of $d_i = d_i(\boldsymbol{\theta}, \mathbf{I}_{\boldsymbol{k}})$ conditional on $\mathbf{U}_i = \mathbf{U}_i(\boldsymbol{\theta}, \mathbf{I}_{\boldsymbol{k}})$, we have

$$\begin{aligned}
\underset{\sim}{\boldsymbol{\mu}}_{\mathbf{T}}^{(n)}(\mathbf{U}^{(n)}) &= n^{-1/2} \sum_{i=1}^{n} \mathrm{E}\left[K\left(\frac{R_i}{n+1}\right)\Big|\mathbf{U}^{(n)}\right] \operatorname{vec}\left(\mathbf{U}_i \mathbf{U}_i' - \frac{1}{k}\mathbf{I}_k\right) \\
&= n^{-1/2} \sum_{i=1}^{n} \left(\mathrm{E}\left[K\left(\frac{R_i}{n+1}\right)\Big|\mathbf{U}^{(n)}\right]\right. \\
&\qquad \left. - \int_0^\infty K\left(n^{-1} \sum_{j=1}^{n} F^{d_j|\mathbf{U}_j}(r)\right) dF^{d_i|\mathbf{U}_i}(r)\right) \\
&\qquad \times \operatorname{vec}\left(\mathbf{U}_i \mathbf{U}_i' - \frac{1}{k}\mathbf{I}_k\right) \\
&\quad + n^{-1/2} \sum_{i=1}^{n} \int_0^\infty K\left(n^{-1} \sum_{j=1}^{n} F^{d_j|\mathbf{U}_j}(r)\right) dF^{d_i|\mathbf{U}_i}(r) \\
&\qquad \times \operatorname{vec}\left(\mathbf{U}_i \mathbf{U}_i' - \frac{1}{k}\mathbf{I}_k\right) \\
&=: \mathbf{E}_1^{(n)} + \mathbf{E}_2^{(n)}, \qquad \text{say.}
\end{aligned}$$

Clearly,

$$\mathbf{E}_2^{(n)} = n^{-1/2} \sum_{i=1}^{n} \mathrm{E}\left[K\left(\frac{1}{n} \sum_{j=1}^{n} \mathrm{P}[d_j \le d_i | d_i, \mathbf{U}_i, \mathbf{U}_j]\right)\Big|\mathbf{U}^{(n)}\right] \operatorname{vec}\left(\mathbf{U}_i \mathbf{U}_i' - \frac{1}{k}\mathbf{I}_k\right).$$

As for $\mathbf{E}_1^{(n)}$, Proposition 2 in [22] implies that for each component $\mathbf{E}_{1;\ell}^{(n)}$ of $\mathbf{E}_1^{(n)}$ and appropriate $r$ and $s$,

$$\begin{aligned}
|\mathbf{E}_{1;\ell}^{(n)}| &\le n^{-1/2} \sum_{i=1}^{n} \left| \mathrm{E}\left[K\left(\frac{R_i}{n+1}\right)\Big|\mathbf{U}^{(n)}\right] \right. \\
&\qquad \left. - \int_0^\infty K\left(n^{-1} \sum_{j=1}^{n} F^{d_j|\mathbf{U}_j}(r)\right) dF^{d_i|\mathbf{U}_i}(r) \right| \left| U_{i,r} U_{i,s} - \frac{1}{k}\delta_{rs} \right| \\
&\le n^{-1/2} \sum_{i=1}^{n} \left| \mathrm{E}\left[K\left(\frac{R_i}{n+1}\right)\Big|\mathbf{U}^{(n)}\right] - \int_0^\infty K\left(n^{-1} \sum_{j=1}^{n} F^{d_j|\mathbf{U}_j}(r)\right) dF^{d_i|\mathbf{U}_i}(r) \right| \\
&\le n^{-1/2} \mathrm{E}[C n^{1/2} J(K)] = C J(K),
\end{aligned}$$

with $J(K) < \infty$ defined in (5.4) and $C \le 8$ (cf. page 359 of [44]). Hence, $\underset{\sim}{\boldsymbol{\mu}}^{(n)}_{\mathbf{T}}(\mathbf{U}^{(n)}) = \mathbf{E}_2^{(n)} + O_{\mathrm{P}}(1)$ and $\underset{\sim}{\boldsymbol{\mu}}^{(n)}_{\mathbf{T}}(\mathbf{U}^{(n)}) \xrightarrow{\mathrm{P}} \infty$ iff $\mathbf{E}_2^{(n)} \xrightarrow{\mathrm{P}} \infty$; (5.5) follows.

(iii) For each $\ell$, convexity of $K$ implies

$$\begin{aligned}
&\mathrm{E}\left[K\left(\frac{1}{n}\sum_{j=1}^{n}\mathrm{P}[d_j \le d_i | d_i, \mathbf{U}_i, \mathbf{U}_j]\right)\Big|\mathbf{U}^{(n)}\right]\\
&\qquad \le \frac{1}{n}\sum_{j=1}^{n}\mathrm{E}[K(\mathrm{P}[d_j \le d_i | d_i, \mathbf{U}_i, \mathbf{U}_j])\mathbf{U}_i, \mathbf{U}_j].
\end{aligned}$$

Similarly, Jensen's inequality implies that

$$\begin{aligned}
&n^{-1/2}\sum_{i=1}^{n}\mathrm{E}\left[K\left(\frac{1}{n}\sum_{j=1}^{n}\mathrm{P}[d_j \le d_i | d_i, \mathbf{U}_i, \mathbf{U}_j]\right)\Big|\mathbf{U}^{(n)}\right]\operatorname{vec}\left(\mathbf{U}_i\mathbf{U}_i' - \frac{1}{k}\mathbf{I}_k\right)_\ell^-\\
&\quad \ge \frac{1}{n}\sum_{m=1}^{n}\operatorname{vec}\left(\mathbf{U}_m\mathbf{U}_m' - \frac{1}{k}\mathbf{I}_k\right)_\ell^-\\
&\qquad \times n^{1/2}K\left(\left(\frac{1}{n}\sum_{m=1}^{n}\operatorname{vec}\left(\mathbf{U}_m\mathbf{U}_m' - \frac{1}{k}\mathbf{I}_k\right)_\ell^-\right)^{-1}\right.\\
&\qquad\qquad \times \frac{1}{n^2}\sum_{i=1}^{n}\sum_{j=1}^{n}\mathrm{E}[\mathrm{P}[d_j \le d_i | d_i, \mathbf{U}_i, \mathbf{U}_j] | \mathbf{U}_i, \mathbf{U}_j]\\
&\qquad\qquad\qquad \left. \times \operatorname{vec}\left(\mathbf{U}_i\mathbf{U}_i' - \frac{1}{k}\mathbf{I}_k\right)_\ell^-\right).
\end{aligned}$$

It follows that $\mathbf{E}_{2;\ell}^{(n)}$ is bounded from above by

$$\begin{aligned}
&n^{1/2}\left\{\frac{1}{n(n-1)}\sum\sum_{1\le i\ne j\le n}\mathrm{E}[K(\mathrm{P}[d_j \le d_i | d_i, \mathbf{U}_i, \mathbf{U}_j]) | \mathbf{U}_i, \mathbf{U}_j]\operatorname{vec}\left(\mathbf{U}_i\mathbf{U}_i' - \frac{1}{k}\mathbf{I}_k\right)_\ell^+\right.\\
&\quad - \mathrm{E}\left[\operatorname{vec}\left(\mathbf{U}_1\mathbf{U}_1' - \frac{1}{k}\mathbf{I}_k\right)_\ell^-\right]\\
&\quad\quad \times K\left(\left(\mathrm{E}\left[\operatorname{vec}\left(\mathbf{U}_1\mathbf{U}_1' - \frac{1}{k}\mathbf{I}_k\right)_\ell^-\right]\right)^{-1}\right.\\
&\quad\quad\quad \left.\left. \times \frac{1}{n(n-1)}\sum\sum_{1\le i\ne j\le n}\mathrm{E}[I[d_j \le d_i] | \mathbf{U}_i, \mathbf{U}_j]\operatorname{vec}\left(\mathbf{U}_i\mathbf{U}_i' - \frac{1}{k}\mathbf{I}_k\right)_\ell^-\right)\right\}\\
&\quad + o_{\mathrm{P}}(1),
\end{aligned}$$

where

$$\frac{1}{n(n-1)}\sum\sum_{1\le i\ne j\le n}\mathrm{E}[K(\mathrm{P}[d_j\le d_i|d_i,\mathbf{U}_i,\mathbf{U}_j])|\mathbf{U}_i,\mathbf{U}_j]\,\mathrm{vec}\biggl(\mathbf{U}_i\mathbf{U}_i'-\frac{1}{k}\mathbf{I}_k\biggr)_\ell^+$$

and

$$\frac{1}{n(n-1)}\sum\sum_{1\le i\ne j\le n}\mathrm{E}[I[d_j\le d_i]|\mathbf{U}_i,\mathbf{U}_j]\,\mathrm{vec}\biggl(\mathbf{U}_i\mathbf{U}_i'-\frac{1}{k}\mathbf{I}_k\biggr)_\ell^-$$

are U-statistics with finite-variance kernels

$$(\mathbf{u},\mathbf{v})\mapsto \mathrm{E}[K(\mathrm{P}[d_2\le d_1|d_1,\mathbf{U}_1=\mathbf{u},\mathbf{U}_2=\mathbf{v}])|\mathbf{U}_1=\mathbf{u},\mathbf{U}_2=\mathbf{v}]$$
$$\times\,\mathrm{vec}\biggl(\mathbf{u}\mathbf{u}'-\frac{1}{k}\mathbf{I}_k\biggr)_\ell^+$$

and

$$(\mathbf{u},\mathbf{v})\mapsto \mathrm{E}[I[d_2\le d_1]|\mathbf{U}_1=\mathbf{u},\mathbf{U}_2=\mathbf{v}]\,\mathrm{vec}\biggl(\mathbf{u}\mathbf{u}'-\frac{1}{k}\mathbf{I}_k\biggr)_\ell^-,$$

respectively. The continuous mapping theorem and standard asymptotic normality results for U-statistics (see, e.g., [21]) imply that

$$\mathbf{E}_{2;\ell}^{(n)}\le n^{1/2}\mathrm{E}\biggl[\mathrm{vec}\biggl(\mathbf{U}_1\mathbf{U}_1'-\frac{1}{k}\mathbf{I}_k\biggr)_\ell^-\biggr]$$
$$\times\biggl\{\frac{\mathrm{E}[K(\mathrm{P}[d_2\le d_1|d_1,\mathbf{U}_1,\mathbf{U}_2])\,\mathrm{vec}(\mathbf{U}_1\mathbf{U}_1'-(1/k)\mathbf{I}_k)_\ell^+]}{\mathrm{E}[\mathrm{vec}(\mathbf{U}_1\mathbf{U}_1'-(1/k)\mathbf{I}_k)_\ell^-]} \tag{A.3}$$
$$-K\biggl(\frac{\mathrm{E}[I[d_2\le d_1]\,\mathrm{vec}(\mathbf{U}_1\mathbf{U}_1'-(1/k)\mathbf{I}_k)_\ell^-]}{\mathrm{E}[\mathrm{vec}(\mathbf{U}_1\mathbf{U}_1'-(1/k)\mathbf{I}_k)_\ell^-]}\biggr)\biggr\}+O_\mathrm{P}(1).$$

A sufficient condition for (5.5) to hold is thus that the quantity in braces in this upper bound be strictly negative, yielding part (5.6) of the claim. Similar arguments imply that

$$\mathbf{E}_{2;\ell}^{(n)}\ge n^{1/2}\mathrm{E}\biggl[\mathrm{vec}\biggl(\mathbf{U}_1\mathbf{U}_1'-\frac{1}{k}\mathbf{I}_k\biggr)_\ell^+\biggr]$$
$$\times\biggl\{K\biggl(\frac{\mathrm{E}[I[d_2\le d_1]\,\mathrm{vec}(\mathbf{U}_1\mathbf{U}_1'-(1/k)\mathbf{I}_k)_\ell^+]}{\mathrm{E}[\mathrm{vec}(\mathbf{U}_1\mathbf{U}_1'-(1/k)\mathbf{I}_k)_\ell^+]}\biggr) \tag{A.4}$$
$$-\frac{\mathrm{E}[K(\mathrm{P}[d_2\le d_1|d_1,\mathbf{U}_1,\mathbf{U}_2])\,\mathrm{vec}(\mathbf{U}_1\mathbf{U}_1'-(1/k)\mathbf{I}_k)_\ell^-]}{\mathrm{E}[\mathrm{vec}(\mathbf{U}_1\mathbf{U}_1'-(1/k)\mathbf{I}_k)_\ell^+]}\biggr\}+O_\mathrm{P}(1),$$

yielding part (5.7) of the claim.

(iv) For Wilcoxon scores, the upper bound (A.3) and the lower bound (A.4) both reduce to

$$
\begin{aligned}
n^{1/2}\Big\{ & \mathrm{E}\Big[I[d_2 \le d_1]\,\mathrm{vec}\Big(\mathbf{U}_1\mathbf{U}_1' - \frac{1}{k}\mathbf{I}_k\Big)_\ell^+\Big] \\
& - \mathrm{E}\Big[I[d_2 \le d_1]\,\mathrm{vec}\Big(\mathbf{U}_1\mathbf{U}_1' - \frac{1}{k}\mathbf{I}_k\Big)_\ell^+\Big]\Big\} + O_{\mathrm{P}}(1) \\
& = n^{1/2}\mathrm{E}\Big[I[d_2 \le d_1]\,\mathrm{vec}\Big(\mathbf{U}_1\mathbf{U}_1' - \frac{1}{k}\mathbf{I}_k\Big)_\ell\Big] + O_{\mathrm{P}}(1).
\end{aligned}
$$

Part (iv) of the proposition follows.

(v) For the sign test, that is, when the score function $K$ reduces to a constant, the necessary and sufficient condition (5.3) takes the form

$$
n^{-1/2}\sum_{i=1}^{n}\mathrm{vec}\Big(\mathbf{U}_i\mathbf{U}_i' - \frac{1}{k}\mathbf{I}_k\Big) \xrightarrow{\mathrm{P}} \infty, \qquad \text{under } \mathcal{K}^{(n)}(\underline{f}), \text{ as } n \to \infty. \tag{A.5}
$$

The central limit theorem implies that this happens iff the summands in (A.5) are incorrectly centered, that is, whenever (5.9) holds. This completes the proof of Proposition 5.1. □

**Acknowledgments.** The authors are grateful to two anonymous referees and an Associate Editor for their critical comments and patient persistence which helped in the extending of the results that had been obtained in the initial version of this paper. They also wish to thank Hannu Oja for stimulating discussion and Alexandre Dutrifoy for insightful remarks on quadratic mean differentiability.

E.C.A.R.E.S., Institute for Research in Statistics
and Département de Mathématique
Université Libre de Bruxelles
Campus de la Plaine CP 210
B-1050 Bruxelles
Belgium
E-mail: mhallin@ulb.ac.be
dpaindav@ulb.ac.be
URL: http://homepages.ulb.ac.be/~dpaindav